\documentclass[english]{article}
\usepackage{mathptmx}
\usepackage{helvet}
\usepackage{courier}
\usepackage[T1]{fontenc}
\usepackage[utf8]{inputenc}
\usepackage{color}
\usepackage{babel}
\usepackage{array}
\usepackage{longtable}
\usepackage{multirow}
\usepackage{dsfont}
\usepackage{amsmath}
\usepackage{amssymb}
\usepackage{mathdots}
\usepackage{stmaryrd}
\usepackage{graphicx}
\usepackage{wasysym}
\usepackage[unicode=true,pdfusetitle,
 bookmarks=true,bookmarksnumbered=false,bookmarksopen=false,
 breaklinks=true,pdfborder={0 0 0},pdfborderstyle={},backref=false,colorlinks=true]
 {hyperref}

\makeatletter

\providecommand{\tabularnewline}{\\}
\newcommand{\lyxdot}{.}

\newcommand{\lyxaddress}[1]{
	\par {\raggedright #1
	\vspace{1.4em}
	\noindent\par}
}

\usepackage{url} 

\makeatother

\begin{document}
\title{On the Beta Transformation}
\author{Linas Vepštas}
\date{December 2017 (Current version: January 2024)\thanks{This text is an extract of a diary of research notes, first published
in December 2017, with major updates published in February 2018, December
2018 and January 2024. The current version expunges all the unfinished,
boring and speculative parts of the diary, and prettifies the rest
into a presentable format. MSC2020 Classification 37E05, 37G10.}}
\maketitle

\lyxaddress{linasvepstas@gmail.com}
\begin{abstract}
The beta transformation is the iterated map $\beta x\,\mod1$; it
generates the base-$\beta$ expansion of a real number $x$. Every
iterated piece-wise monotonic map is topologically conjugate to the
beta transformation. For all but a countable subset of $\beta$, the
orbits of $x$ are ergodic; yet it is the finite orbits that determine
overall behavior.

This is a large text; it splits into four parts. The first part provides
a review of general concepts and properties associated with the beta
shift. The second part examines the spectrum of the Ruelle–Frobenius–Perron
operator, and gives explicit expressions for a set of bounded eigenfunctions.
These form a discrete spectrum, accumulating on a circle of radius
$1/\beta$ in the complex plane.

The third part examines the finite and the periodic orbits. These
are in one-to-one correspondence with monic integer polynomials. They
are ``quasi-cyclotomic'' and can be counted with Moreau's necklace-counting
function; curiously, they do not have any obvious relation to other
systems countable by the necklace function. The positive real roots
are dense in the reals; they include the Golden and silver ratios,
the Pisot numbers, the n-bonacci (tribonacci, tetranacci, \emph{etc.})
numbers. The beta-polynomials yoke all of these together into a regular
structure. An explicit bijection to the rationals is presented.

The fourth part of this text examines small perturbations. These introduce
Arnold tongues, which inflate the finite orbits, a set of measure
zero, to finite size.

This text assumes very little mathematical sophistication on the part
of the reader, and should be approachable for any enthusiast with
minimal or no prior experience in ergodic theory. Most of the development
is casual. As a side effect, the introductory sections are perhaps
a fair bit longer than strictly needed to present the new results.

\pagebreak{}
\end{abstract}

\section{Introduction}

Given a real number $x$, it can be written down in base-ten, decimal
notation, as a string of digits 0-9. Each successive digit is given
by an iterated formula. Starting with $y_{0}=x$, let

\[
y_{n+1}=\left(10y_{n}\right)\mod1
\]
The $n$'th digit in the decimal expansion is $d_{n}=\left\lfloor y_{n}\right\rfloor $.
Here, mod is the modulo (or remainder) function, so that $a\mod b$
is the remainder of $a$ after dividing by $b$. The bracket $\left\lfloor z\right\rfloor $
is the floor function, returning the largest integer less than $z$.
The decimal expansion of $x$ is then $d_{0}.d_{1}d_{2}d_{3}\cdots$
with all digits $d_{n}$ (except for the first) running between 0
and 9. The iteration has the form of a shift; it shifts over by one
digit, and repeats the process. The shift is ergodic over the reals:
any and every decimal expansion will occur.

The same idea can be applied by replacing ten by any number $\beta$;
this gives the base-$\beta$ expansion. The $\beta$ does not need
to be an integer; the rest of this text is focused on exploring the
expansions coming from $1<\beta\le2$, so ``base two or less''.

The beta transformation is the map 
\[
t_{\beta}\left(x\right)=\beta x\mod1
\]
iterated on the unit interval $\left[0,1\right]$ of the reals, and
$\beta$ taken to be a real positive constant. For the special case
of $\beta=2$, this is the Bernoulli process, the ergodic process
of random coin flips. The Bernoulli shift is solvable in more ways
than one; it is well-understood, and provides a canonical model of
ergodic behavior. The beta transformation is less well-known and less
carefully scrutinized. A full review of previous results is presented
at the end of this introduction. The goal of this text is two-fold:
to expand the collection of known properties, and to present them
in the simplest, most accessible fashion possible.

One primary driver of interest in iterated functions is that they
generate fractals, and that they provide simple models of chaotic
dynamics. The $\beta$-transform is the simplest possible map that
exhibits all of the hallmarks of chaotic dynamical systems: a non-trivial,
nonuniform invariant measure, with both periodic and non-periodic
orbits. It is also ``generic'', in that any other piecewise-monotonic
map of the unit interval is topologically conjugate to it. If one
can understand the $\beta$-transform, one has gone a long ways towards
understanding iterated maps in general.

There are two broad approaches for studying iterated functions. One
is to examine the point dynamics and orbits: ``where does the point
$x$ go, when iterated?''. The other is in terms of distributions:
``how does a scattered dust of points evolve over time?''. Within
the context of physics, these give two broad philosophical views of
reality. The first is of microscopic, time-reversible systems whose
future is deterministic and known with infinite precision. The second
is of macroscopic, time-irreversible thermodynamics, where time can
only go forward, and the future is unknown and unknowable. Of these
two approaches, the first is commonplace and inescapable; the second
remains obscure, poorly-recognized and opaque. Thus, a large part
of this text is devoted to this second approach.

If iterating a map $x\mapsto f\left(x\right)\mapsto f\left(f\left(x\right)\right)\mapsto\cdots$
pulls a point $x$ through $f$, through time, then the action of
the map $f$ on a distribution $\rho$ is a pushforward: 
\[
\rho\left(A\right)\mapsto\rho\left(f^{-1}\left(A\right)\right)\mapsto\rho\left(f^{-1}\left(f^{-1}\left(A\right)\right)\right)\mapsto\cdots
\]
The proper definition of a pushforward requires a significant development
of the concepts of measurable spaces and Borel sigma algebras, topics
that will be gently reviewed a bit further in this introductory section.
For the present, it is enough to take $\rho:\left[0,1\right]\to\mathbb{R}$
to be some function defined on the unit interval. In the above, $A\subset\left[0,1\right]$
is a subset of the unit interval, so that $\rho\left(A\right)=\int_{A}\rho\left(x\right)dx$
is an ordinary integral, the ``size'' of the set $A$ with respect
to the distribution.

The challenge is to find an explicit expression for the pushforward
$\rho\left(f^{-1}\left(A\right)\right)$. This can be obtained as
a change of variable $y=f\left(x\right)$ under integration. Start
with any function $h$; it's integral over the set $A$ is as above:
$h\left(A\right)=\int_{A}h\left(y\right)dy$. Under the change of
variable, this becomes

\[
h\left(A\right)=\int_{A}h\left(y\right)dy=\int_{f^{-1}\left(A\right)}h\left(f\left(x\right)\right)\left|f^{\prime}\left(x\right)\right|dx
\]
Writing the integrand as $\rho\left(x\right)=h\left(f\left(x\right)\right)\left|f^{\prime}\left(x\right)\right|$
and working backwards, one recovers
\[
h\left(y\right)=\frac{\rho\left(f^{-1}\left(y\right)\right)}{\left|f^{\prime}\left(f^{-1}\left(y\right)\right)\right|}
\]
Plugging this back through gives the identity
\[
\int_{A}\frac{\rho\left(f^{-1}\left(y\right)\right)}{\left|f^{\prime}\left(f^{-1}\left(y\right)\right)\right|}dy=\int_{f^{-1}\left(A\right)}\rho\left(x\right)dx=\rho\left(f^{-1}\left(A\right)\right)
\]
The right-hand-side is the desired pushforward; the left-hand side
is an explicit expression for it. There was a minor sleight-of-hand
in the above derivation: the map $y=f\left(x\right)$ may not be one-to-one.
Thus, there may be several distinct points $x=f^{-1}\left(y\right)$.
In this case, the above needs to be amended as 

\[
h\left(y\right)=\sum_{x\in f^{-1}\left(y\right)}\frac{\rho\left(x\right)}{\left|f^{\prime}\left(x\right)\right|}
\]
As $h$ depends only on $\rho$ and $f$, the sum construction on
the right-hand side can be thought of as an operation $\mathcal{L}$,
defined by $f$, acting on $\rho$; in short-hand, $h=\mathcal{L}_{f}\rho$.

The symbol $\mathcal{L}$ is used to remind that this is a linear
operator: $\mathcal{L}\left(a\rho+b\sigma\right)=a\mathcal{L}\rho+b\mathcal{L}\sigma$
for any pair of real numbers $a,b$ and any functions $\rho,\sigma$.
The pushforward sequence now becomes
\[
\rho\left(A\right)\mapsto\left[\mathcal{L}_{f}\rho\right]\left(A\right)\mapsto\left[\mathcal{L}_{f}\mathcal{L}_{f}\rho\right]\left(A\right)\mapsto\cdots
\]
Thus, we've defined a linear operator $\mathcal{L}_{f}$ that depends
only on the iterated function $f$, and has the property of mapping
distributions to other distributions as it is iterated. It is the
result of commuting with function composition: $\mathcal{L}_{f}\circ\rho=\rho\circ f^{-1}$;
it's a kind of a trick with function composition. Indeed, one can
define an analogous operator $\mathcal{K}_{f}\circ\rho=\rho\circ f$,
the ``composition operator'' or ``Koopman operator'', that acts
as a kind of (one-sided) inverse to $\mathcal{L}_{f}$.

Formally, the pushforward $\mathcal{L}_{f}$ is called the ``transfer
operator'' or the ``Ruelle–\-Frobenius–\-Perron operator''. As
a linear operator, the full force of operator theory comes into play.
The primary task is to describe it's spectrum (it's eigenfunctions
and eigenvalues). Two aspects of this spectrum are interesting. The
first is the so-called ``invariant measure'', the distribution $\mu:\left[0,1\right]\to\mathbb{R}$
that defines a density on the unit interval that is invariant under
the application of the pushforward: $\mathcal{L}_{f}\mu=\mu$. An
informal example of such an invariant measure are the rings of Saturn:
an accumulation of dust and gravel, orbiting Saturn, coupled by gravitation
to both Saturn and orbiting moons, yet in a stable dynamical distribution.
This is the physical meaning and importance of the invariant measure;
more generally, it appears as the ``ground state'' or ``thermodynamic
equilibrium state'' in a vast variety of dynamical systems.

Aside from the invariant measure, there is also the question of the
rest of the spectrum. These are described by the eigenfunctions $\rho$
satisfying $\mathcal{L}\rho=\lambda\rho$. By the theorem of Frobenius–Perron,
all these other solutions have eigenvalue $\left|\lambda\right|<1$.
In physics, these correspond to the decaying modes, to the distributions
that disappear over time. For the example of Saturn, these are anything
not orbiting in the plane of the rings: tidal forces and perturbations
from the moons will force such orbits either into the ring, or crash
into the planet, or possibly fly away to infinity. The other orbits
are not stable. Thus, a characterization of the decaying spectrum
is of general interest.

A much stronger conception is that the decaying spectrum has something
to do with the irreversibility of time. In the macroscopic world,
this is plainly obvious. In the microscopic world, the laws of physics
are manifestly time-reversible. Somehow, complex dynamical systems
pass through a region of chaotic and turbulent motion, culminating
in thermodynamic equilibrium. The decaying spectrum provides a conceptual
framework in which one can ponder this transition.

This is where the mathematical fun begins. The spectrum is not a ``fixed
thing'', but depends strongly on the space of functions on which
$\mathcal{L}_{f}$ is allowed to act. If one limits oneself to $\rho$
drawn from the space of piece-wise continuous and smooth functions,
\emph{i.e.} polynomials, then $\mathcal{L}_{f}$ will in general have
a discrete spectrum. If instead, $\rho\in L^{2}\left[0,1\right]$
the space of functions that are square-integrable on the unit interval,
then the spectrum will often be continuous, and perhaps may have a
large kernel. Larger spaces exhibit even wilder behavior: if one asks
only that $\rho$ be $L^{1}$-integrable (not square-integrable),
then it is possible for continuous-nowhere functions to appear as
eigenfunctions of $\mathcal{L}_{f}$. An explicit example of the latter
is the Minkowski measure for the transfer operator of the Minkowski
Question Mark function: it vanishes on the rationals, but can be integrated
just fine; it's integral is the Question Mark function. In short,
a rich variety can often be found. In the present case, it seems,
nothing quite this rich, but getting there.

\subsection{Summary}

This text is quite long and large, as it is a summary of the results
from a multi-year investigation into the $\beta$-transform. It ranges
over a variety of disparate topics.

\textbf{\textit{Part one}} continues with a general overview of basic
concepts, including the Bernoulli shift, the definition of shift spaces
in a relatively abstract fashion, some basics about the $\beta$-shift,
and a collection of pretty visualizations to anchor and motivate.
Part one concludes with a short review of prior research into the
$\beta$-shift.

\textbf{\textit{Part two}} defines and works with the transfer operator
$\mathcal{L}_{\beta}$ for the $\beta$-shift. The concept of ``analytic
algorithmics'' is introduced, by analogy to the idea of ``analytic
combinatorics''. The Heaviside theta function $\Theta\left(x\right)$
can be viewed as the algorithmic snippet ``if $x>0$ then 1 else
0''; eigenfunctions of the transfer operator are series summations
over $\Theta\left(\beta t_{n}-1\right)$ for a sequence of iterates
$t_{n}$. A countable number of eigenvalues $\lambda$ of $\mathcal{L}_{\beta}$
are shown to accumulate on the circle $\left|\lambda\right|=1/\beta$
in the complex plane.

The description is incomplete. Examples of eigenvalues with $\left|\lambda\right|<1/\beta$
can be found, but a comprehensive description is elusive. Iterating
$\mathcal{L}_{\beta}$ shows that there are almost-resonances: distributions
that are almost, but not quite stable under repeated application of
$\mathcal{L}_{\beta}$. These remain uncharacterized.

\textbf{\textit{Part three}} takes a close look at those $\beta$
values which are characterized by finite orbits. These $\beta$ values
have the curious property of being ``self-describing'': they occur
as roots of a polynomial generated by the orbit of $x=1$ under the
iteration of $t_{\beta}\left(x\right)=\beta x\mod1$. The orbit produces
a bit-string $b_{0}b_{1}\cdots b_{k-1}$ of binary bits. When these
are arranged into a polynomial $p_{n}\left(z\right)=z^{k+1}-b_{0}z^{k}-b_{1}z^{k-1}-\cdots-b_{k-1}z-1$,
the unique real positive root of $p_{n}\left(z\right)$ is then the
$\beta$ that generated the orbit: thus, ``self-describing''.

The polynomial $p_{n}\left(z\right)$ is an odd beast. It is ``quasi-cyclotomic'',
in that the (complex) roots are approximately evenly distributed in
an approximately circular fashion on the complex plane. It is a generator
of number sequences that generalize the Fibonacci numbers; indeed,
the first polynomial in the series, $p_{1}\left(z\right)=z^{2}-z-1$
has the Golden mean $\varphi=\left(1+\sqrt{5}\right)/2$ as the (real,
positive) root. The assortment of popular ``generalized Fibonacci
numbers'' are all special cases of this class. There seems to be
a rich associated number theory to go with this, but it does not seem
to align with any known ``conventional'' combinatorics. For example,
the $p_{n}\left(z\right)$ vaguely resemble polynomials over the field
$\mathbb{F}_{2}$, but that is as far as things seem to go: a vague
resemblance. It is tempting to guess that ideas, concepts and theorems
from Galois theory might go through, but they don't: the polynomials
don't fit that mold in any obvious way. It seems to be a large class,
but there is no existing theory that it maps to.

Such self-describing $\beta$ values are dense in the reals, and can
be placed in one-to-one correspondence with the dyadic rationals.
This is not as easy as it sounds: not all possible bit-sequences occur,
and thus much of the work is to describe the ``valid'' bit-sequences.
An explicit formula is presented: it provides a bijection between
the set of valid bit-sequences and the set of all bit-sequences. This
bijection is articulated; since the bit-sequences are countable and
dense, and the reals are seperable, the mapping extends to a monotonically
increasing continuous, nowhere-differentiable function on the reals.

The non-terminating orbits fall into two classes: the ultimately-periodic
orbits, and the ergodic orbits. The ultimately-periodic orbits correspond
to rationals, but again, not all rationals are ``valid''; they don't
form self-describing orbits. The above-mentioned bijection saves the
day: it maps all rationals to valid rationals.

\textbf{\textit{Part four}} explores perturbations of the $\beta$-map.
The exploration is driven by the idea that the periodic orbits re-materialize
as the ``islands of stability'' in other iterated systems. In the
$\beta$-map, the periodic orbits are dense in $\beta$, but they
form a set of measure zero: they are countable. Small perturbations
of the $\beta$-map give rise to Arnold tongues. These are now finite-sized
regions occupied by the periodic orbits: these are the ``mode-locked''
regions. As noted earlier, the $\beta$-map is topologically conjugate
to \emph{all} piece-wise continuous maps of the unit interval; but
this is a very abstract statement. Explicit examples of mode-locking
and Arnold tongues illustrates exactly what this conjugacy actually
looks like. Literally; this section is filled with pretty pictures.

\subsection{Bernoulli shift}

The Bernoulli shift, also known as the bit-shift map, the dyadic transform
and the full shift, is an iterated map on the unit interval, given
by 
\begin{equation}
b(x)=\begin{cases}
2x & \mbox{ for }0\le x<\frac{1}{2}\\
2x-1 & \mbox{ for }\frac{1}{2}\le x\le1
\end{cases}\label{eq:Bernoulli shift}
\end{equation}
It can be written much more compactly as $b(x)=2x\mod1$. The symbolic
dynamics of this map gives the binary digit expansion of $x$. That
is, write 
\[
b^{n}(x)=(b\circ b\circ\cdots\circ b)(x)=b(b(\cdots b(x)\cdots))
\]
to denote the $n$-fold iteration of $b$ and let $b^{0}(x)=x$. The
symbolic dynamics is given by the bit-sequence
\begin{equation}
b_{n}\left(x\right)=\begin{cases}
0 & \mbox{ if }0\le b^{n}(x)<\frac{1}{2}\\
1 & \mbox{ if }\frac{1}{2}\le b^{n}(x)\le1
\end{cases}\label{eq:bernoulli-bits}
\end{equation}
Attention: $n$ is a subscript on the left, and a superscript on the
right! The left is a sequence, the right is an iteration. Using the
letter $b$ one both sides is a convenient abuse of notation. Notation
will be abused a lot in this text, except when it isn't. The symbolic
dynamics recreates the initial real number:
\begin{equation}
x=\sum_{n=0}^{\infty}b_{n}\left(x\right)2^{-n-1}\label{eq:homomorphism}
\end{equation}
All of this is just a fancy way of saying that a real number can be
written in terms of it's base-2 binary expansion. That is, the binary
digits for $x$ are the $b_{n}=b_{n}\left(x\right)$, so that
\[
x=0.b_{0}b_{1}b_{2}\cdots
\]
is a representation of a real number with a bit-string.

\subsection{Bijections}

A variety of mathematical objects that can be placed into a bijection
with collections of bit-strings, and much of this text is an exploration
of what happens when this is done. There will be several recurring
themes; these are reviewed here.

The collection of all infinitely-long bit-strings $\left\{ 0,1\right\} ^{\omega}=\left\{ 0,1\right\} \times\left\{ 0,1\right\} \times\cdots$
is known as the Cantor space; $\omega$ denotes countable infinity,
so this is a countable product of repeated copies of two things. Closely
related is the Cantor set, which is famously the collection of points
$y=\sum_{n=0}^{\infty}b_{n}\left(x\right)3^{-n-1}$ that results from
taking the binary expansion of a real number, and re-expressing it
as a base-three expansion. The Cantor set can also be constructed
by repeatedly removing the middle-third. If one is careful that the
middle-third is always an open set, what remains after a single subtraction
is a closed set. What remains after infinite repetition is a ``perfect
set'', and a key theorem is that this perfect set is identical to
the collection of points obtained with the sum above. Bouncing between
these two distinct constructions requires the definition of the product
topology on Cantor space, and thence the Borel sigma algebra, so that
one can work in a consistent way with set complements. These ideas
will be reviewed as the need arises.

Associated with Cantor space is the infinite binary tree. Any given
location in the tree can be specified by giving a sequence of left-right
moves, down the tree, starting at the root. Such left-right moves
can (of course) be interpreted as bit-strings. After a finite number
of moves, one arrives at a node, and under that node extends another
infinite binary tree, just like the original. If one has a function
$f\left(b\right)$ that is defined on every node $b$ of the binary
tree, then one can compare this function to $f\left(Lb\right)$ and
$f\left(Rb\right)$ that result from a left move and a right move.
If these are equal to each other, or scale in some way when compared
to the original, or if there exist two other functions $g_{L}$ and
$g_{R}$ such that $f\left(Lb\right)=g_{L}\left(f\left(b\right)\right)$
for all $b$, and likewise, that $f\left(Rb\right)=g_{R}\left(f\left(b\right)\right)$,
then one has fractal self-similarity. More explicitly, whenever one
has a pair of commuting diagrams, $f\circ L=g_{L}\circ f$ and , $f\circ R=g_{R}\circ f$,
then one has a dyadic monoid self-symmetry. This is the symmetry of
a large class of fractals.

\subsubsection{Formalities}

The last paragraph is a bit glib, and so some formalities and examples
are in order. These are all very straightforward and conventional,
almost trivial, belaboring the obvious. Despite the seeming triteness
of the next handful of paragraphs, these formal definitions will be
needed, so as to avoid future ambiguities and confusions.

Let $\mathbb{M}=\left\{ L,R\right\} ^{<\omega}$ denote the collection
of finite-length binary strings. These can be graded by the length
$\nu$ of the string, so that 
\[
\mathbb{M}=\left\{ L,R\right\} ^{<\omega}=\bigcup_{\nu=0}^{\infty}\left\{ L,R\right\} ^{\nu}=\epsilon\cup\left\{ L,R\right\} \cup\left(\left\{ L,R\right\} \times\left\{ L,R\right\} \right)\cup\cdots
\]
with $\epsilon$ denoting the empty (zero-length) string. This can
be turned into a monoid by defining multiplication as string concatenation:
given $\gamma\in\left\{ L,R\right\} ^{n}$ some sequence of $L,R$
moves of length $n$, and $\gamma^{\prime}\in\left\{ L,R\right\} ^{m}$
some other sequence of length $m$, then $\gamma\gamma^{\prime}\in\left\{ L,R\right\} ^{n+m}$
is some other string of length $m+n$.

This set is in bijection to the integers in a straight-forward way.
This bijection, written as $\kappa:\left\{ L,R\right\} ^{<\omega}\to\mathbb{N}$,
can be defined recursively, by a simple commuting diagram. Define
$\kappa\left(\epsilon\right)=1$ and ask that $\kappa\left(L\gamma\right)=2\kappa\left(\gamma\right)$
and that $\kappa\left(R\gamma\right)=2\kappa\left(\gamma\right)+1$
for every $\gamma\in\left\{ L,R\right\} ^{n}$. Thus, $\kappa\left(L\right)=2$
and $\kappa\left(R\right)=3$ and $LL,LR,RL,RR$ map to $4,5,6,7$
respectively.

This bijection commutes with the canonical $L,R$ moves on the natural
numbers. Write these as a pair of functions $L:\mathbb{N}\to\mathbb{N}$
and $R:\mathbb{N}\to\mathbb{N}$, defined as $L:m\mapsto2m$ and $R:m\mapsto2m+1$.
Interpreting these as strings, one promptly has a pair of commuting
diagrams $\kappa\circ L=L\circ\kappa$ and $\kappa\circ R=R\circ\kappa$.
True formality would have required writing $\kappa\circ L_{s}=L_{\mathbb{N}}\circ\kappa$
or even that $\kappa\left(L_{s}\gamma\right)=L_{\mathbb{N}}\left(\kappa\left(\gamma\right)\right)$
in order to remind us that $L_{s}$ is a string (of length one) concatenated
onto some other string $\gamma$, while $L_{\mathbb{N}}$ is a map
of the natural numbers. It is convenient to drop these labels, as
they mostly serve to clutter the text. The intended meaning is always
clear from context.

Associated with the set $\left\{ L,R\right\} ^{<\omega}$ is a binary
tree $\mathbb{B}$. It can be defined as a graph of vertexes $v_{j}$
and edges $e_{ij}$ connecting vertex $v_{i}$ to vertex $v_{j}$.
Formally, it is the graph $\mathbb{B}=\left\{ v_{i},e_{ij}:i\in\mathbb{N},j\in\left\{ 2i,2i+1\right\} \right\} $.
Every vertex $v_{i}\in\mathbb{B}$ can be given an integer label:
it is just the integer $i$ itself. To formalize this, there is a
map $\eta:\mathbb{N}\to\mathbb{B}$ that provides this labeling. The
canonical labeling gives the root node a label of 1, the left and
right sub-nodes 2,3, and so on.

The canonical moves on this binary tree are $L:\mathbb{B}\to\mathbb{B}$
and $R:\mathbb{B}\to\mathbb{B}$ defined by $L:v_{i}\mapsto v_{2i}$
and $L:v_{i}\mapsto v_{2i+1}$. Just as above, these commute with
the left and right moves on the integers. So, $\eta\circ L=L\circ\eta$
and $\eta\circ R=R\circ\eta$.

The pattern repeats. Consider the set $\mathbb{D}$ of dyadic rationals
between zero and one. These are fractions that can be written as $\left(2n+1\right)/2^{m}$
for some non-negative integers $m,n$. These are in one-to-one correspondence
with the integers: there is a canonical bijection is $\delta:\mathbb{D}\to\mathbb{N}$
given by $\delta:\left(2n+1\right)/2^{m}\mapsto2^{m-1}+n$. There
are obvious left and right moves, given by $\delta\circ L=L\circ\delta$
and $\delta\circ R=R\circ\delta$. So, for example, $L\left(1/2\right)=1/4$
and $R\left(1/2\right)=3/4$. Viewed as a tree, this places 1/2 at
the root of the tree, and 1/4 and 3/4 as the nodes to the left and
right.

All of these maps were bijections: they are all invertable. They place
elements of all four objects $\mathbb{M},\mathbb{N},\mathbb{B},\mathbb{D}$
in one-to-one correspondence with each other. The commutation of the
$L,R$ moves guarantees that the multiplication on $\mathbb{M}$,
i.e. string concatenation, allows elements $\gamma\in\mathbb{M}$
to act on $\mathbb{N},\mathbb{B},\mathbb{D}$ in the obvious, intended
way.

\subsubsection{Example: Julia Sets}

As a practical example of the above machinery, consider the Julia
set of the Mandelbrot map. Recall the Mandelbrot map is an iterated
map on the complex plane, given by $z\mapsto z^{2}+c$. The Julia
set is the set of points of ``where things came from'' in the Mandelbrot
map; it is the inverse ``map'' $z\mapsto\pm\sqrt{z-c}$. The word
``map'' is in scare quotes, as the plus-minus in front of the square
root indicate that each $z$ maps to either one of two distinct predecessors.
The choice of the plus-minus signs can be interpreted as left-right
moves, and so the Julia set can be interpreted as a representation
of the binary tree, with it's elements labeled by integers, or dyadic
fractions, or nodes in the binary tree, or strings of ones and zeros.

We take a moment to make this explicit. Fix a point $c\in\mathbb{C}$
in the complex plane. Define a function $j_{c}:\mathbb{N}\to\mathbb{C}$
recursively, by writing $j_{c}\left(1\right)=0$ and then $j_{c}\left(2m\right)=-\sqrt{j_{c}\left(m\right)-c}$
and $j_{c}\left(2m+1\right)=+\sqrt{j_{c}\left(m\right)-c}$. Equivalently,
there are a pair of moves on the complex plane, $L_{c}:\mathbb{C}\to\mathbb{C}$
and $R_{c}:\mathbb{C}\to\mathbb{C}$ given by $L_{c}:z\to-\sqrt{z-c}$
and $R_{c}:z\to+\sqrt{z-c}$ which commute with the Julia map: $j_{c}\circ L=L_{c}\circ j_{c}$
and $j_{c}\circ R=R_{c}\circ j_{c}$.

The skeleton of the Julia set is a set of points in the complex plane:
\[
\mathbb{J}_{c}=\left\{ z\in\mathbb{C}:z=j_{c}\left(n\right),n\in\mathbb{N}\right\} 
\]
The full Julia set is the closure $\overline{\mathbb{J}}_{c}$ which
includes all of the limit points of $\mathbb{J}_{c}$. The set $\mathbb{J}_{c}$
is countable; the closure $\overline{\mathbb{J}}_{c}$ is uncountable.

\subsubsection{Example: de Rham Curves}

The above construction generalizes. Given any pair of functions $f:X\to X$
and $g:X\to X$ and a point $x_{0}\in X$ for any space $X$, the
dyadic monoid induces a map $j:\mathbb{N}\to X$ given recursively
as $j\left(1\right)=x_{0}$ and $j\left(2m\right)=f\left(j\left(m\right)\right)$
and $j\left(2m+1\right)=g\left(j\left(m\right)\right)$. As a commuting
diagram, one has that $j\circ L=f\circ j$ and $j\circ R=g\circ j$.
This in turn induces a set $\mathbb{J}=\left\{ x\in X:x=j\left(n\right),n\in\mathbb{N}\right\} $.
If the space $X$ is a topological space, so that limits can be meaningfully
taken, then one also has the closure $\overline{\mathbb{J}}$.

In 1957, Georges de Rham notes that if $X$ is a topological space,
so that continuity can be defined, and if there exist two points $x_{a},x_{b}\in X$
such that $g\left(x_{a}\right)=f\left(x_{b}\right)$, then the set
$\overline{\mathbb{J}}$ is a continuous curve.\cite{DeRham57} The
original proof also requires that $X$ be a metric space, and that
the two functions $f,g$ be suitably contracting, so that the set
$\mathbb{J}$ remains compact. More precisely, so that the Banach
fixed-point theorem can be applied, giving two fixed points $x_{a}=f\left(x_{a}\right)$
and $x_{b}=g\left(x_{b}\right)$. Without compactness, the curve wanders
off to infinity, where conceptions of continuity break down. It is
no longer a curve, ``out there''.

The continuity condition $g\left(x_{a}\right)=f\left(x_{b}\right)$
and the fixed points have a direct interpretation from the viewpoint
of the binary tree $\mathbb{B}$. Pick a point $x\in X$, any point
at all, and make nothing but left moves: the infinite string $LLL\cdots$.
The map converts this to the iteration $f\left(f\left(\cdots\left(x\right)\right)\right)$
which converges onto the fixed point $x_{a}=f\left(x_{a}\right)$.
Likewise, all right-moves converge onto $x_{b}=g\left(x_{b}\right)$.
In between, there all other branches in $\mathbb{B}$; but there are
also the ``gaps'' in between the branches.

Consider the two paths $RLLL\cdots$ and $LRRR\cdots$ down the tree.
Both start at the root, but end up at different places. Yet, they
are immediate neighbors: there are no other branches ``in between''
these two. Such immediate neighbors always lie at either end of a
``gap''. Each gap is headed up by the root that sits immediately
above them, so that each gap can be labeled by the node from which
these two distinct branches diverged. The continuity condition asks
that these gaps be closed up: the requirement that $g\left(x_{a}\right)=f\left(x_{b}\right)$
is the requirement that the two sides of the central gap converge
to the same point. The curve becomes continuous at this point. By
self-similarity, each gap in the tree closes up; the curve is continuous
at all such gaps.

As a specific example, consider $X=\mathbb{R}$ and $f\left(x\right)=x/2$
and $g\left(x\right)=\left(x+1\right)/2$. The fixed points are $f\left(0\right)=0$
and $g\left(1\right)=1$ and the continuity condition is satisfied:
$g\left(0\right)=f\left(1\right)$. Iteration produces a curve that
is just all of the real numbers of the unit interval. This curve is
just the standard mapping of the Cantor space to the unit interval:
it is one-to-one for all points that are not dyadic rationals, and
it is two-to-one at the dyadic rationals, as the continuity condition
explicitly forces the two-to-one mapping.

Note that Julia sets are not de Rham curves: they don't satisfy the
continuity criterion.

\subsubsection{Shifts}

Adjoint to the left and right moves is a shift $\tau$ that undoes
what the L and R moves do. It cancels them out, so that $\tau\circ L=\tau\circ R=e$
with $e$ the identity function. 

Given a string $\gamma\in\mathbb{M}$ of length $\nu$, consisting
of letters $a_{0}a_{1}\cdots a_{\nu-1}$ so that each $a_{k}\in\left\{ L,R\right\} $
is a single letter, define the shift $\tau:\mathbb{M}\to\mathbb{M}$
as the function that lops off a single letter from the front, so that
$\tau:a_{0}a_{1}\cdots a_{\nu-1}\mapsto a_{1}a_{2}\cdots a_{\nu-1}$
is a string that is one letter shorter. This shift is adjoint to the
moves $L,R$, which prepend either L or R to the string. That is,
$L:\mathbb{M}\to\mathbb{M}$ which acts on strings as $L:\gamma\mapsto L\gamma$,
and likewise, $R:\gamma\mapsto R\gamma$. Then, taken as functions,
$\tau\circ L=\tau\circ R=e$ with $e$ being the identity function
on $\mathbb{M}$, the function that does nothing. The shift $\tau$
is only an adjoint, not an inverse, since there is no way to reattach
what was lopped off, at least, not without knowing what it was in
the first place. Thus $L\circ\tau\ne e\ne R\circ\tau$. The maps $L,R$
were one-to-one but not onto; the map $\tau$ is onto but not one-to-one.

The shift can be composed with either $\kappa$ or with $\delta$,
to have the obvious effect. It's handy to introduce a new letter and
a new function $T:\mathbb{N}\to\mathbb{N}$ so that one has $T\circ\kappa=\kappa\circ\tau$
acting as $T\left(2m\right)=T\left(2m+1\right)=m$. Since $\tau$
applied to the empty string returns the empty string, so also $T\left(1\right)=1$.
Recycling the same letter $T:\mathbb{D}\to\mathbb{D}$ and defining
it so that $\delta\circ T=T\circ\delta$ one can infer that $Tx=2x\mod1$
for any $x\in\mathbb{D}$. So, for example, $T\left(1/4\right)=T\left(3/4\right)=1/2$.

On the binary tree $\mathbb{B}$, the shift moves back up the tree,
from either the left or the right side. That is, given a vertex $v\in\mathbb{B}$,
it is the map $\tau:v_{2m}\mapsto v_{m}$ and likewise $\tau:v_{2m+1}\mapsto v_{m}$.

For the Julia set example, it has a meaningful form: $\tau:z=j_{c}\left(n\right)\mapsto z^{2}+c$.
It re-does what the two Julia set maps undid. It is onto: it maps
$\mathbb{J}_{c}$ into all of $\mathbb{J}_{c}$, and likewise $\overline{\mathbb{J}}_{c}$
onto $\overline{\mathbb{J}}_{c}$. For the de Rham curve example,
it maps the curve back onto itself. In all three examples, these sets
are fixed points of $\tau$. Taking $\mathbb{J}_{c}\subset\mathbb{C}$
as a subset of the complex plane, it is invariant under the action
of $\tau$, so that one has $\tau\left(\mathbb{J}_{c}\right)=\mathbb{J}_{c}$
and likewise $\tau\left(\mathbb{\overline{\mathbb{J}}}_{c}\right)=\mathbb{\overline{\mathbb{J}}}_{c}$.
Likewise, the de Rham curve stays fixed in $X$. These are all examples
of invariant subspaces.

\subsubsection{Completions}

The previous section defined a binary tree $\mathbb{B}$, but this
tree is not the ``infinite binary tree'' alluded to in the opening
paragraphs. It is incomplete, in that it does not go ``all the way
down'' to it's leaves. It is not compact, in the same sense that
the dyadic fractions $\mathbb{D}$ are not compact: the limit points
are absent. The Cantor tree $\overline{\mathbb{B}}$ is the closure
or completion ore compactification of $\mathbb{B}$; it contains all
infinitely-long branches, all the way down to the ``leaves'' of
the tree. The Cantor tree $\overline{\mathbb{B}}$ is in one-to-one
correspondence with the Cantor space $\left\{ 0,1\right\} ^{\omega}$,
and both can be mapped down to the reals on the unit interval, using
eqn \ref{eq:homomorphism}. None of this is particularly deep, but
a few paragraphs to articulate these ideas will help avoid later confusion
and imprecision.

To convert letter strings to binary strings, define a function $b:\mathbb{M}\to\left\{ 0,1\right\} $
such that it returns 0 if the first letter of a string is L, and otherwise
it returns 1. If the string is of zero length, then one has a choice:
one can take $b$ to be undefined, or let it be 0, or 1, or introduce
a wild-card character $*=0\vee1$ denoting ``either zero or one''.
For the present, any of these choices is satisfactory. The wildcard
is appealing when working with the product topology; but, at the moment,
we have no topologies at play.

To extract the $n$'th letter from a string, define $b_{n}:\mathbb{M}\to\left\{ 0,1\right\} $
as $b_{n}=b\circ\tau^{n}$. Thus, given $\gamma\in\mathbb{M}$ of
length $\nu$, one can create a bitstring $b_{0}b_{1}\cdots b_{\nu-1}$.
It can be assigned the obvious numerical value:

\[
\left[\delta^{-1}\circ\kappa\right]\left(\gamma\right)=\sum_{n=0}^{\nu-1}b_{n}\left(\gamma\right)2^{-n-1}
\]
Comparing this to eqn \ref{eq:homomorphism}, the completion is obvious:
$\mathbb{M}=\left\{ L,R\right\} ^{<\omega}$ is completed as $\overline{\mathbb{M}}=\left\{ L,R\right\} ^{\omega}$
so that it contains all strings of infinite length. This is consistent
with the completion of the dyadic rationals $\mathbb{D}$ being the
entire real unit interval: $\overline{\mathbb{D}}=\left[0,1\right]$.
There is no completion $\overline{\mathbb{N}}$ of the countable numbers,
at least, not unless one wishes to say that it is the uncountable
infinity. This could be done, but then the games gets even more circular,
as this completion is just the Cantor set, and we already have that.
It seems best to leave $\overline{\mathbb{N}}$ undefined, to avoid
circular confusions. The completion $\overline{\mathbb{B}}$ engenders
similar confusion. In the original definition, $\mathbb{B}$ was defined
as a graph with a countable number of vertexes, each labeled with
an integer. This labeling must be abandoned: $\overline{\mathbb{B}}$
is a graph with an uncountable number of vertexes, each labeled by
an element from $\overline{\mathbb{M}}$.

The distinction between $\overline{\mathbb{M}}$, $\overline{\mathbb{B}}$
and $\left\{ 0,1\right\} ^{\omega}$ becomes hard to maintain at this
point: they are all isomorphic. The distinction between $\overline{\mathbb{M}}$
and $\left\{ 0,1\right\} ^{\omega}$ is particularly strained: they
are both collections of strings in two symbols. The primary purpose
of trying to maintain this distinction is to remind that $\overline{\mathbb{M}}$
should be though of as a collection of actions that can be applied
to sets, while $\left\{ 0,1\right\} ^{\omega}$ is a set, a collection
of points that sometimes act as labels for things. This distinction
is useful for avoiding off-by-one mistakes during calculations; it
is a notational convenience. This is a variant of common practice
in textbooks: after showing that two things are isomorphic, only rarely
is the notation collapsed into one big tangle. One maintains a Rosetta
Stone of different ways of writing the same thing. And so here: a
distinction without a difference.

\subsection{Shift space}

The shift $\tau$ was defined above as an operator that takes a sequence
of letters, and lops off the left-most symbol, returning a new sequence
that is the remainder of the string. A shift space $S$ is any subset
of the set of all infinite strings that remains invariant under the
shift: $\tau S=S$.

In general settings, one considers a vocabulary of $N$ letters, and
the set of infinite sequences $N^{\omega}$, so that a shift space
$S\subseteq N^{\omega}$ is a subset of the ``full shift'' $N^{\omega}$
(which is trivially invariant under $\tau$). Shifts that are proper
subsets of a full shift will be called subshifts. For the Bernoulli
shift, there are $N=2$ letters, and the Bernoulli shift is by definition
the full shift $2^{\omega}=\left\{ 0,1\right\} ^{\omega}$. A trivial
example of a subshift that is not a full shift is the set $S=\left\{ 0^{\omega},1^{\omega}\right\} $;
it has two elements, and is obviously invariant; both $0^{\omega}$
and $1^{\omega}$ are fixed points of $\tau$. Another example is
$S=\left\{ \left(01\right)^{\omega},\left(10\right)^{\omega}\right\} $,
where $\left(01\right)^{\omega}=0101\cdots$ is a repeating periodic
string. Any collection of such periodic strings forms a subshift.

Clearly, the union of two subshifts is a subshift, and so, to classify
subshifts, one wants to find all of the indecomposable pieces, and
characterize those. Factors include periodic strings of fixed period;
but not all of these are unique: so, the period-4 string $\left(0101\right)^{\omega}$
is really just the period-two string in disguise.

Subshifts consisting entirely of periodic strings can be characterized
in terms of Lyndon words. Lyndon words are fixed length strings that
are not decomposable into shorter sequences. Thus, each one characterizes
a periodic subshift. Cyclic permutations of a Lyndon word give the
same subshift; for example, both $\left(01\right)^{\omega}$ and $\left(10\right)^{\omega}$
belong to the same subshift. The number of distinct, unique subshifts
of length $\nu$ is given by Moreau's necklace counting function:
it counts the number of distinct sequences of a given length, modulo
cyclic permutations thereof.

Characterizing the subshifts that do not consist of periodic orbits
is considerably harder. For example, consider the string $s=010^{\omega}$.
It has an orbit: $\tau s=10^{\omega}$ and $\tau^{2}s=0^{\omega}$
and so one can write down a set $\left\langle s\right\rangle =\left\{ 010^{\omega},10^{\omega},0^{\omega}\right\} $
which has the property that $\tau\left\langle s\right\rangle \subset\left\langle s\right\rangle $.
However, it is not a subshift, because $\tau\left\langle s\right\rangle \ne\left\langle s\right\rangle $.
The first two points ``wander away'' under the application of the
shift; they are part of the ``wandering set''. What remains is the
fixed point $\tau0^{\omega}=0^{\omega}$. The ergodic decomposition
theorem states that all such sets $X$ having the property that $\tau X\subset X$
can be decomposed into two pieces: $X=S\cup W$ with $S$ a subshift,
$\tau S=S$ and $W$ the wandering set or dissipative set, that eventually
dissipates into the empty set: $\lim_{n\to\infty}\tau^{n}W=\varnothing$.
Subshifts are fixed-points; everything else disappears.

Given some real number $x\in\left[0,1\right]$, and it's binary expansion
$x=0.b_{0}b_{1}\cdots$, defined in eqn \ref{eq:bernoulli-bits},
what is the nature of $\left\langle s\right\rangle $ for $s=b_{0}b_{1}\cdots$
? That is, defining
\[
\left\langle s\right\rangle =\left\{ \gamma=\tau^{n}s:n\in\mathbb{N}\right\} 
\]
what portion of $\left\langle s\right\rangle $ is wandering, and
which part is a subshift? If $x=p/q$ is a rational number, the answer
is easy: rational numbers have binary expansions that are eventually
periodic. They consist of some finite-length prefix of non-repeating
digits, followed by an infinite-length cyclic orbit. The finite-length
prefix is the wandering set; the cyclic part is a subshift. If the
period of the cyclic part is $\nu$, then the subshift contains precisely
$\nu$ elements.

For the Bernoulli shift $\tau x=2x\mod1$, for the real numbers, the
answer is provided by the ergodic theorem. For all real numbers $x\in\left[0,1\right]\backslash\mathbb{Q}$,
that is, the unit interval excluding the rationals, the orbit of $x$
is ergodic: given any real number $\varepsilon>0$ and any $y\in\left[0,1\right]$
there exists some $n\in\mathbb{N}$ such that $\left|y-\tau^{n}x\right|<\varepsilon$.
Iteration takes $x$ arbitrarily close to any location on the unit
interval. In terms of symbolic dynamics, the binary expansion of $y$
and the binary expansion of $\tau^{n}x$ will have $m=\left\lfloor \log_{2}\varepsilon\right\rfloor $
digits in common. The number $m$ can be made arbitrarily large; the
subsequence will occur somewhere in the expansion. Put differently,
every finite-length string $\gamma\in\left\{ 0,1\right\} ^{<\omega}$
occurs as a prefix of (uncountably many) of the strings in $\left\langle x\right\rangle $.

In essence, that takes care of that, for the Bernoulli shift, at least,
if one is looking at it from the point of view of point dynamics.
As long as one thinks with the mind-set of points and their orbits,
there is not much more to be said. The above is, in a sense, a complete
description of the Bernoulli shift. But the introduction to this text
gave lie to this claim. If one instead looks at the shift in terms
of it's transfer operator acting on distributions, then much more
can be said. The spectrum of the transfer operator is non-trivial,
and the eigenfunctions are fractal, in general. This will be examined
more carefully, later; but for now, the topic of point dynamics in
the Bernoulli shift is exhausted. This is the end of the line.

The Bernoulli shift is not the only shift on Cantor space. And so,
onward ho.

\subsection{Beta shift}

The beta shift is similar to the Bernoulli shift, replacing the number
2 by a constant real-number value $1<\beta\le2$. It can be defined
as 

\begin{equation}
T_{\beta}(x)=\begin{cases}
\beta x & \mbox{ for }0\le x<\frac{1}{2}\\
\beta\left(x-\frac{1}{2}\right) & \mbox{ for }\frac{1}{2}\le x\le1
\end{cases}\label{eq:downshift}
\end{equation}
This map, together with similar maps, is illustrated in figure \ref{fig:Iterated-piece-wise-linear}
below. 

Just as the Bernoulli shift generates a sequence of digits, so does
the beta shift: write
\begin{equation}
k_{n}=\begin{cases}
0 & \mbox{ if }0\le T_{\beta}^{n}\left(x\right)<\frac{1}{2}\\
1 & \mbox{ if }\frac{1}{2}\le T_{\beta}^{n}\left(x\right)\le1
\end{cases}\label{eq:down-bits}
\end{equation}
Given the symbolic dynamics, one can reconstruct the original value
whenever $1<\beta$ as
\[
x=\frac{k_{0}}{2}+\frac{1}{\beta}\left(\frac{k_{1}}{2}+\frac{1}{\beta}\left(\frac{k_{2}}{2}+\frac{1}{\beta}\left(\frac{k_{3}}{2}+\frac{1}{\beta}\left(\cdots\right)\right)\right)\right)
\]
Written this way, the $T_{\beta}(x)$ clearly acts as a shift on this
sequence:
\[
T_{\beta}\left(x\right)=\frac{k_{1}}{2}+\frac{1}{\beta}\left(\frac{k_{2}}{2}+\frac{1}{\beta}\left(\frac{k_{3}}{2}+\frac{1}{\beta}\left(\frac{k_{4}}{2}+\frac{1}{\beta}\left(\cdots\right)\right)\right)\right)
\]
In this sense, this shift is ``exactly solvable'': the above provides
a closed-form solution for iterating and un-iterating the sequence.

Multiplying out the above sequence, one obtains the so-called ``$\beta$-expansion''
of a real number $x$, namely the series
\begin{equation}
x=\frac{1}{2}\,\sum_{n=0}^{\infty}\frac{k_{n}}{\beta^{n}}\label{eq:shift series}
\end{equation}
The bit-sequence that was extracted by iteration can be used to reconstruct
the original real number. Setting $\beta=2$ in eqn \ref{eq:down-bits}
gives the Bernoulli shift: $T_{2}\left(x\right)=b\left(x\right)$.

Unlike the Bernoulli shift, not every possible bit-sequence occurs
in this system. It is a subshift of the full shift: it is a subset
of $\left\{ 0,1\right\} ^{\omega}$ that is invariant under the action
of $T_{\beta}$. The structure of this shift is explored in detail
in a later section. 

\subsection{Associated polynomial}

The iterated shift can also be written as a finite sum. A later section
will be devoted entirely to the properties of this sum. Observe that
\[
T_{\beta}(x)=\beta\left(x-\frac{k_{0}}{2}\right)
\]
and that
\[
T_{\beta}^{2}(x)=\beta^{2}x-\frac{\beta}{2}\left(\beta k_{0}+k_{1}\right)
\]
and that 
\[
T_{\beta}^{3}(x)=\beta^{3}x-\frac{\beta}{2}\left(\beta^{2}k_{0}+\beta k_{1}+k_{2}\right)
\]
The general form is then:
\begin{equation}
T_{\beta}^{p}(x)=\beta^{p}x-\frac{\beta}{2}\sum_{m=0}^{p-1}k_{m}\beta^{p-m-1}\label{eq:iterated shift}
\end{equation}
Since the $k_{m}$ depend on both $\beta$ and on $x$, and are only
piece-wise continuous functions, this is not a true polynomial. It
does provide a polynomial-like representation with a range of interesting
properties.

\subsection{Density Visualizations}

\begin{figure}
\caption{Beta-shift Bifurcation Diagram\label{fig:Undershift-Bifurcation-Diagram}}

\medskip{}

\begin{centering}
\includegraphics[width=1\columnwidth]{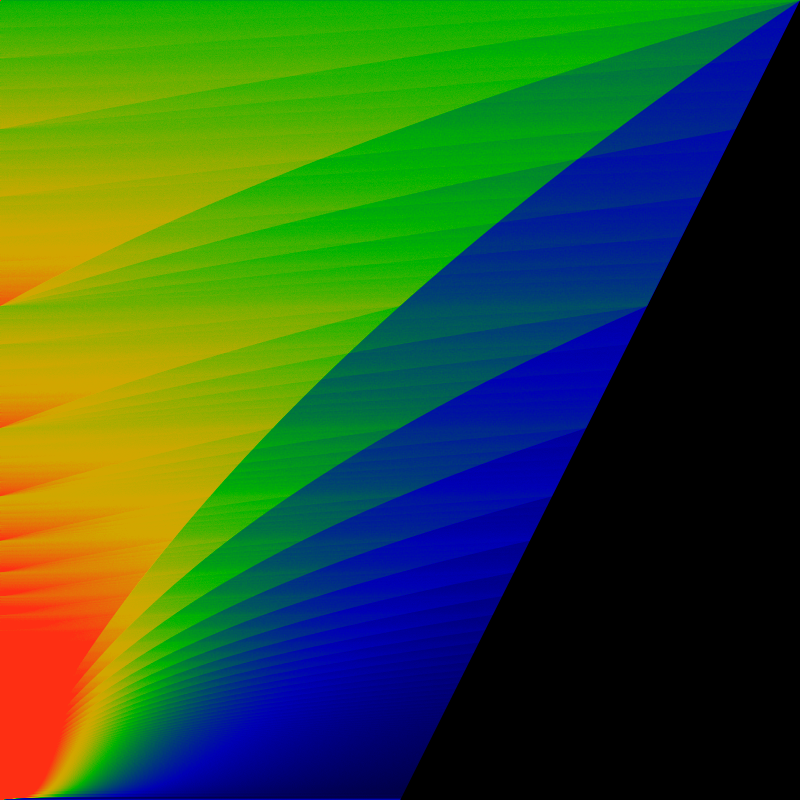}
\par\end{centering}
\medskip{}

This figure shows the density $\rho_{\beta}(y)$, rendered in color.
The constant $\beta$ is varied from 1 at the bottom to 2 at the top;
whereas $y$ runs from 0 on the left to 1 on the right. Thus, a fixed
value of $\beta$ corresponds to a horizontal slice through the diagram.
The color green represents values of $\rho_{\beta}(y)\approx1.0$,
while red represents $\rho_{\beta}(y)\gtrsim2$ and blue-to-black
represents $\rho_{\beta}(y)\apprle0.5$. The lines forming the fan
shape are not straight. The precise form is examined in a later section,
given by a variant of the polynomial in eqn \ref{eq:iterated shift}.
The discontinuities in this figure are more clearly highlighted in
the next figure, \ref{fig:Midpoint-Trace}.

\rule[0.5ex]{1\columnwidth}{1pt}
\end{figure}

The long-term dynamics of the $\beta$-shift can be visualized by
means of a bifurcation diagram. The idea of a bifurcation diagram
gained fame with the Feigenbaum map (shown in figure \ref{fig:Logistic-Map-Bifurcation}).
The same idea is applied here: track orbits over long periods of time,
and see where they go. This forms a density, which can be numerically
explored by histogramming. This is shown in figure \ref{fig:Undershift-Bifurcation-Diagram}.

When this is done for the $\beta$-shift, one thing becomes immediately
apparent: there are no actual ``bifurcations'', no ``islands of
stability'', no extended period-doubling regions, none of the stuff
normally associated with the Feigenbaum map. Although there are periodic
orbits, these form a set of measure zero: the iteration produces purely
chaotic motion for almost all values of $x$ and all values of $\beta>1$.
In this sense, the beta transform provides a clean form of ``pure
chaos'',\footnote{Formal mathematics distinguishes between many different kinds of chaotic
number sequences: those that are ergodic, those that are weakly or
strongly Bernoulli, weakly or strongly mixing. The beta transform
is known to be ergodic,\cite{Renyi57} weakly mixing\cite{Parry60}
and weakly Bernoulli.\cite{Dajani97}} without the pesky ``islands of stability'' popping up intermittently.

The visualization of the long-term dynamics is done by generating
a histogram and then taking the limit. The unit interval is divided
into a fixed sequence of equal-sized bins; say, a total of $N$ bins,
so that each is $1/N$ in width. Pick a starting $x$, and then iterate:
if, at the $n$'th iteration, one has that $j/N\le b_{\beta}^{n}(x)<(j+1)/N$,
then increment the count for the $j$'th bin. After a total of $M$
iterations, let $c(j;M)$ be the count in the $j$'th bin. This count
is the histogram. In the limit of a large number of iterations, as
well as small bin sizes, one obtains a distribution: 
\[
\rho(y;x)=\lim_{N\to\infty}\lim_{M\to\infty}\frac{c(j;M)}{M}\mbox{ for }\frac{j}{N}\le y<\frac{j+1}{N}
\]
This distribution depends on the initial value $x$ chosen for the
point to be iterated; a ``nice'' distribution results when one averages
over all starting points:
\[
\rho(y)=\int_{0}^{1}\rho(y;x)\,dx
\]
Numerically, this integration can be achieved by randomly sampling
a large number of starting points. By definition, $\rho(y)$ is a
probability distribution:
\[
1=\int_{0}^{1}\rho(x)\,dx
\]

Probability distributions are the same thing as measures; they assign
a density to each point on the unit interval. It can be shown that
this particular distribution is invariant under iteration, and thus
is often called the invariant measure, or sometimes the Haar measure.

\begin{figure}
\caption{Beta-shift Density Distribution\label{fig:Undershift-Density-Distribution}}

\begin{centering}
\includegraphics[width=1\columnwidth]{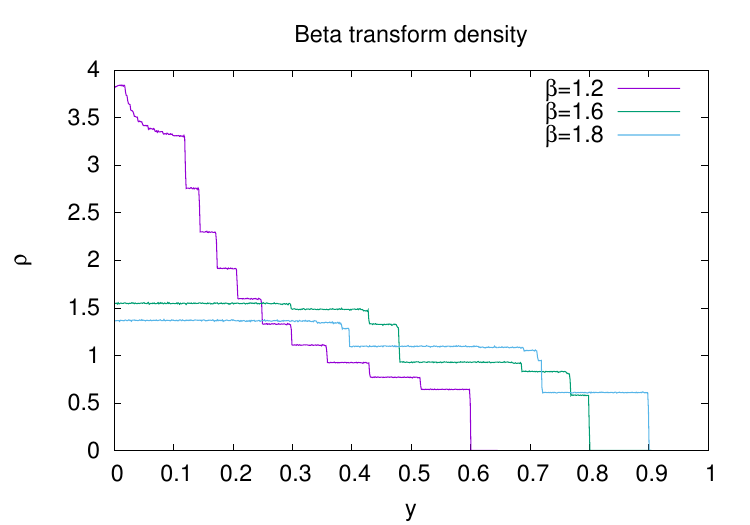}
\par\end{centering}
The above figure shows three different density distributions, for
$\rho_{1.2}(y)$, $\rho_{1.6}(y)$ and $\rho_{1.8}(y)$, calculated
numerically. These are obtained by histogramming a large number of
point trajectories, as described in the text. The small quantities
of jitter are due to a finite number of samples. To generate this
figure, a total of $M=4000$ iterations were performed, using randomly
generated arbitrary-precision floats (using the Gnu GMP package),
partitioned into $N=800$ bins, and sampled 24000 times (or 30 times
per bin) to perform the averaging integral. It will later be seen
that the discontinuities in this graph occur at the ``iterated midpoints''
$m_{p}=T_{\beta}^{p}\left(\beta/2\right)$. The flat plateaus are
not quite flat, but are filled with microscopic steps. There is a
discontinuous step at every $p$; these are ergodically distributed,
\emph{i.e.} dense in the interval, so that there are steps everywhere.
This is the general case; for special cases, when the midpoint has
a finite orbit, then there are a finite number of perfectly flat plateaus.
The first such example occurs at $\beta=\left(1+\sqrt{5}\right)/2=\varphi$
the Golden Ratio. In this case, there are only two such plateaus.

\rule[0.5ex]{1\columnwidth}{1pt}
\end{figure}

For each fixed $\beta$, one obtains a distinct distribution $\rho_{\beta}(y)$.
The figure \ref{fig:Undershift-Density-Distribution} illustrates
some of these distributions for a selection of fixed $\beta$. Note
that, for $\beta<1$, the distribution is given by $\rho_{\beta}(y)=\delta(y)$,
a Dirac delta function, located at $y=0$.

\begin{figure}
\caption{Midpoint Trace\label{fig:Midpoint-Trace}}

\includegraphics[width=1\columnwidth]{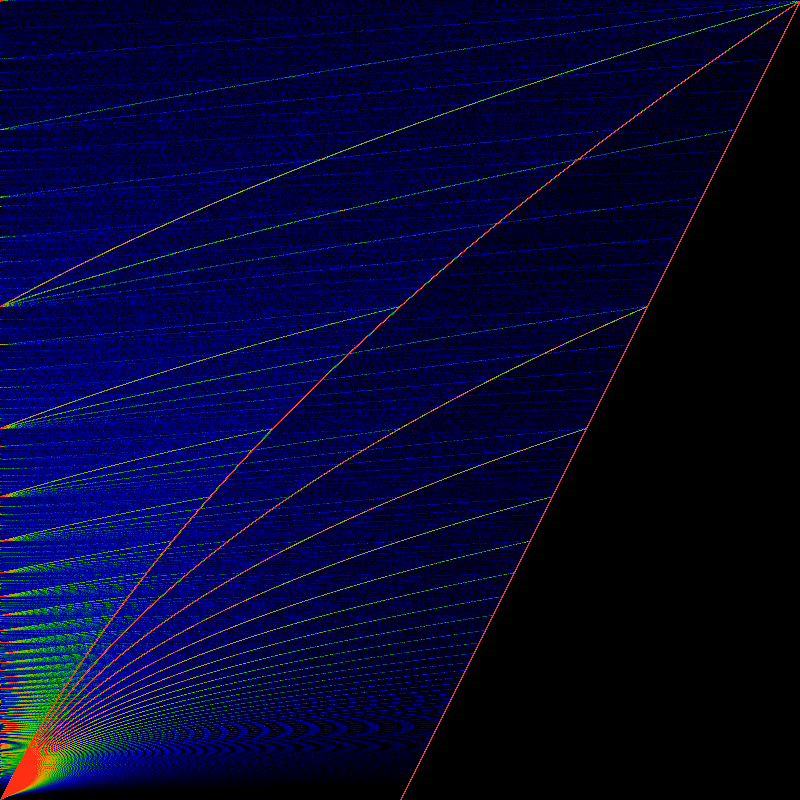}

Traces of midpoint iteration. Each horizontal line corresponds to
a fixed $\beta$, with $\beta$ running fom 1 at the bottom, to 2
at the top. At each fixed $\beta$, the midpoint $x=1/2$ is iterated
to generate $T_{\beta}^{n}(1/2)$. At each such location (from left
to right, of 0 to 1), the corresponding pixel is given a color assignment,
fading from red, through a rainbow, to black, as $n$ increases. This
is a variant of \ref{fig:Undershift-Bifurcation-Diagram}, highlighting
the edges. A formal analysis of the traceries begins with figure \ref{fig:Location-of-Midpoints}

\rule[0.5ex]{1\columnwidth}{1pt}
\end{figure}

The general trend of the distributions, as a function of $\beta$,
can be visualized with a Feigenbaum-style ``bifurcation diagram'',
shown in figure \ref{fig:Undershift-Bifurcation-Diagram}. This color-codes
each distribution $\rho_{\beta}(y)$ and arranges them in a stack;
a horizontal slice through the diagram corresponds to $\rho_{\beta}(y)$
for a fixed value of $\beta$. A related visualization is in \ref{fig:Midpoint-Trace},
which highlights the discontinuities in \ref{fig:Undershift-Bifurcation-Diagram}.
Periodic orbits appear wherever the traceries in this image intersect.
A characterization of these orbits occupies a large portion of this
text. 

\subsection{Tent Map}

\begin{figure}
\caption{Tent Map Bifurcation Diagram\label{fig:Tent-Map-Bifur}}

\begin{centering}
\includegraphics[width=1\columnwidth]{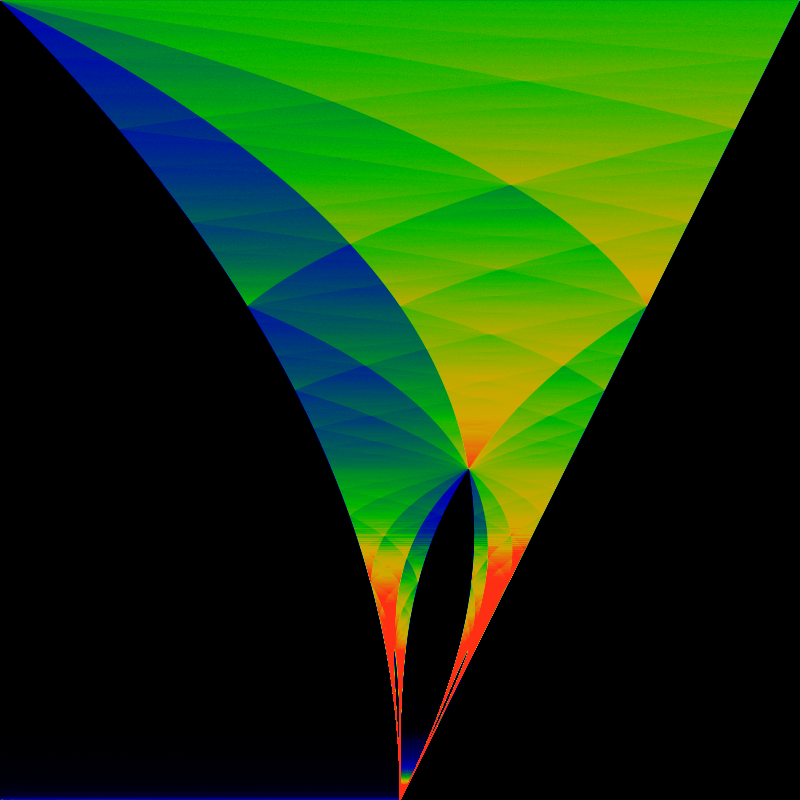}
\par\end{centering}
The bifurcation diagram for the tent map. The value of $\beta$ runs
from 1 at the bottom of the image, to 2 at the top. The color scheme
is adjusted so that green represents the average value of the distribution,
red represents areas of more than double the average value, while
blue shows those values that are about half the average value. Note
that this is a different color scheme than that used in figure \ref{fig:Undershift-Bifurcation-Diagram};
that scheme would obliterate the lower half of this figure in red.

\qquad{}The black areas represent parts of the iterated range that
are visited at most a finite number of times. To the right, the straight
boundary indicates that after one iteration, points in the domain
$\beta/2\le x\le1$ are never visited. To the left, points in the
domain $0\le x\le\beta\left(1-\beta/2\right)$ are never visited more
than a finite number of times.

\qquad{}The figure can be imagined to be a superposition of a countable
number of copies of figure \ref{fig:Undershift-Bifurcation-Diagram},
each drawn out so as to terminate in a point, but separated into distinct
arms. Each copy is recapitulated in the Feigenbaum bifurcation diagram.

\rule[0.5ex]{1\columnwidth}{1pt}
\end{figure}

The tent map is a closely related iterated map, given by iteration
of the function
\[
v_{\beta}(x)=\begin{cases}
\beta x & \mbox{ for }0\le x<\frac{1}{2}\\
\beta\left(1-x\right) & \mbox{ for }\frac{1}{2}\le x\le1
\end{cases}
\]
Its similar to the beta shift, except that the second arm is reflected
backwards, forming a tent. The bifurcation diagram is shown in figure
\ref{fig:Tent-Map-Bifur}. Its is worth contemplating the similarities
between this, and the corresponding beta shift diagram. Clearly, there
are a number of shared features.

\subsection{Logistic Map}

The logistic map is related to the tent map, and is given by iteration
of the function
\[
f_{\beta}\left(x\right)=2\beta x(1-x)
\]
It essentially replaces the triangle forming the tent map with a parabola
of the same height. That is, the function is defined here so that
the the same value of $\beta$ corresponds to the same height for
all three maps. Although the heights of the iterators have been aligned
so that they match, each exhibits rather dramatically different dynamics.
The $\beta$-transform has a single fixed point for $\beta<1$, and
then explodes into a fully chaotic regime above that. By contrast,
the logistic map maintains a single fixed point up to $\beta=3/2$,
where it famously starts a series of period-doubling bifurcations.
The onset of chaos is where the bifurcations come to a limit, at $\beta=3.56995/2=1.784975$.
Within this chaotic region are ``islands of stability'', which do
not appear in either the $\beta$-transform, or in the tent map. The
tent map does show a period-doubling regime, but in this region, there
are no fixed points: rather, the motion is chaotic, but confined to
multiple arms. At any rate, the period doubling occurs at different
values of $\beta$ than for the logistic map.

\begin{figure}
\caption{Logistic Map Bifurcation Diagram\label{fig:Logistic-Map-Bifurcation}}

\includegraphics[width=1\columnwidth]{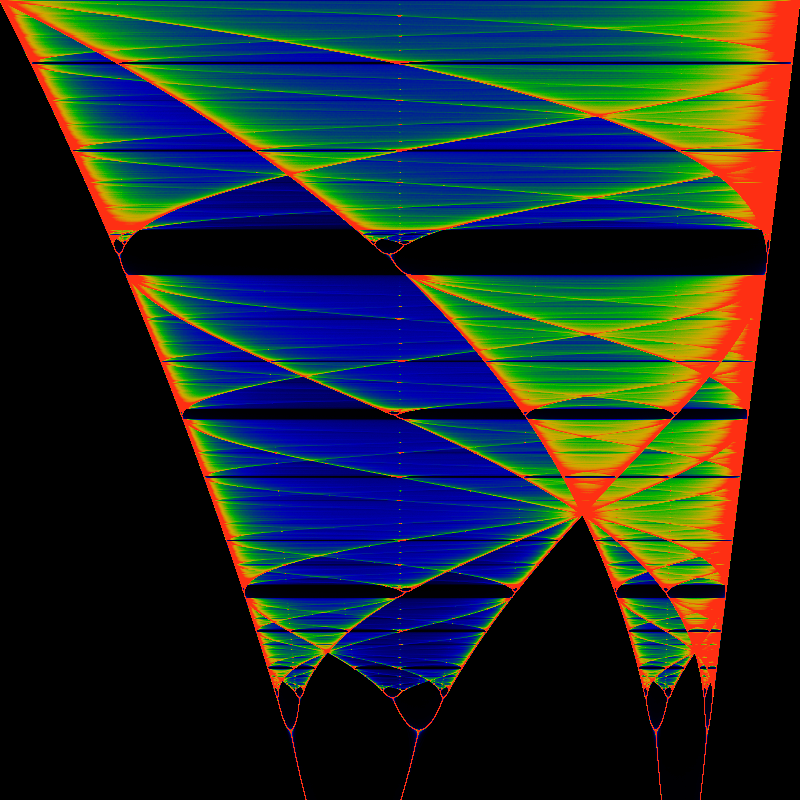}

The logistic map bifurcation diagram. The value of $\beta$ runs from
1.75 at the bottom of the image, to 2 at the top. The color scheme
is adjusted so that green represents the average value of the distribution,
red represents areas of more than double the average value, while
blue shows those values that are about half the average value. Clearly,
the orbits of the iterated points spend much of their time near the
edges of the diagram.

\qquad{}This is a very widely reproduced diagram. The goal here is
not to waste space reproducing it yet again, but to draw attention
to the similarities between this diagram, and figure \ref{fig:Tent-Map-Bifur}.
The tent can be smoothly deformed into the parabola; doing so separates
the multiple legs visible at the bottom of \ref{fig:Tent-Map-Bifur},
with each leg becoming a bifurcation branch. The central point, where
all legs come together, is preserved in both images. The curving arcs,
highlighted in figure \ref{fig:Midpoint-Trace}, are retained as well.
The legs of \ref{fig:Tent-Map-Bifur} are chaotic regions, and remain
so as the tent is deformed into a parabola. They are merely separated
from one-another. The legs do \emph{not} shrink to zero width. 

\rule[0.5ex]{1\columnwidth}{1pt}
\end{figure}

The bifurcation diagram is shown in figure \ref{fig:Logistic-Map-Bifurcation}.
Again, it is worth closely examining the similarities between this,
and the corresponding tent-map diagram, as well as the $\beta$-transform
diagram. Naively, it would seem that the general structure of the
chaotic regions are shared by all three maps. Thus, in order to understand
chaos in the logistic map, it is perhaps easier to study it in the
$\beta$-transform.

The general visual similarity between the figures \ref{fig:Undershift-Bifurcation-Diagram},
\ref{fig:Tent-Map-Bifur} and \ref{fig:Logistic-Map-Bifurcation}
should be apparent, and one can pick out and find visually similar
regions among these three illustrations. Formalizing this similarity
is a bit harder, but it can be done: all three of these maps are topologically
conjugate to one-another. This is perhaps surprising, but is based
on the observation that the ``islands of stability'' in the logistic
map are countable, and are in one-to-one correspondence with certain
``trouble points'' in the iterated beta transformation. These are
in turn in one-to-one correspondence with rational numbers. With a
slight distortion of the beta transformation, the ``trouble points''
can be mapped to the islands of stability, in essentially the same
way that phase locking regions (Arnold tongues) appear in the circle
map. This is examined in a later section; it is mentioned here only
to whet the appetite.

\subsection{Beta Transformation}

After exactly one iteration of the beta shift, all initial points
$\beta/2\le x\le1$ are swept up into the domain $0\le x<\beta/2$,
and never leave. Likewise, the range of the iterated beta-shift is
$0\le x<\beta/2$. Thus, an alternative representation of the beta
shift, filling the entire unit square, can be obtained by dividing
both the domain and range by $\beta/2$ to obtain the function
\begin{equation}
t_{\beta}(u)=\begin{cases}
\beta u & \mbox{ for }0\le u<\frac{1}{\beta}\\
\beta u-1 & \mbox{ for }\frac{1}{\beta}\le u\le1
\end{cases}\label{eq:beta transform}
\end{equation}
This can be written more compactly as $t_{\beta}\left(x\right)=\beta x\mod1$.
In this form, the function is named ``the beta-transform'', written
as the $\beta$-transformation, presenting a typesetting challenge
to search engines when used in titles of papers. The orbit of a point
$x$ in the beta-shift is identical to the orbit of a point $u=2x/\beta$
in the beta-transformation. Explicitly comparing to the beta-shift
of eqn \ref{eq:downshift}:
\begin{equation}
T_{\beta}^{n}\left(x\right)=\frac{\beta}{2}t_{\beta}^{n}\left(\frac{2x}{\beta}\right)\label{eq:xform-shift equivalence}
\end{equation}
The beta-shift and the $\beta$-transformation are essentially the
same function; this text works almost exclusively with the beta-shift,
and is thus idiosyncratic, as it flouts the more common convention
of working with the $\beta$-transformation. The primary reason for
this is a historical quirk, as this text was started before the author
became aware of the $\beta$-transformation.

After a single iteration of the tent map, a similar situation applies.
After one iteration, all initial points $\beta/2\le x\le1$ are swept
up into the domain $0\le x<\beta/2$. After a finite number of iterations,
all points $0<x\le\beta\left(1-\beta/2\right)$ are swept up, so that
the remaining iteration takes place on the domain $\beta\left(1-\beta/2\right)<x<\beta/2$.
It is worth defining a ``sidetent'' function, which corresponds
to the that part of the tent map in which iteration is confined. It
is nothing more than a rescaling of the tent map, ignoring those parts
outside of the above domain that wander away. The sidetent is given
by
\[
s_{\beta}(u)=\begin{cases}
\beta\left(u-1\right)+2 & \mbox{ for }0\le u<\frac{\beta-1}{\beta}\\
\beta\left(1-u\right) & \mbox{ for }\frac{\beta-1}{\beta}\le u\le1
\end{cases}
\]
Performing a left-right flip on the side-tent brings it closer in
form to the beta-transformation. The flipped version, replacing $u\to1-u$
is 
\[
f_{\beta}(u)=\begin{cases}
\beta u & \mbox{ for }0\le u<\frac{1}{\beta}\\
2-\beta u & \mbox{ for }\frac{1}{\beta}\le u\le1
\end{cases}
\]

\begin{figure}
\caption{Iterated piece-wise linear maps\label{fig:Iterated-piece-wise-linear}}

\begin{centering}
\includegraphics[width=1\columnwidth]{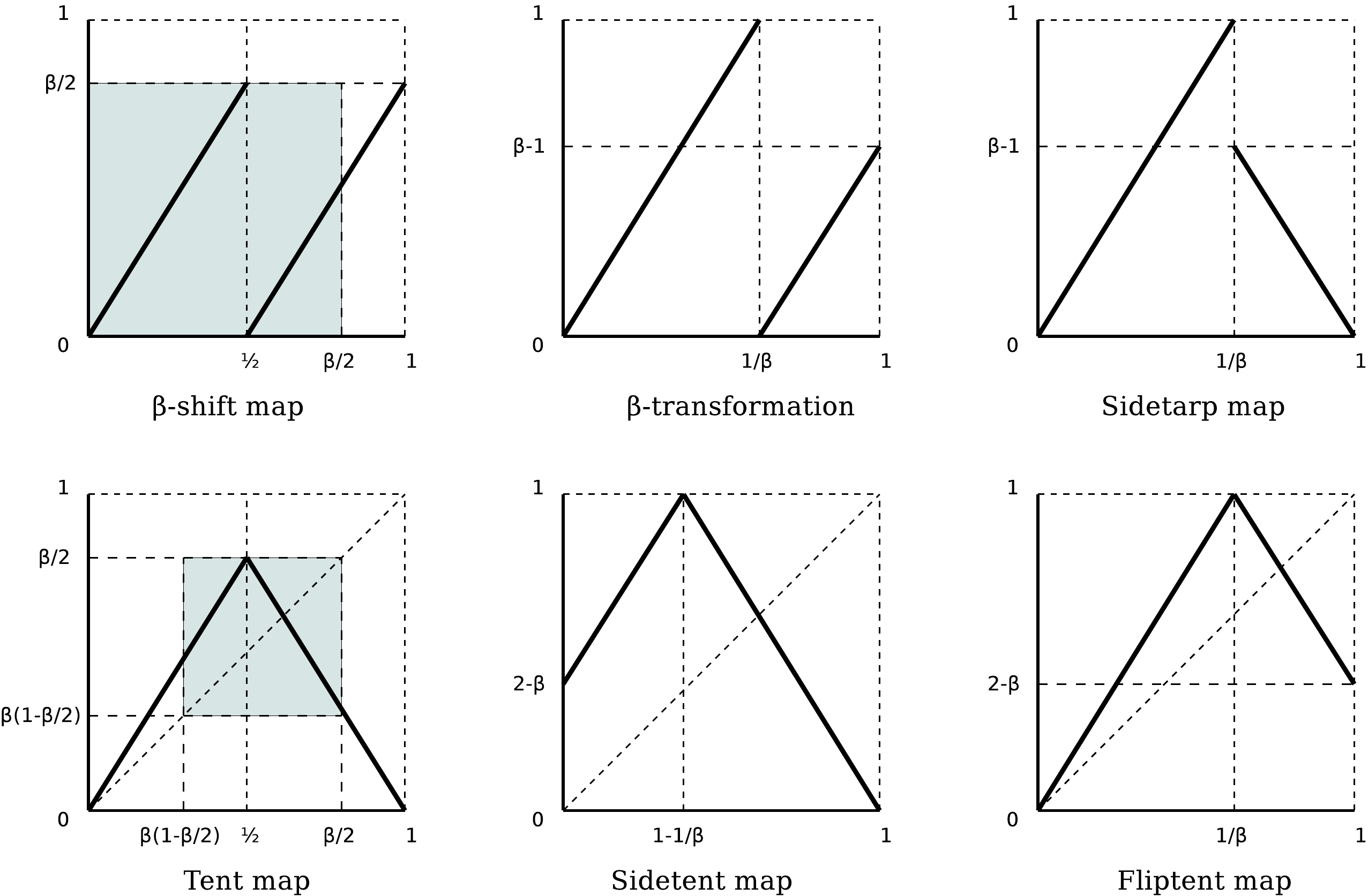}
\par\end{centering}
The beta shift map, shown in the upper left, generates orbits that
spend all of their time in the shaded area: a box of size $\frac{\beta}{2}\times\frac{\beta}{2}$.
Enlarging this box to the unit square gives the $\beta$-transformation.
The tent map resembles the beta shift, except that one arm is flipped
to make a tent-shape. After a finite number of iterations, orbits
move entirely in the shaded region; enlarging this region to be the
unit square gives the sidetent map. Flipping it left-right gives the
fliptent map. Although it is not trivially obvious, the fliptent map
and the sidetent map have the same orbits, and thus the same bifurcation
diagram.

\qquad{}The bottom three maps all have prominent fixed points and
periodic orbits, essentially because the diagonal intersects the map.
The top three maps have periodic orbits, but these occur only for
a countable number of $\beta$ values. General orbits are purely chaotic,
essentially because the diagonal does not intersect them. Note that
the slopes and the geometric proportions of all six maps are identical;
they are merely rearrangents of the same basic elements.

\rule[0.5ex]{1\columnwidth}{1pt}
\end{figure}

\begin{figure}
\caption{Beta transform and Side-tent\label{fig:Side-Shift-and-Side-tent}}

\begin{centering}
\includegraphics[width=0.32\columnwidth]{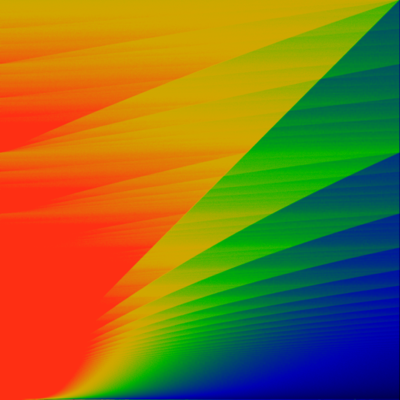}\,\includegraphics[width=0.32\columnwidth]{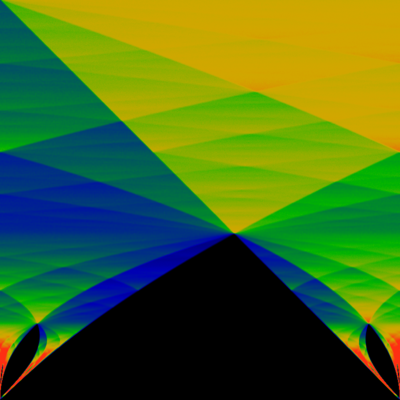}\,\includegraphics[width=0.32\columnwidth]{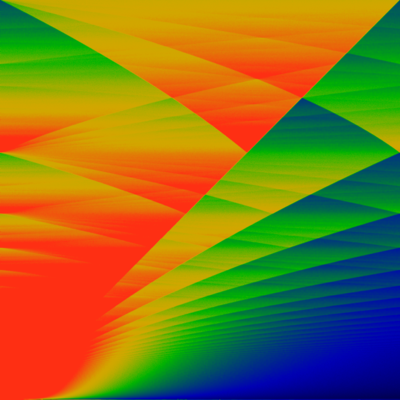}
\par\end{centering}
The left figure shows the bifurcation diagram for the $\beta$-transform,
as it is normally defined as the $\beta x\mod1$ map. It is the same
map as the beta shift, just rescaled to occupy the entire unit square.
In all other respects, it is identical to \ref{fig:Undershift-Bifurcation-Diagram}.

\qquad{}The middle figure is a similarly-rescaled tent map, given
the name ``side tent'' in the main text. It is essentially identical
to \ref{fig:Tent-Map-Bifur}, with the middle parts expanded and the
sides removed. In both figures, $\beta$ runs from 1 at the bottom
to 2 at the top. The right-hand-side figure is the ``sidetarp'',
clearly its an oddly-folded variant of the beta transform. 

\rule[0.5ex]{1\columnwidth}{1pt}
\end{figure}

The tent map (and the flipped tent) exhibits fixed points (periodic
orbits; mode-locking) for the smaller values of $\beta$. These can
be eliminated by shifting part of the tent downwards, so that the
diagonal is never intersected. This suggests the ``sidetarp'':
\[
a_{\beta}(u)=\begin{cases}
\beta u & \mbox{ for }0\le u<\frac{1}{\beta}\\
\beta\left(1-u\right) & \mbox{ for }\frac{1}{\beta}\le u\le1
\end{cases}
\]
The six different maps under consideration here are depicted in figure
\ref{fig:Iterated-piece-wise-linear}. It is interesting to compare
three of the bifurcation diagrams, side-by-side. These are shown in
figure \ref{fig:Side-Shift-and-Side-tent}.

\subsection{Dynamical Systems}

A brief review of dynamical systems is in order, as it provides a
coherent language with which to talk about and think about the beta-shift.
The technical reason for this is that a subshift $S\subset\left\{ 0,1\right\} ^{\omega}$
provides a more natural setting for the theory, and that a lot of
the confusion about what happens on the unit interval is intimately
entangled with the homomorphism \ref{eq:homomorphism} (or \ref{eq:shift series}
as the case may be). Disentangling the subshift from the homomorphism
provides a clearer insight into what phenomena are due to which component.

The review of dynamical systems here is more-or-less textbook-standard
material; it is included here only to provide a firm grounding for
later discussion.

The Cantor space $\left\{ 0,1\right\} ^{\omega}$ can be given a topology,
the product topology. The open sets of this topology are called ``cylinder
sets''. These are the infinite strings in three symbols: a finite
number of 0 and 1 symbols, and an infinite number of {*} symbols,
the latter meaning ``don't care''. Set union is defined location-by-location,
with $0\cup*=1\cup*=*$ and set intersection as $0\cap*=0$ and $1\cap*=1$.
Set complement exchanges 0 and 1 and leaves {*} alone: $\overline{0}=1$,
$\overline{1}=0$ and $\overline{*}=*$. The topology is then the
collection of all cylinder sets. Note that the intersection of any
finite number of cylinder sets is still a cylinder set, as is the
union of an infinite number of them. The product topology does \emph{not}
contain any ``points'': strings consisting solely of just 0 and
1 are \emph{not} allowed in the topology. By definition, topologies
only allow finite intersections, and thus don't provide any way of
constructing ``points''. Of course, points can always be added ``by
hand'', but doing so tends to generate a topology (the ``box topology'')
that is ``too fine''; in particular, the common-sense notions of
a continuous function are ruined by fine topologies. The product topology
is ``coarse''.

The Borel algebra, or sigma-algebra, takes the topology and also allows
set complement. This effectively changes nothing, as the open sets
are still the cylinder sets, although now they are ``clopen'', as
they are both closed and open.

Denote the Borel algebra by $\mathcal{B}$. A shift is now a map $T:\mathcal{B}\to\mathcal{B}$
that lops off the leading symbol of a given cylinder set. This is
provides strong theoretical advantages over working with ``point
dynamics'': confusions about counting points and orbits and defining
densities go away. This is done by recasting discussion in terms of
functions $f:\mathcal{B}\to\mathbb{R}$ from Borel sets to the reals
(or the complex numbers $\mathbb{C}$ or other fields, when this is
interesting). An important class of such functions are the measures.
These are functions $\mu:\mathcal{B}\to\mathbb{R}$ that are positive,
and are ``compatible'' with the sigma algebra, in that $\mu\left(A\cup B\right)=\mu\left(A\right)+\mu\left(B\right)$
whenever $A\cap B=\varnothing$ and (for product-space measures) that
$\mu\left(A\cap B\right)=\mu\left(A\right)\mu\left(B\right)$ for
all $A,B\in\mathcal{B}$. The measure of the total space $\Omega=\left\{ 0,1\right\} ^{\omega}$
is by convention unity: $\mu\left(\Omega\right)=1$.

The prototypical example of a measure is the Bernoulli measure, which
assigns probability $p$ to any string containing a single 0 and the
rest all {*}'s. By complement, a string containing a single 1 and
the rest all {*}'s has probability $1-p$. The rest follows from the
sigma algebra: a cylinder set consisting of $m$ zeros and $n$ ones
has measure $p^{m}\left(1-p\right)^{n}$. It is usually convenient
to take $p=1/2$, the ``fair coin''; the Bernoulli process is a
sequence of coin tosses.

The map given in equation \ref{eq:homomorphism} is a homomorphism
from the Cantor space to the unit interval. It extends naturally to
a map from the Borel algebra $\mathcal{B}$ to the algebra of intervals
on the unit interval. It is not an isomorphism: cylinder sets are
both open and closed, whereas intervals on the real number line are
either open, or closed (or half-open). It is convenient to take the
map as a map to closed intervals, so that it's a surjection onto the
reals, although usually, this detail does not matter. What does matter
is if one takes $p=1/2$, then the Bernoulli measure is preserved:
it is mapped onto the conventional measure on the real-number line.
Thus, the cylinder set $0***\cdots$ is mapped to the interval $\left[0,1/2\right]$
and $1***\cdots$ is mapped to $\left[1/2,1\right]$ and both have
a measure of $1/2$ and this extends likewise to all intersections
and unions. Points have a measure of zero. That is, the homomorphism
\ref{eq:homomorphism} preserves the fair-coin Bernoulli measure.

Much of what is said above still holds for subshifts. Recall, a subshift
$S$ is a subspace $S\subset\left\{ 0,1\right\} ^{\omega}$ that is
invariant under the shift $T$, so that $TS=S$. The space $S$ inherits
a topology from $\left\{ 0,1\right\} ^{\omega}$; this is the subspace
topology. The Borel algebra $\mathcal{B}$ is similarly defined, as
are measures. One can now (finally!) give a precise definition for
an invariant measure: it is a measure $\mu$ such that $\mu\circ T^{-1}=\mu$,
or more precisely, for which $\mu\left(T^{-1}\left(\sigma\right)\right)=\mu\left(\sigma\right)$
for almost all cylinder sets $\sigma\in S$. This is what shift invariance
looks like. Note carefully that $T^{-1}$ and not $T$ is used in
the definition. This is because $T^{-1}$ is a surjection while $T$
is not: every cylinder set $\sigma$ in the subshift ``came from
somewhere''; we want to define invariance for all $\sigma$ and not
just for some of them.

The $T^{-1}$ is technically called a ``pushforward'', and it defines
a linear operator $\mathcal{L}_{T}$ on the space $\mathcal{F}$ of
all functions $f:\mathcal{B}\to\mathbb{R}$. It is defined as $\mathcal{L}_{T}:\mathcal{F\to\mathcal{F}}$
by setting $\mathcal{L}_{T}:f\mapsto f\circ T^{-1}$. It is obviously
linear, in that $\mathcal{L}_{T}\left(af+bg\right)=a\mathcal{L}_{T}\left(f\right)+b\mathcal{L}_{T}\left(g\right)$.
This pushforward is canonically called the ``transfer operator''
or the ``Ruelle-Frobenius-Perron operator''. Like any linear operator,
it has a spectrum. The precise spectrum depends on the space $\mathcal{F}$.

The canonical example is again the Bernoulli shift. For this, we invoke
the inverse of the mapping of eqn \ref{eq:homomorphism} so that $f:\left[0,1\right]\to\mathbb{R}$
is a function defined on the unit interval, instead of $f:\mathcal{B}\to\mathbb{R}$.
When $\mathcal{F}$ is the space of real-analytic functions on the
unit interval, that is, the closure of the space of all polynomials
in $x\in\left[0,1\right]$, then the spectrum of $\mathcal{L}_{T}$
is discrete. It consists of the Bernoulli polynomials $B_{n}\left(x\right)$
corresponding to an eigenvalue of $2^{-n}$. That is,$\mathcal{L}_{T}B_{n}=2^{-n}B_{n}$.
Note that $B_{0}\left(x\right)=1$ is the invariant measure on the
full shift. For the space of square-integrable functions $f:\left[0,1\right]\to\mathbb{R}$,
the spectrum of $\mathcal{L}_{T}$ is continuous, and consists of
the unit disk in the complex plane; the corresponding eigenfunctions
are fractal. Even more interesting constructions are possible; the
Minkowski question mark function provides an example of a measure
on $\left\{ 0,1\right\} ^{\omega}$ that is invariant under the shift
defined by the Gauss map $h\left(x\right)\mapsto\frac{1}{x}-\left\lfloor \frac{1}{x}\right\rfloor $.
That is, as a measure, it solves $\mathcal{L}_{T}?^{\prime}=?^{\prime}$
with $?$ the Minkowski question mark function, and $?^{\prime}$
it's derivative; note that the derivative is ``continuous nowhere''.
This rather confusing idea (of something being ``continuous nowhere'')
can be completely dispelled by observing that it is well-defined on
all cylinder sets in $\mathcal{B}$ and is finite on all of them –
not only finite, but less than one, as any good measure must obey.

These last examples are mentioned so as to reinforce the idea that
working with $\mathcal{B}$ instead of the unit interval $\left[0,1\right]$
really does offer some strong conceptual advantages. They also reinforce
the idea that the Bernoulli shift is not the only ``full shift''.
In the following text, we will be working with subshifts, primarily
the beta-shift, but will draw on ideas from the above so as to make
rigorous statements about measurability and invariance, without having
to descend into either \emph{ad hoc} hand-waving or provide painfully
difficult (and confusing) reasoning about subsets of the real-number
line.

\subsection{Beta Transformation Literature Review and References}

The $\beta$-transformation, in the form of $t_{\beta}\left(x\right)=\beta x\mod1$
has been well-studied over the decades. The beta-expansion \ref{eq:down-bits}
was introduced by A. Renyi\cite{Renyi57} in 1957, who demonstrates
the existence of the invariant measure. The ergodic properties of
the transform were proven by W. Parry\cite{Parry60} in 1960, who
also shows that the system is weakly mixing. 

An explicit expression for the invariant measure was obtained independently
by A.O. Gel'fond\cite{Gelfond59} in 1959, and by W. Parry\cite{Parry60},
as a summation of step functions
\begin{equation}
\nu_{\beta}\left(y\right)=\frac{1}{F}\sum_{n=0}^{\infty}\frac{\varepsilon_{n}\left(y\right)}{\beta^{n}}\label{eq:invariant measure}
\end{equation}
where $\varepsilon_{n}$ is the digit sequence 
\begin{equation}
\varepsilon_{n}\left(y\right)=\begin{cases}
1 & \mbox{ if }y<t_{\beta}^{n}\left(1\right)\\
0 & \mbox{ otherwise}
\end{cases}\label{eq:Gelfond-Parry digits}
\end{equation}
and $F$ is a normalization constant. By integrating $\varepsilon_{n}\left(y\right)$
under the sum, the normalization is given by
\[
F=\sum_{n=0}^{\infty}\frac{t_{\beta}^{n}\left(1\right)}{\beta^{n}}
\]

Analogous to the way in which a dyadic rational $p/2^{n}$ has two
different binary expansions, one ending in all-zeros, and a second
ending in all-ones, so one may also ask if and when a real number
$x$ might have more than one $\beta$-expansion (for fixed $\beta$).
In general, it can; N. Sidorov shows that almost every number has
a continuum of such expansions!\cite{Sidorov03} This signals that
the beta shift behaves rather differently from the Cantor set in it's
embedding.

Conversely, the ``univoke numbers'' are those values of $\beta$
for which there is only one, unique expansion for $x=1$. These are
studied by De Vries.\cite{DeVries06}

The $\beta$-transformation has been shown to have the same ergodicity
properties as the Bernoulli shift.\cite{Dajani97} The fact that the
beta shift, and its subshifts are all ergodic is established by Climenhaga
and Thompson.\cite{Clim10}

An alternative to the notion of ergodicity is the notion of universality:
a $\beta$-expansion is universal if, for any given finite string
of bits, that finite string occurs somewhere in the expansion. This
variant of universality was introduced by Erdös and Komornik\cite{Erdos98}.
Its is shown by N. Sidorov that almost every $\beta$-expansion is
universal.\cite{Sidorov02} Conversely, there are some values of $\beta$
for which rational numbers have purely periodic $\beta$-expansions;\cite{Adamczewski10}
all such numbers are Pisot numbers.\cite{Schmidt80}

The symbolic dynamics of the beta-transformation was analyzed by F.
Blanchard\cite{Blanchard89}. A characterization of the periodic points
are given by Bruno Maia\cite{Maia07}. A discussion of various open
problems with respect to the beta expansion is given by Akiyama.\cite{Akiyama09} 

When the beta expansion is expanded to the entire real-number line,
one effectively has a representation of reals in a non-integer base.
One may ask about arithmetic properties, such as the behavior of addition
and multiplication, in this base - for example, the sum or product
of two $\beta$-integers may have a fractional part! Bounds on the
lengths of these fractional parts, and related topics, are explored
by multiple authors.\cite{Guimond01,Julien06,Hbaib13}

Certain values of $\beta$ – generally, the Pisot numbers, generate
fractal tilings,\cite{Thurston89,Berthe05,Ito05,Adamczewski10,Akiyama09}
which are generalizations of the Rauzy fractal. An overview, with
common terminology and definitions is provided by Akiyama.\cite{Akiyama17}
The tilings, sometimes called (generalized) Rauzy fractals, can be
thought of as living in a direct product of Euclidean and $p$-adic
spaces.\cite{Berthe04}

The set of finite beta-expansions constitutes a language, in the formal
sense of model theory and computer science. This language is recursive
(that is, decidable by a Turing machine), if and only if $\beta$
is a computable real number.\cite{Simonsen11}

The zeta function, and a lap-counting function, are given by Lagarias\cite{Lagarias94}.
The Hausdorff dimension, the topological entropy and general notions
of topological pressure arising from conditional variational principles
is given by Daniel Thompson\cite{Thompson10}. A proper background
on this topic is given by Barreira and Saussol\cite{Barreira01}. 

None of the topics or results cited above are made use of, or further
expanded on, or even touched on in the following. This is not intentional,
but rather a by-product of different goals.\pagebreak{}

\section{Transfer operators}

Given any map from a space to itself, mapping points to points, the
pushforward maps distributions to distributions. The pushforward is
a linear operator, called the transfer operator or the Ruelle–Frobenius–Perron
operator. The spectrum of this operator, broken down into eigenfunctions
and eigenvalues, can be used to understand the time evolution of a
given density distribution. The invariant measure is an eigenstate
of this operator, it is the eigenstate with eigenvalue one. There
are other eigenstates; these are explored in this section.

Restricting to the unit interval, given an iterated map $f:\left[0,1\right]\to\left[0,1\right]$,
the transfer operator acting on a distribution $\rho:\left[0,1\right]\to\mathbb{R}$
is defined as
\[
\left[\mathcal{L}_{f}\rho\right]\left(y\right)=\sum_{x=f^{-1}(y)}\frac{\rho(x)}{\left|f^{\prime}(x)\right|}
\]
The next subsection gives an explicit expression for this, when $f$
is the $\beta$-transform. After that, a subsection reviewing the
invariant measure, and then a discussion of some other eigenfunctions.

\subsection{The $\beta$-shift Transfer Operator}

This text works primarily with the $\beta$-shift, instead of the
more common $\beta$-transform. These two are more-or-less the same
thing, differing only by scale factors, as given in eqn. \ref{eq:xform-shift equivalence}.
The transfer operators are likewise only superficially different,
being just rescalings of one-another; both are given below.

The transfer operator the beta-shift map $T_{\beta}(x)$ is
\[
\left[\mathcal{L}_{\beta}f\right]\left(y\right)=\begin{cases}
\frac{1}{\beta}\left[f\left(\frac{y}{\beta}\right)+f\left(\frac{y}{\beta}+\frac{1}{2}\right)\right] & \mbox{ for }0\le y\le\beta/2\\
0 & \mbox{ for }\beta/2<y\le1
\end{cases}
\]
or, written more compactly
\begin{equation}
\left[\mathcal{L}_{\beta}f\right]\left(y\right)=\frac{1}{\beta}\left[f\left(\frac{y}{\beta}\right)+f\left(\frac{y}{\beta}+\frac{1}{2}\right)\right]\Theta\left(\frac{\beta}{2}-y\right)\label{eq:xfer oper}
\end{equation}
where $\Theta$ is the \href{https://en.wikipedia.org/wiki/Heaviside_step_function}{Heaviside step function}.
The transfer operator for the beta-transform map $t_{\beta}(x)$ is
\[
\left[\mathcal{M}_{\beta}f\right]\left(y\right)=\frac{1}{\beta}\left[f\left(\frac{y}{\beta}\right)+f\left(\frac{y}{\beta}+\frac{1}{\beta}\right)\Theta\left(\beta-1-y\right)\right]
\]

The density distributions graphed in figure \ref{fig:Undershift-Density-Distribution}
are those functions satisfying 
\begin{equation}
\left[\mathcal{L}_{\beta}\mu\right]\left(y\right)=\mu\left(y\right)\label{eq:FP-eigenvector}
\end{equation}
That is, the $\mu\left(y\right)$ satisfies 
\begin{equation}
\mu\left(y\right)=\frac{1}{\beta}\left[\mu\left(\frac{y}{\beta}\right)+\mu\left(\frac{y}{\beta}+\frac{1}{2}\right)\right]\Theta\left(\frac{\beta}{2}-y\right)\label{eq:eigen-eqn}
\end{equation}
Likewise, the Gelfond-Parry measure of eqn \ref{eq:invariant measure}
satisfies
\[
\left[\mathcal{M}_{\beta}\nu_{\beta}\right]\left(y\right)=\nu_{\beta}\left(y\right)
\]
Recall that $\mu\left(y\right)=\frac{2}{\beta}\nu_{\beta}\left(\frac{2y}{\beta}\right)\Theta\left(\frac{\beta}{2}-y\right)$;
the two invariant measures are just scaled copies of one-another.
Both are normalized so that $\int_{0}^{1}\mu\left(y\right)dy=\int_{0}^{1}\nu_{\beta}\left(y\right)dy=1$.

Both of these invariant measures are the Ruelle-Frobenius-Perron (RFP)
eigenfunctions of the corresponding operators, as they correspond
to the largest eigenvalues of the transfer operators, in each case
being the eigenvalue one.

More generally, one is interested in characterizing the spectrum
\[
\left[\mathcal{L}_{\beta}\rho\right]\left(y\right)=\lambda\rho\left(y\right)
\]
for eigenvalues $\left|\lambda\right|\le1$ and eigenfunctions $\rho\left(y\right)$.
Solving this equation requires choosing a space of functions in which
to work. Natural choices include piece-wise continuous smooth functions
(piece-wise polynomial functions), any of the Banach spaces, and in
particular, the space of square-integrable functions. In general,
the spectrum will be complex-valued: eigenvalues will be complex numbers.

If a distribution $\rho\left(y\right)$ is nonzero on the interval
$\left[\beta/2,1\right]$, the operator $\mathcal{L}_{\beta}$ will
map it to one that is zero on this interval. Thus, it makes sense
to restrict oneself to densities that are nonzero only on $\left[0,\beta/2\right]$.
When this is done, eqn \ref{eq:xfer oper} has the slightly more convenient
form
\[
\left[\mathcal{L}_{\beta}f\right]\left(y\right)=\frac{1}{\beta}\left[f\left(\frac{y}{\beta}\right)\Theta\left(m_{0}-y\right)+f\left(\frac{y}{\beta}+\frac{1}{2}\right)\Theta\left(m_{1}-y\right)\right]
\]
with $m_{0}=\beta/2$ and $m_{1}=\beta\left(\beta-1\right)/2=T_{\beta}\left(m_{0}\right)$.
It is always the case that $m_{1}<m_{0}$ and so the second term above
vanished on the interval $\left[m_{1},m_{0}\right]$. This can be
gainfully employed in a variety of settings; typically to write $\mathcal{L}_{\beta}f$
on $\left[m_{1},m_{0}\right]$ as a simple rescaling of $\mathcal{L}_{\beta}f$
on $\left[0,m_{1}\right]$.

This equation can be treated as a recurrence relation; setting $\mathcal{L}_{\beta}f=f$
gives the $\lambda=1$ eigenstate. Performing this recursion gives
exactly the densities shown in figure \ref{fig:Undershift-Density-Distribution}.
Computationally, these are much cheaper to compute than trying to
track a scattered cloud of points; the result is free of stochastic
sampling noise. This density is the Ruelle–Frobenius–Perron eigenstate;
an explicit expression was given by Gelfond and by Parry, as described
in the next section. 

\subsection{The Gelfond–Parry measure}

An explicit expression for the solution to $\mathcal{M}_{\beta}\nu=\nu$
was given by Gelfond\cite{Gelfond59} and by Parry\cite{Parry60}.
It is the expression given by eqn \ref{eq:invariant measure}. Unfortunately,
I find the Russian original of Gelfond's article unreadable, and Parry's
work, stemming from his PhD thesis, is not available online. Therefore,
it is of some interest to provide a proof suitable for the current
text. A generalization of this proof, stated in terms of a Borel algebra,
is used in the subsequent section to construct general eigenfunctions.

There are two routes: either a direct verification that eqn \ref{eq:invariant measure}
is correct, or a derivation of eqn \ref{eq:invariant measure} from
geometric intuition. The direct verification is useful for practical
purposes; the geometric construction, as a stretch-cut-stack map,
provides insight. Both are given.

\begin{figure}
\caption{Gelfond–Parry Normalization\label{fig:Gelfond=002013Parry-Normalization}}

\includegraphics[width=1\columnwidth]{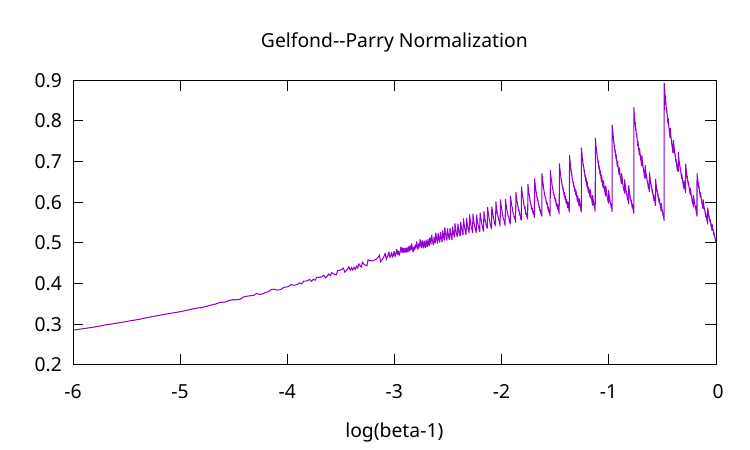}

The above figure shows $1/F$ for the normalization constant $F=\sum_{n=0}^{\infty}t_{n}\beta^{-n}$
as a function of $\beta$. The horizonal axis is stretched out using
$\log\left(\beta-1\right)$ so as to amplify the behavior as $\beta\to1$.
One has that $F\to\infty$ in this limit (so $1/F\to0$); the curve
suggests just how catastrophic that limit is. A graph of $F$ vs.
$\beta$, without the rescaling of the horizontal axis, is shown in
figure \ref{fig:Normalization-Integral}.

\rule[0.5ex]{1\columnwidth}{1pt}
\end{figure}

The Gelfond–Parry measure includes a normalization factor. It will
be of recurring interest, and so a graph f it is presented in figure
\ref{fig:Gelfond=002013Parry-Normalization}.

\subsubsection{Direct verification}

A direct verification of correctness is done below, explicitly showing
all steps in laborious detail. It's not at all difficult; just a bit
hard on the eyes.

As before, let $t\left(x\right)\equiv t_{\beta}\left(x\right)=\beta x\mod1$
be the $\beta$-transformation of eqn \ref{eq:beta transform}, and
$t^{n}\left(x\right)$ the iterated transformation. Let $\Theta(x)$
be the Heaviside step function as always, and to keep notation brief,
let $t_{n}\equiv t^{n}\left(1\right)$. The Gelfond-Parry measure
is then 
\[
\nu\left(y\right)=\frac{1}{F}\sum_{n=0}^{\infty}\frac{\Theta\left(t_{n}-y\right)}{\beta^{n}}
\]
where the normalization $F$ is given by 
\[
F=\sum_{n=0}^{\infty}\frac{t_{n}}{\beta^{n}}
\]
The transfer operator $\mathcal{M}$ for the beta-transformation is
slightly more convenient to work with than $\mathcal{L}$ for this
particular case. It is given by 
\[
\left[\mathcal{M}f\right]\left(y\right)=\frac{1}{\beta}\left[f\left(\frac{y}{\beta}\right)+f\left(\frac{y+1}{\beta}\right)\Theta\left(\beta-1-y\right)\right]
\]
and we wish to verify that $\mathcal{M}\nu=\nu$. Plugging in directly,
\begin{align*}
P=\nu\left(\frac{y}{\beta}\right) & +\nu\left(\frac{y+1}{\beta}\right)\Theta\left(\beta-1-y\right)=\\
 & =\frac{1}{F}\sum_{n=0}^{\infty}\frac{1}{\beta^{n}}\left(\Theta\left(t_{n}-\frac{y}{\beta}\right)+\Theta\left(t_{n}-\frac{y+1}{\beta}\right)\Theta\left(\beta-1-y\right)\right)\\
 & =\frac{1}{F}\sum_{n=0}^{\infty}\frac{1}{\beta^{n}}\left(\Theta\left(\beta t_{n}-y\right)+\Theta\left(\beta t_{n}-1-y\right)\right)
\end{align*}
The second line follows from the first, since $\Theta\left(ax\right)=\Theta\left(x\right)$
for all constants $a$, and the $\Theta\left(\beta-1-y\right)$ can
be safely dropped, since $t_{n}\le1$ for all $n$. These terms simplify,
depending on whether $t_{n}$ is small or large. Explicitly, one has
\begin{align*}
\Theta\left(\beta t_{n}-1-y\right)=0 & \mbox{ if }\beta t_{n}-1<0\\
\Theta\left(\beta t_{n}-y\right)=1 & \mbox{ if }\beta t_{n}-1>0
\end{align*}
and so
\[
P=\frac{1}{F}\sum_{n=0}^{\infty}\frac{1}{\beta^{n}}\left(\Theta\left(1-\beta t_{n}\right)\Theta\left(\beta t_{n}-y\right)+\Theta\left(\beta t_{n}-1\right)\left(1+\Theta\left(\beta t_{n}-1-y\right)\right)\right)
\]
These can be collapsed by noting that 
\begin{align*}
\beta t_{n}=t_{n+1} & \mbox{ if }\beta t_{n}-1<0\\
\beta t_{n}-1=t_{n+1} & \mbox{ if }\beta t_{n}-1>0
\end{align*}
and so
\begin{align*}
P & =\frac{1}{F}\sum_{n=0}^{\infty}\frac{1}{\beta^{n}}\left(\Theta\left(1-\beta t_{n}\right)\Theta\left(t_{n+1}-y\right)+\Theta\left(\beta t_{n}-1\right)\left(1+\Theta\left(t_{n+1}-y\right)\right)\right)\\
 & =\frac{1}{F}\sum_{n=0}^{\infty}\frac{1}{\beta^{n}}\left(\Theta\left(t_{n+1}-y\right)\left[\Theta\left(1-\beta t_{n}\right)+\Theta\left(\beta t_{n}-1\right)\right]+\Theta\left(\beta t_{n}-1\right)\right)\\
 & =\frac{1}{F}\sum_{n=0}^{\infty}\frac{1}{\beta^{n}}\left(\Theta\left(t_{n+1}-y\right)+\Theta\left(\beta t_{n}-1\right)\right)\\
 & =\beta\nu\left(y\right)-\frac{\beta}{F}+\frac{1}{F}\sum_{n=0}^{\infty}\frac{\Theta\left(\beta t_{n}-1\right)}{\beta^{n}}\\
 & =\beta\nu\left(y\right)
\end{align*}
The last sum on the right is just the $\beta$-expansion for 1. That
is, the $\beta$-expansion of $x$ is
\[
x=\sum_{n=0}^{\infty}\frac{\Theta\left(\beta t^{n}\left(x\right)-1\right)}{\beta^{n+1}}
\]
This is just eqn \ref{eq:shift series} written in a different way
(making use of the equivalence \ref{eq:xform-shift equivalence}).
Thus $P=\beta\nu\left(y\right)$ and so $\mathcal{M}\nu=\nu$ as claimed.

\subsubsection{Stretch–Cut–Stack Map}

The measure can be geometrically constructed and intuitively understood
as the result of the repeated application of a bakers-map-style stretch
and squash operation. The idea is to iterate $\nu_{n+1}=\mathcal{M}_{\beta}\nu_{n}$,
starting with $\nu_{0}=1$ and then take the limit $n\to\infty$;
the result is $\nu_{n}\to\nu$ in the limit, given as above. The basic
stretch-cut-stack operation provides intuition as to why $\nu$ is
invariant. The proof that it converges as desired is nearly identical
to the (shorter) proof above.

Consider a first approximation that $\nu_{0}$ is constant on the
interval $\left[0,1\right]$, so that $\nu_{0}\left(y\right)=1$.
The operation of $\mathcal{M}$ acting on this is to stretch it out
to the interval $\left[0,\beta\right]$, chop off the $\left[1,\beta\right]$
part, move it to $\left[0,\beta-1\right]$, stack it on top, thus
doubling the density in this region. The doubling, though, is partly
counteracted by the stretching, which thins out the density to $1/\beta$
uniformly over the entire interval $\left[0,\beta\right]$. This operation
preserves the grand-total measure on the unit interval. Writing $\beta-1=t_{1}=t\left(t_{0}\right)=t\left(1\right)$
for the first iterate of the endpoint $t_{0}=1$, this stretch, cut
and stack operation should result in 
\[
\nu_{1}\left(y\right)=\begin{cases}
\frac{2}{\beta} & \mbox{for }y\in\left[0,t_{1}\right]\\
\frac{1}{\beta} & \mbox{for }y\in\left[t_{1},t_{0}\right]
\end{cases}
\]
The stretch–cut–stack operation is the intuitive, geometrical explanation
for obtaining the first iterate. The same result is obtained algebraically,
by writing $\nu_{1}=\mathcal{M}\nu_{0}$ and then plugging and chugging:
\[
\nu_{1}\left(y\right)=\frac{1}{\beta}\left[1+\Theta\left(t_{1}-y\right)\right]
\]

The derivation of a recursive formula for $\nu_{n}$ follows along
the line of the algebraic proof given in the previous section, this
time, limiting the sums to finite order. The general form is 

\[
\nu_{n}=\rho_{n}-\sum_{k=0}^{n-1}c_{n-k}\nu_{k}
\]
which is preserved under iteration as $\nu_{n+1}=\mathcal{M}\nu_{n}$.
The two parts are $\rho_{n}$ defined as
\[
\rho_{n}=\sum_{k=0}^{n}\frac{1}{\beta^{k}}\Theta\left(t_{k}-y\right)
\]
This converges, up to overall normalization, to the invariant measure:
$\rho_{n}\to F\nu$ with $F$ the normalization constant from before.
The $c_{n}$ are constants, defined as
\begin{equation}
c_{n}=1-\sum_{k=0}^{n-1}\frac{\Theta\left(\beta t_{k}-1\right)}{\beta^{k+1}}\label{eq:iteration const}
\end{equation}
The sum in $c_{n}$ is a partial $\beta$-expansion for 1. Each term
can be recognized as a bit from the bit-expansion:
\[
\Theta\left(\beta t_{k}-1\right)=d_{n}\left(\frac{1}{2}\right)=k_{n}\left(\frac{\beta}{2}\right)=\varepsilon_{n}\left(\frac{1}{\beta}\right)
\]
with $k_{n}$ as defined in eqn \ref{eq:down-bits}, $\varepsilon_{n}$
defined as in eqn \ref{eq:Gelfond-Parry digits} and $d_{n}$ as in
eqn \ref{eq:rescaled-bitseq} (apologies for the variety of notations;
each is ``natural'' in a specific context.) The $\beta$-expansion
of a real number $x$ is as given in eqn \ref{eq:shift series}. In
the present context, this is the expansion for $x=\beta/2$, or, after
rescaling, the expansion for $x=1$:
\[
1=\sum_{k=0}^{\infty}\frac{\Theta\left(\beta t_{k}-1\right)}{\beta^{k+1}}
\]
Thus, $\left|c_{n}\right|<\beta^{-n+1}$ and $c_{n}\to0$ as $n\to\infty$.
The sum in $c_{n}$ is a partial $\beta$-expansion for 1. Thus, $\left|c_{n}\right|<\beta^{-n+1}$
and $c_{n}\to0$ as $n\to\infty$.

Each iteration preserves the measure; this is built into the construction.
That is, $\int_{0}^{1}\nu_{n}\left(y\right)dy=1$ for all $n$. Plugging
through gives a curious identity:
\[
1=\int_{0}^{1}\nu_{n}\left(y\right)dy=\sum_{k=0}^{n}\frac{t_{k}}{\beta^{k}}-\sum_{k=1}^{n}c_{k}
\]
which holds for all $n$. This is readily verified by further plugging
through; the underlying identity is that
\[
t_{n}=\beta^{n}c_{n}
\]
This follows by noticing that each $\Theta\left(\beta t_{k}-1\right)$
records a decision to decrement, or not, the product $\beta t_{k}$
occurring in the iteration $t_{k+1}=\beta t_{k}\mod1$.

The Gelfond–Parry normalization was the $n\to\infty$ limit of the
first sum: $F=\sum_{k=0}^{\infty}\beta^{-k}t_{k}$.

\subsubsection{Stacking Generic Functions}

A slightly different calculation results if one works with an an arbitrary
function $\nu_{0}\left(y\right)$. The iteration to obtain $\nu_{n}=\mathcal{M}^{n}\nu_{0}$
proceeds much as before, with only minor modifications needed to obtain
the general form.

The trick is to track two distinct travellers: one that travels with
$\Theta\left(t_{k}-y\right)$ and another that travels with the bitsequence
$b_{k}=\Theta\left(\beta t_{k}-1\right)$. The general solution has
the form
\begin{align}
\nu_{n}\left(y\right) & =a_{n}\left(y\right)+\sum_{k=1}^{n}h_{nk}\left(y\right)\Theta\left(t_{k}-y\right)\nonumber \\
a_{n}\left(y\right) & =\sum_{k=0}^{n-1}b_{k}e_{nk}\left(y\right)\label{eq:iterated xfer oper}
\end{align}
The functions $h_{nk}$ and $e_{nk}$ can be solved for recursively;
the recursive relations are simple enough that they can be rolled
up as series summations. These are given by

\begin{align*}
h_{nk}\left(y\right) & =\frac{1}{\beta^{k}}a_{n-k}\left(t_{1}-\frac{\beta t_{k}-1}{\beta^{k}}+\frac{y}{\beta^{k}}\right)\quad\mbox{ for }n\ge k
\end{align*}
which express $h_{nk}$ in terms of the function $a_{n}$. This is
given as a recursive series, obtained by iterating

\[
e_{nk}\left(y\right)=\frac{1}{\beta^{n}}\left[\nu_{0}\left(1-\frac{t_{k}}{\beta^{k}}+\frac{y}{\beta^{n}}\right)+\Theta\left(k\right)\sum_{m=1}^{n-k-1}\beta^{m}a_{m}\left(t_{1}-\frac{\beta t_{k}-1}{\beta^{k}}+\frac{\beta^{m}y}{\beta^{n}}\right)\right]
\]
The starting point for iteration is $a_{0}=\nu_{0}$.

The function $a_{n}$ is polynomial if $\nu_{0}$ is, and of the same
degree; $a_{n}$ is analytic, if $\nu_{0}$ is, and so on. The function
$h_{nk}$ is likewise, on the domain $0\le y\le t_{k}$; the discontinuities
in $\nu_{n}$ are entirely due to the $\Theta\left(t_{k}-y\right)$
term.

The proper calculation of the $n\to\infty$ limit of $a_{n}$ remains
a mystery. This is the primary obstacle to constructing general eigenfunctions
from this series.

\paragraph{Verification}

When $\nu_{0}=1$, all three functions $a_{n},e_{nk}$ and $h_{nk}$
become constants. In this case, $a_{n}\to1/F$ as $n\to\infty$, with
$F$ the normalization constant as before. Defining $\alpha_{n}=\beta^{n}a_{n}$,
the recursion relation takes the curious form 
\[
\alpha_{n}=1+\sum_{k=1}^{n-1}b_{k}\sum_{m=0}^{n-k-1}\alpha_{m}
\]
This is an integer sequence: each $\alpha_{n}$ is an integer, starting
with $\alpha_{0}=\alpha_{1}=1$. This has the form of a generalized
Fibonacci sequence. Denote the partial sum as $s_{i}=\sum_{m=0}^{i}\alpha_{m}$.
This implies $\alpha_{i}=s_{i}-s_{i-1}$ and so
\begin{align*}
\alpha_{n} & =1+b_{1}s_{n-2}+b_{2}s_{n-3}+\cdots\\
s_{n} & =1+s_{n-1}+b_{1}s_{n-2}+b_{2}s_{n-3}+\cdots
\end{align*}
If the orbit is of finite length $p$, so that $b_{k}=0$ for $k>p$,
then this can be recast directly into generalized Fibonacci form,
by defining $f_{n}=s_{n}+\sum_{k=1}^{p}b_{k}$. The recursion relation
is then
\[
f_{n}=f_{n-1}+b_{1}f_{n-2}+\cdots+b_{p}f_{n-p}
\]
For example, the orbit generated by $\beta=\varphi=1.618\cdots$ the
golden mean has $p=1,b_{1}=1$ and the sequence is $\alpha_{m}=1,1,2,3,5,8,13,\cdots$.
Finite orbits and generalized Fibonacci sequences will be treated
at length in the next chapter. One of the interesting properties is
that $\beta=\lim_{n\to\infty}f_{n}/f_{n-1}$ holds in the general
case, and not just for the golden mean. The bitsequences are self-describing;
this is byproduct from the identity $\beta=\sum_{k=0}^{\infty}b_{k}\beta^{-k}$.

\subsubsection{Scaled iteration\label{subsec:Scaled-iteration}}

It is convenient to rescale the iterate $\mathcal{M}^{n}$ by a constant
factor, so that instead, one is examining a sequence of functions
$\mathcal{M}\mu_{n}=\lambda\mu_{n+1}$. The iterated operator then
has the form
\begin{align*}
\mu_{n}\left(y\right) & =c_{n}^{\prime}\left(y\right)+\sum_{k=1}^{n}\frac{\Theta\left(t_{k}-y\right)}{\left(\lambda\beta\right)^{k}}c_{n-k}^{\prime}\left(t_{1}-\frac{\beta t_{k}-1}{\beta^{k}}+\frac{y}{\beta^{k}}\right)\\
c_{n}^{\prime}\left(y\right) & =\nu_{0}\left(y\right)+\frac{1}{\left(\lambda\beta\right)^{n}}\sum_{k=1}^{n-1}b_{k}\left[\nu_{0}\left(1-\frac{t_{k}}{\beta^{k}}+\frac{y}{\beta^{n}}\right)+\sum_{m=1}^{n-k-1}\left(\lambda\beta\right)^{m}c_{m}^{\prime}\left(t_{1}-\frac{\beta t_{k}-1}{\beta^{k}}+\frac{\beta^{m}y}{\beta^{n}}\right)\right]
\end{align*}
This is identical to the earlier forms when $\lambda=1$. The starting
point for iteration is $c_{0}^{\prime}=\nu_{0}$. For For $\nu_{0}=1$
and $\lambda=1$, one regains the Gelfond–Parry distribution in the
limit of $n\to\infty$. In this limit, $c_{n}^{\prime}\to1/F$.

\subsection{Analytic Gelfond–Parry function\label{subsec:Rotated-Renyi-Parry-function}}

The technique above can be repeated verbatim for a ``rotated'' or
``coherent'' function
\begin{equation}
\nu_{\beta;z}\left(y\right)=\sum_{n=0}^{\infty}z^{n}\frac{\Theta\left(t_{n}-y\right)}{\beta^{n}}\label{eq:generalized-eigenfunction}
\end{equation}
for a given complex-valued $z$. No changes are required, and the
result can be read off directly:
\begin{align*}
\left[\mathcal{M}\nu_{\beta;z}\right]\left(y\right) & =\frac{\nu_{\beta;z}\left(y\right)}{z}-\frac{1}{z}+\sum_{n=0}^{\infty}z^{n}\frac{\Theta\left(\beta t_{n}-1\right)}{\beta^{n+1}}\\
 & =\frac{\nu_{\beta;z}\left(y\right)}{z}+C\left(\beta;z\right)
\end{align*}
with $C\left(\beta;z\right)$ being a constant independent of $y$.
If there are values of $\beta$ and/or $z$ at which $C\left(\beta;z\right)=0$,
then this becomes the eigenequation for $\mathcal{M}$.

The eigenfunction for $\mathcal{L}$ is the same, up to rescaling
of $y\mapsto\beta x/2$. Recycling notation slightly, write

\begin{equation}
v_{\beta;z}\left(x\right)=\sum_{n=0}^{\infty}\frac{d_{n}\left(x\right)}{\beta^{n}}z^{n}\label{eq:coherent-first-order}
\end{equation}
where the $d_{n}\left(x\right)$ are exactly the same digits as defined
by Parry, just rescaled for the beta-shift. That is, 
\begin{equation}
d_{n}\left(x\right)=\varepsilon_{n}\left(\frac{2x}{\beta}\right)=\Theta\left(\frac{\beta}{2}t_{n}-x\right)=\Theta\left(T^{n}\left(\frac{\beta}{2}\right)-x\right)=\begin{cases}
1 & \mbox{if }x<T^{n}\left(\frac{\beta}{2}\right)\\
0 & \mbox{otherwise}
\end{cases}\label{eq:rescaled-bitseq}
\end{equation}
where $T$ the beta-shift map of eqn \ref{eq:downshift} and eqn \ref{eq:xform-shift equivalence}
is used. The iterated end-point becomes the iterated midpoint:
\[
t^{n}\left(1\right)=\frac{2}{\beta}T^{n}\left(\frac{\beta}{2}\right)
\]

Holding both $\beta$ and $x$ fixed, the summation is clearly convergent
(and holomorphic in $z$) for complex numbers $z$ within the disk
$\left|z\right|<\beta$. The eigenequation has the same form:

\[
\left[\mathcal{L}_{\beta}v_{\beta;z}\right]\left(x\right)=\frac{1}{z}v_{\beta;z}\left(x\right)+C\left(\beta;z\right)
\]
where $C\left(\beta;z\right)$ is a constant independent of $x$.
Numeric verification reveals we were a bit glib: $C\left(\beta;z\right)$
is a constant for $x\le\beta/2$ and is zero otherwise! (This is normal;
$\mathcal{L}$ was defined in such a way that it is always exactly
zero for $x>\beta/2$.)

The interesting limit is where $\left|z\right|\to\beta$ and so its
convenient to re-express $C$ in terms of $\zeta=z/\beta$, so that
everything is mapped to the unit disk. With some rearrangements, one
obtains
\begin{equation}
E\left(\beta;\zeta\right)\equiv\zeta\beta C\left(\beta;\zeta\beta\right)=-1+\zeta\sum_{n=0}^{\infty}\zeta^{n}d_{n}\left(\frac{1}{2}\right)\label{eq:holomorphic-disk}
\end{equation}
Given that $d_{n}\left(1/2\right)=\Theta\left(\beta t_{n}-1\right)$,
the above can be recognized as the rotated/coherent form of eqn \ref{eq:iteration const}
in the $n\to\infty$ limit. The primary task is to characterize the
zeros of $E\left(\beta;\zeta\right)$. This is a straight-forward
sum to examine numerically; results will be presented in the section
after the next.

\subsection{Analytic ergodics}

This section proposes that the entire $\beta$-subshift can be tied
together with a single holomorphic equation. The holomorphic equation
effectively provides a continuum (\emph{i.e.} uncountable number)
of distinct relationships between different parts of the subshift.
This can be interpreted either as a form of interaction across the
subshift, or as a kind of mixing. Given the nature of the relationship,
the moniker of ``fundamental theorem of analytic ergodics'' is an
amusing name to assign to the result.

The constant term can be independently derived through a different
set of manipulations. Explicitly plugging in eqn \ref{eq:coherent-first-order}
into the transfer operator yields

\begin{align*}
C\left(\beta;z\right)= & \left[\mathcal{L}_{\beta}v_{\beta;z}\right]\left(y\right)-\frac{v_{\beta;z}\left(y\right)}{z}\\
= & \sum_{n=0}^{\infty}\frac{z^{n}}{\beta^{n}}\left[\frac{1}{\beta}\left[d_{n}\left(\frac{y}{\beta}\right)+d_{n}\left(\frac{y}{\beta}+\frac{1}{2}\right)\right]\Theta\left(\frac{\beta}{2}-y\right)-\frac{1}{z}d_{n}\left(y\right)\right]
\end{align*}
Replacing $z$ by $\zeta=z/\beta$ gives
\begin{align*}
E\left(\beta;\zeta\right)= & \zeta\beta C\left(\beta;\zeta\beta\right)\\
= & \sum_{n=0}^{\infty}\zeta^{n}\left[\zeta\left[d_{n}\left(\frac{y}{\beta}\right)+d_{n}\left(\frac{y}{\beta}+\frac{1}{2}\right)\right]\Theta\left(\frac{\beta}{2}-y\right)-d_{n}\left(y\right)\right]
\end{align*}
This is holomorphic on the unit disk $\left|\zeta\right|<1$, as each
individual $d_{n}$ is either zero or one; there won't be any poles
inside the unit disk. Note that $d_{n}\left(y\right)=0$ for all $y>\beta/2$,
and so one may pull out the step function to write 
\[
E\left(\beta;\zeta\right)=\sum_{n=0}^{\infty}\zeta^{n}\left[\zeta d_{n}\left(\frac{y}{\beta}\right)+\zeta d_{n}\left(\frac{y}{\beta}+\frac{1}{2}\right)-d_{n}\left(y\right)\right]\Theta\left(\frac{\beta}{2}-y\right)
\]
confirming the earlier observation that $E\left(\beta;\zeta\right)$
vanishes for all $y>\beta/2$.

The bottom equation holds without assuming that $E\left(\beta;\zeta\right)$
is independent of $y$. However, we've already proven that it is;
and so a simplified expression can be given simply by picking a specific
$y$. Setting $y=0$, noting that $d_{n}\left(0\right)=1$ and canceling
terms, one obtains eqn. \ref{eq:holomorphic-disk} again.

Staring at the right-hand side of the sum above, it is hardly obvious
that it should be independent of $y$. In a certain sense, this is
not ``one equation'', this holds for a continuum of $y$, for all
$0\le y\le1$. It is an analytic equation tying together the entire
subshift. For each distinct $y$, it singles out three completely
different bit-sequences out of the subshift, and ties them together.
It is a form of mixing. Alternately, a form of interaction: the bit-sequences
are not independent of one-another; they interact. This section attempts
to make these notions more precise.

The tying-together of seemingly unrelated sequences seems somehow
terribly important. It is amusing to suggest that this is a kind of
``fundamental theorem of analytic ergodics''.

For such a claim, it is worth discussing the meaning at length, taking
the effort to be exceptionally precise and verbose, perhaps a bit
repetitive. Equation \ref{eq:downshift} defined a map, the $\beta$-shift.
Equation \ref{eq:down-bits} defined a bit-sequence, the $\beta$-expansion
of a real number $0\le x\le1$, where equation \ref{eq:shift series}
is the definition of the $\beta$-expansion. The set of all such bit-sequences
defines the shift. To emphasize this point, its best to compare side-by-side.
Copying equation \ref{eq:down-bits}, one bitsequence records the
orbit of $x$ relative to the midpoint: 
\[
k_{n}\left(x\right)=\begin{cases}
0 & \mbox{ if }0\le T_{\beta}^{n}(x)<\frac{1}{2}\\
1 & \mbox{ if }\frac{1}{2}\le T_{\beta}^{n}(x)\le1
\end{cases}
\]
 while a different bitsequence records the orbit of the midpoint,
relative to $x$:
\[
d_{n}\left(x\right)=\begin{cases}
1 & \mbox{if }x<T^{n}\left(\frac{\beta}{2}\right)=T^{n+1}\left(\frac{1}{2}\right)\\
0 & \mbox{otherwise}
\end{cases}
\]
The iterations are running in opposite directions; this is as appropriate,
since the the transfer operator was a pushforward.

It is useful to return to the language of sigma algebras and cylinder
sets, as opposed to point dynamics. Recall that the Borel algebra
$\mathcal{B}$ was defined as the sigma algebra, the collection of
all cylinder sets in the product topology of $\left\{ 0,1\right\} ^{\omega}$.
A subshift is a subset $\mathcal{S}\subset\mathcal{B}$ together with
a map $T:\mathcal{S}\to\mathcal{S}$ that lops off the leading symbol
of a given cylinder set but otherwise preserves the subshift: $T\mathcal{S}=\mathcal{S}$.
The inverted map $T^{-1}$ is a pushforward, in that it defines the
transfer operator, a linear operator $\mathcal{L}_{T}:\mathcal{F\to\mathcal{F}}$
on the space $\mathcal{F}$ of all functions $f:\mathcal{S}\to\mathbb{R}$;
explicitly, it is given by $\mathcal{L}_{T}:f\mapsto f\circ T^{-1}$.
Insofar as the $d_{n}$ arose in the exploration of the transfer operator,
it is not surprising that the shift is acting ``backwards''.

The problem with the language of point dynamics is that one cannot
meaningfully write $T^{-n}\left(x\right)$ for a real number, a point
$x$, at least, not without severe contortions that lead back to the
Borel algebra. Not for lack of trying; the $T^{-n}\left(x\right)$
is called the ``Julia set'' (to order $n$) of $x$: it is the preimage,
the set of all points that, when iterated, converge onto $x$.

Can the analytic relation be restated in terms of cylinder sets? Yes,
and it follows in a fairly natural way. The first step is to extend
$d_{n}$ to a map $d_{n}:\mathcal{S}\to\left[0,1\right]$ into the
unit interval, as opposed to being a single bit. Let $\mu:\mathcal{B}\to\left[0,1\right]$
be the canonical Bernoulli measure. Using the Bernoulli mapping \ref{eq:homomorphism},
the interval $\left[0,T^{n}\left(\beta/2\right)\right]$ maps to some
cylinder; call it $\Delta_{n}$. Then, given some cylinder $A\in\mathcal{S}$,
define

\[
d_{n}\left(A\right)=\mu\left(A\cap\Delta_{n}\right)
\]
The rotated (pre-)measure is extended likewise: 
\[
\nu\left(A\right)=\sum_{n=0}^{\infty}\zeta^{n}d_{n}\left(A\right)
\]
with $\zeta=z/\beta$ as before, recovering the Parry measure by setting
$z=1$.

The Parry measure should be invariant under the action of $T^{-1}:\mathcal{S}\to\mathcal{S}$,
and otherwise yield eqn \ref{eq:holomorphic-disk}. Let's check. The
proof will mirror the one of the previous section. Again, it is convenient
to use the $\beta$-transform $t_{\beta}$ instead of the $\beta$-shift
$T_{\beta}$. This is primarily a conceptual convenience; the subshift
is more easily visualized in terms of the mod 1 map. Otherwise, the
same notation is used, but rescaled, so that $\Delta_{n}$ is the
cylinder corresponding to the interval $\left[0,t^{n}\left(1\right)\right]$.

Recall that for every $A\in\mathcal{S}$ and every $y\in A$, one
will find that $y/\beta\in t^{-1}\left(A\right)$ and, whenever $y\le\beta-1$
that also $\left(y+1\right)/\beta\in t^{-1}\left(A\right)$. Thus,
$t^{-1}\left(A\right)$ naturally splits into two parts: the cylinder
that maps to $\left[0,1/\beta\right]$, call it $D$, and the complement
$\overline{D}$.

The pushforward action is then 
\begin{align*}
\nu\left(t^{-1}\left(A\right)\right)= & \sum_{n=0}^{\infty}\zeta^{n}\mu\left(\Delta_{n}\cap t^{-1}\left(A\right)\right)\\
= & \sum_{n=0}^{\infty}\zeta^{n}\left[\mu\left(\Delta_{n}\cap D\cap t^{-1}\left(A\right)\right)+\mu\left(\Delta_{n}\cap\overline{D}\cap t^{-1}\left(A\right)\right)\right]
\end{align*}
Two distinct cases emerge. When $t^{n}\left(1\right)<1/\beta$ then
one has that $\Delta_{n}\cap\overline{D}=\varnothing$. Thus, the
second term can be written as 
\begin{align*}
\mu\left(\Delta_{n}\cap\overline{D}\cap t^{-1}\left(A\right)\right) & =\Theta\left(t_{n}-\frac{1}{\beta}\right)\mu\left(\Delta_{n}\cap\overline{D}\cap t^{-1}\left(A\right)\right)\\
 & =\Theta\left(t_{n}-\frac{1}{\beta}\right)\frac{1}{\beta}\mu\left(\Delta_{n+1}\cap A\right)
\end{align*}
where the second line follows from the first by linearity, and that
$\overline{D}$ selected out one of the two branches of $t^{-1}\left(A\right)$.
Meanwhile, when $t^{n}\left(1\right)>1/\beta$, then $D\subset\Delta_{n}$
so that $D\cap\Delta_{n}=D$. Thus, the first term splits into two:
\begin{align*}
\Theta\left(t_{n}-\frac{1}{\beta}\right)\mu\left(\Delta_{n}\cap D\cap t^{-1}\left(A\right)\right) & =\Theta\left(t_{n}-\frac{1}{\beta}\right)\mu\left(D\cap t^{-1}\left(A\right)\right)\\
 & =\Theta\left(t_{n}-\frac{1}{\beta}\right)\frac{1}{\beta}\mu\left(A\right)
\end{align*}
while
\[
\Theta\left(\frac{1}{\beta}-t_{n}\right)\mu\left(\Delta_{n}\cap D\cap t^{-1}\left(A\right)\right)=\Theta\left(\frac{1}{\beta}-t_{n}\right)\frac{1}{\beta}\mu\left(\Delta_{n+1}\cap A\right)
\]
Reassembling these pieces and making use of $\Delta_{0}\cap A=A$
one gets
\begin{align*}
\nu\left(t^{-1}\left(A\right)\right) & =\sum_{n=0}^{\infty}\frac{\zeta^{n}}{\beta}\left[\mu\left(A\right)\Theta\left(t_{n}-\frac{1}{\beta}\right)+\mu\left(\Delta_{n+1}\cap A\right)\right]\\
 & =\frac{1}{z}\nu\left(A\right)-\frac{\mu\left(A\right)}{z}+\mu\left(A\right)\sum_{n=0}^{\infty}\frac{\zeta^{n}}{\beta}\Theta\left(t_{n}-\frac{1}{\beta}\right)\\
 & =\frac{1}{z}\nu\left(A\right)+\frac{\mu\left(A\right)}{z}E\left(\beta;z\right)
\end{align*}
with the constant term as before, in eqn \ref{eq:holomorphic-disk}:
\begin{align*}
E\left(\beta;z\right) & =-1+\zeta\sum_{n=0}^{\infty}\zeta^{n}\Theta\left(t_{n}-\frac{1}{\beta}\right)\\
 & =-1+\zeta\sum_{n=0}^{\infty}\zeta^{n}d_{n}\left(\frac{1}{2}\right)
\end{align*}
As before, one has for $z=1$ that $E\left(\beta;1\right)=0$ and
so $\nu\circ t^{-1}=\nu$ is indeed the measure invariant under $t^{-1}$.
Other eigenvalues can be found for those values of $z$ for which
$E\left(\beta;z\right)=0$. The task at hand is then to characterize
$E\left(\beta;z\right)$.

\subsection{Exploring $E\left(\beta;z\right)$}

\begin{figure}
\caption{Disk \label{fig:Disk-of-analytic}of $E\left(\beta;z\right)$}

\includegraphics[width=1\columnwidth]{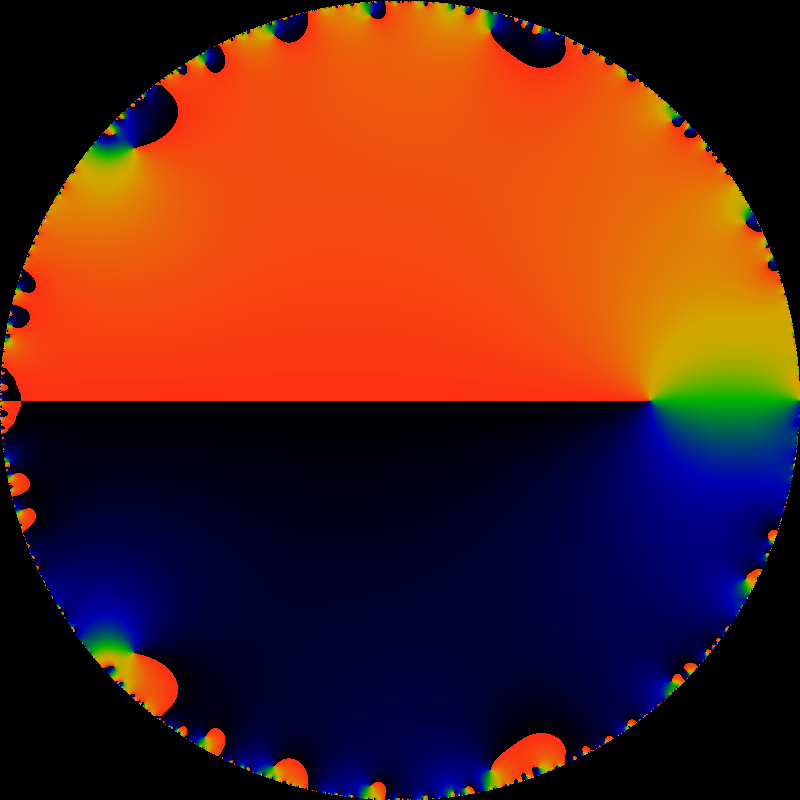}

The above shows a visualization of $E\left(\beta;z\right)$ for $\beta=1.6$
in the complex $\zeta=z/\beta$ plane. The disk consists of all $\zeta$
values with $\left|\zeta\right|\le1$. The plot is a phase plot, showing
the argument $\arg E\left(\beta;z\right)\in\left[-\pi,\pi\right]$.
The color-coding of blue-black puts the phase just above $-\pi$,
green near 0, and red just below $+\pi$. Locations where the phase
wraps around counter-clockwise (right-handed) are zeros of $E\left(\beta;z\right)$;
this follows from Cauchy's principle. The most prominent zero at the
right-hand side of the image corresponds to $z=1$, located at $\zeta=1/\beta$
so well inside the circle. This corresponds to the Gelfond–Parry invariant
measure. Other zeros are seen at the end of black-red whiskers, wrapping
through green. These accumulate at the $\zeta=1$ boundary; in general,
there are a countable number of such zeros. The corresponding eigenvalues
are located at $\lambda=1/z=1/\beta\zeta$, so the accumulation ring
of zeros occurs at $\left|\lambda\right|=1/\beta$. 

\rule[0.5ex]{1\columnwidth}{1pt}
\end{figure}

The function is easily explored numerically. It is clearly convergent
in the unit disk $\left|\zeta\right|<1$ and has no poles in the disk.
For almost all $\beta$, there seem to be a countable number of zeros
within the disk, accumulating uniformly on the boundary as $\left|\zeta\right|\to1$.
An example is shown in figure \ref{fig:Disk-of-analytic}. The notion
of ``uniformly'' will be made slightly more precise in the next
section, where it is observed that, for certain special values of
$\beta$, the bit-sequence $d_{n}\left(\frac{1}{2}\right)$ is periodic,
and thus is a polynomial. When it is polynomial, there are a finite
number of zeros (obviously; the degree of the polynomial), which are
distributed approximately uniformly near the circle $\left|\zeta\right|=1$.
As the degree of the polynomial increases, so do the number of zeros;
but they remain distributed approximately evenly. In this sense, the
limit of infinite degree seems to continue to hold.

A handful of selected zeros are listed in the table below. The numbers
are accurate to about the last decimal place.

\medskip{}

\begin{center}
\begin{tabular}{|c|c|c|c|}
\hline 
$\beta$ & $z$ & $\left|z\right|$ & $1/z$\tabularnewline
\hline 
\hline 
1.8 & -1.591567859 & 1.59156785 & -0.6283112558\tabularnewline
\hline 
1.8 & -1.1962384 +i 1.21602223 & 1.70578321 & -0.41112138 - i 0.41792066\tabularnewline
\hline 
1.8 & 0.99191473 +i 1.44609298 & 1.75359053 & 0.32256553 -i 0.47026194\tabularnewline
\hline 
1.6 & -1.06365138 +i 1.00895989 & 1.46606764 & -0.49487018 -i 0.46942464\tabularnewline
\hline 
1.4 & 0.55083643 +i 1.17817108 & 1.30057982 & 0.32564816 -i 0.69652119\tabularnewline
\hline 
1.2 & 0.95788456 +i 0.60733011 & 1.13419253 & 0.74462841 -i -0.47211874\tabularnewline
\hline 
\end{tabular}
\par\end{center}

\medskip{}

\begin{figure}
\caption{Typical Eigenfunction\label{fig:Typical-Eigenfunction}}

\begin{centering}
\includegraphics[width=0.7\columnwidth]{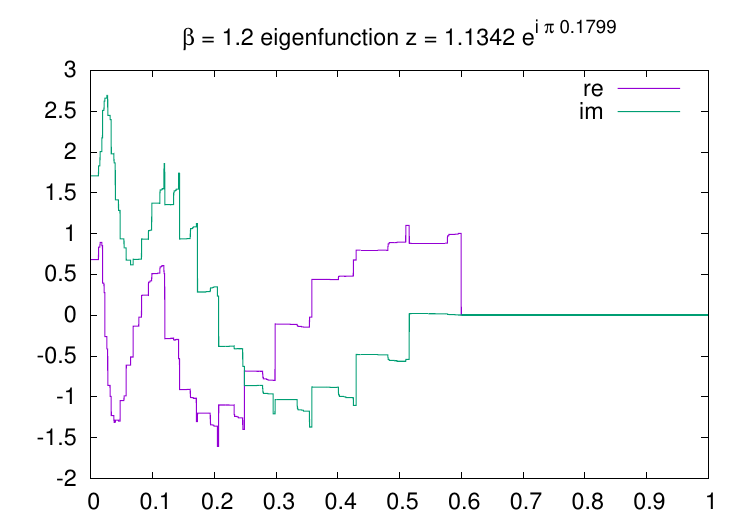}
\par\end{centering}
A typical eigenfunction $v_{\beta;z}\left(x\right)$ solving $\mathcal{L}_{\beta}v=\lambda v$
for $\beta=1.2$. This is the eigenfunction is as defined in eqn \ref{eq:coherent-first-order}.
This eigenfunction corresponds to the zero $z=1.1342\exp i\pi0.1799$
of $E\left(\beta;z\right)$, alternately of $C\left(\beta;z\right)$,
as defined in eqn \ref{eq:holomorphic-disk}. The eigenvalue is $\lambda=1/z$.
Since the eigenvalue is complex, so is the eigenfunction. The real
and imaginary parts are paired in a way that vaguely resembles sine
and cosine; such phased offsets are generic. 

\rule[0.5ex]{1\columnwidth}{1pt}
\end{figure}

These are not particularly meaningful numbers; they just give a flavor
for some locations of eigenvalues. Given a zero, the corresponding
eigenfunction is also very easily computed. A typical eigenfunction
is shown in figure \ref{fig:Typical-Eigenfunction}; this is for the
zero listed in the last row of the table above. Although it is unlike
the figure \ref{fig:Undershift-Density-Distribution}, in that it
is not strictly decreasing, it does have the same general plateau-like
regions. Note that all such eigenfunctions are bounded and generally,
differentiable-nowhere.

As the zeros accumulate onto the circle $\left|\zeta\right|\to1$,
there appears to be no way to holomorphically continue the function
$E\left(\beta;z\right)$ outside of the unit circle. This indicates
that there is a lower bound on the possible eigenvalues that can be
reached via the eigenfunctions of eqn \ref{eq:coherent-first-order}.
The lower bound is $\left|\lambda\right|=1/\beta$; the eigenfunctions
with smaller eigenvalue do not arise from eqn \ref{eq:coherent-first-order}.

The apparent reason for this is that the coherent state was constructed
as a deformed version of the repeated iteration of $\nu_{0}=1$ a
constant, as described in earlier sections. Iterating on $\nu_{0}\left(y\right)$
a polynomial appears to generate eigenfunctions with eigenvalue bounded
by $1/\beta^{p}$ for a polynomial of degree $p$. One explicit example,
generated by a parabola, is given in the next section. A simple, direct
construction remains elusive.

A different derivation of this same holomorphic form will be given
in the next chapter. There, the beta shift is shown to be conjugate
to an infinite shift matrix, in explicit upper-Hessenberg form. Thus,
the shift is explicitly a shift on the monomials $\zeta^{n}$ as the
basis elements.

\subsection{Other eigenfunctions}

Eigenfunctions with eigenvalues $\lambda$ less than $1/\beta$ clearly
exist. Below is an example with $\lambda=1/\beta^{2}$ for $\beta=\varphi=1.618\cdots$
the golden mean. It appears that similar constructions can be performed
for any value of $\beta$ corresponding to a finite length orbit.
Despite this, exhibiting an explicit form remains difficult, verbose
and unsatisfactory. A precise and compact description is desirable,
but remains undiscovered. Turgid details can be found in the research
diary accompanying this text.

\subsubsection{Piecewise quadratic eigenfuntions}

To portray the problem, a piece-wise parabolic presentation is pursued.
For $\beta=\varphi$, the function

\begin{align*}
\rho\left(y\right) & =\begin{cases}
\varphi y^{2}-y+\frac{1}{8} & \mbox{for }0\le y<\frac{1}{2}\\
y^{2}-y+\frac{\varphi}{8} & \mbox{for }\frac{1}{2}\le y<\frac{\varphi}{2}\\
0 & \mbox{for }\frac{\varphi}{2}\le y\le1
\end{cases}
\end{align*}
satisfies $\mathcal{L}_{\varphi}\rho=\lambda\rho$ for $\lambda=\varphi^{-2}$.

\subsubsection{Quasi-resonances}

Returning to the expression for $\mathcal{M}\mu_{n}=\mu_{n+1}$ given
in section \ref{subsec:Scaled-iteration}, one might search for resonances
by hypothesizing $\mu_{0}$ to be some polynomial, and then look for
sequences $\lambda^{n}\mu_{n}\to\mu$ that converge uniformly in the
$n\to\infty$ limit. Performing such a search reveals surprises. The
following was observed.

Fix $\beta=1.6$. Then:
\begin{itemize}
\item If $\mu_{0}\left(x\right)=x-1/2$ then $\mu_{n}$ converges to $\approx-0.0869229\,\nu$
with $\nu$ the Gelfond–Parry invariant measure.
\item If $\mu_{0}\left(x\right)=x-1/2+0.0869229$ (so that the convergent
above is subtracted), then the $\mu_{n}$ appear to bounce around
ergodically, maintaining a bounded norm. Three behaviors are apparent:
$\beta^{n}\int\mu_{n}=\mbox{const}\approx0.2767$. The $L_{1}$ and
$L_{2}$ norms bounce around but remain bounded: $0.3<\beta^{n}\int\left|\mu_{n}\right|<0.8$
and $0.1<\beta^{n}\int\left|\mu_{n}\right|^{2}<1.0$. Each iterate
is approximately orthogonal to the prior one: $\beta^{2n+1}\int\mu_{n}\mu_{n+1}\approx0$.
Here, $\int f=\int_{0}^{1}f\left(x\right)dx$ is just short-hand notation.
\end{itemize}
So there are two questions: what is this magic constant $0.0869229$
? What is this ergodic bouncing? One obvious hypothesis for an answer
to the second question does not bear out. From general principles,
it is reasonable to assume that the bouncing behavior is due to a
complex eigenvalue and eigenfunction. Assume $\mathcal{M}\left(a+ib\right)=\lambda e^{i\phi}\left(a+ib\right)$
for some unknown real functions $a,b$ and unknown real phase angle
$\phi$. Assume $\int ab=0$ (the real and imaginary parts are orthogonal.)
Assume further that $\mu_{n}\to K\left(a\cos\theta+b\sin\theta\right)$
for some unknown normalization $K$ and unknown mixing angle $\theta$.
These assumptions imply that $\int\mu_{n}\mu_{n+1}=K^{2}\lambda\cos\phi$
and so one should observe that $\int\mu_{n}\mu_{n+1}\to\mbox{const}$.
This is not born out by numerics.

If one considers instead a polynomial $\mu_{0}\left(x\right)=x^{n}+c_{n-1}x^{n-1}+\cdots+c_{0}$,
one can find an analogous behavior: $\mu_{n}\to K\nu$. The leading
divergence can be subtracted, so that iterating on $\mu_{0}\left(x\right)=x^{n}+c_{n-1}x^{n-1}+\cdots+c_{0}-K$
gives a sequence of $\mu_{n}$ whose length diminishes uniformly,
in that $\beta^{n}\int\mu_{n}=\mbox{const}$, but otherwise bounce
around ergodically. The relationship between $K$ and the polynomial
is opaque. Of course, it can be computed exactly; the polynomial can
be inserted into eqn \ref{eq:iterated xfer oper}, binomials can be
expanded with binomial coefficients, and one obtains $K$. However,
these formulas are very tangled; an intuitive description of the manifold
or variety is absent.

It has been shown that the beta transform is weak mixing. Is the above
a symptom of this? Is this how mixing manifests itself? If so, then
it seems that thre should be a descriptive framework that makes this
explicit.

\subsection{Conclusion}

It appears that a complete description of the transfer operator remains
elusive. For larger eigenvalues, the analytic approach given above
appears to be sufficient, at least for the case where the transfer
operator acts on a sufficiently tame space of functions. A cohesive
description of the spectrum for $\left|\lambda\right|<1/\beta$ remains
absent. Worse, there appear to be a variety of almost-resonances that
are relatively easy to discover and explore numerically, but defy
coherent explanation. Particularly annoying is the lack of a theoretical
framework that articulates how mixing manifests itself in the structure
of the transfer operator.

\pagebreak{}

\section{Finite Orbits}

The iteration of the midpoint $m_{0}=\beta/2$, that is, the iterated
series $m_{n}=T_{\beta}^{n}\left(\beta/2\right)$ is ergodic in the
unit interval, for almost all values of $\beta$. However, for certain
values of $\beta$, the midpoint iterates will hit the point $x=1/2$
where the $\beta$-shift map has a discontinuity. Here, iteration
stops: at the next step, this point is defined to iterate to zero,
in eqn \ref{eq:downshift}. Zero is a fixed point, and so there is
nowhere further to go. This section explores these special values
of $\beta$.

Aside from the definition in eqn \ref{eq:downshift}, one can consider
the modified map, where the less-than sign has been altered to a less-than-or-equals:
\[
T_{\beta}^{\le}\left(x\right)=\begin{cases}
\beta x & \mbox{ for }0\le x\le\frac{1}{2}\\
\beta\left(x-\frac{1}{2}\right) & \mbox{ for }\frac{1}{2}<x\le1
\end{cases}
\]
In this map, the point $x=1/2$ iterates to $\beta/2$, which is just
the initial midpoint itself. In this case, the halted orbits become
periodic orbits. 

The $\beta$ values at which the midpoint has a periodic or terminating
orbit will be called ``trouble spots'', for lack of a better term.
They can be imagined to be prototypes of a bifurcation point, ``depending
delicately on initial conditions'': where two choices are possible,
depending on infinitesimally small perturbations of $m_{0}$, or,
alternately, of $\beta$.

Before the end of the section, it will be shown that these ``trouble
spots'' are dense in the interval $1\le\beta\le2$, that they can
be placed in a one-to-one bijection with the dyadic rationals, and
that this bijection can be extended to the rationals (describing ultimately-periodic
orbits) and then the reals. The map is monotonically increasing, continuous
but not differentiable; it vaguely resembles the blancmange curve.
Of course, it has a peculiar fractal self-similarity, as all maps
to the infinite binary tree do. 

\subsection{The $\beta$-generalized Golden Ratio}

Trouble spots occur whenever the $p$'th iterate $m_{p}=T_{\beta}^{p}\left(m_{0}\right)$
lands at the discontinuity, so that one may take either $m_{p}=0$
or $m_{p}=m_{0}$. The iteration immediately before corresponds to
$m_{p-1}=1/2$. The length of the orbit is $p$.

The first finite orbit can be found when $\beta=\varphi=\left(1+\sqrt{5}\right)/2$
the Golden Ratio. In this situation, one has that $m_{0}=\varphi/2$
and $m_{1}=1/2$. The length of the orbit is $\nu=2$. For $\nu=3$,
there are two such trouble spots, which occur when either $\beta^{3}-\beta^{2}-1=0$
or when $\beta^{3}-\beta^{2}-\beta-1=0$. These correspond to the
values of $\beta=1.465571231876768\cdots$ and $\beta=1.839286755214161\cdots$.

\begin{figure}
\caption{Location of Midpoints\label{fig:Location-of-Midpoints}}

\begin{centering}
\includegraphics[width=1\columnwidth]{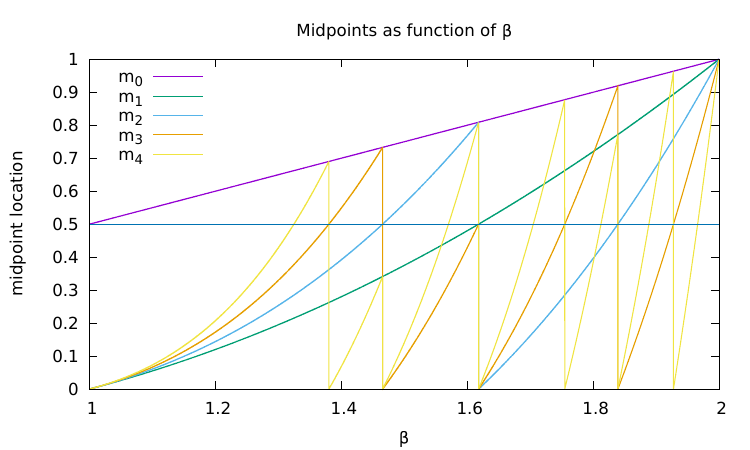}
\par\end{centering}
This chart illustrates the location of the first five midpoints, $m_{0},m_{1},\cdots,m_{4}$
as a function of $\beta$. When $m_{k}=0.5$, further iteration is
ambiguous, as this is the location of the discontinuity in the shift
map, and the next iteration leaves the midpoint bifurcated. These
are the ``trouble spots''. The first trouble spot is visible for
$m_{1}=1/2$, corresponding to $\beta=\varphi$ and a length of $\nu=2$.
The midpoint $m_{2}$ crosses $1/2$ (in the ascending direction)
at $\beta=1.465\cdots$ and $\beta=1.839\cdots$, corresponding to
orbits of length $\nu=3$. It also crosses discontinuously downwards
at $\beta=\varphi$. This crossing point has already been assigned
to a shorter orbit. The midpoint $m_{3}$ has three new crossings.
It also rises to touch $1/2$ at $\beta=\varphi$; but this has already
been assigned to shorter orbits. The midpoint $m_{4}$ has six new
crossings. Crossings will generally fall to the left and right of
earlier crossings, and so are related in a bracketing relationship.
The bracketing is not the full binary tree; it is pruned, as shorter
orbit assignments can knock out longer one. For example, $m_{4}$
falls back down at the $m_{2}$ crossing at $\beta=1.465\cdots$,
and so, here, $m_{4}$ never even gets close to $1/2$; it won't bracket
until later. This is a formalized version of figure \ref{fig:Midpoint-Trace},
which shows midpoints of all orders.

\rule[0.5ex]{1\columnwidth}{1pt}
\end{figure}

Where else are such spots located? Consider, for example, $m_{4}=T_{\beta}^{4}\left(m_{0}\right)$,
and consider the movement of $m_{4}$ as $\beta$ is swept through
the range $1<\beta<2$. This is shown in figure \ref{fig:Location-of-Midpoints}.
As made clear in the image, three new degenerate points appear. These
are located at $\beta=1.380327757\cdots$ and $\beta=1.754877666\cdots$
and $\beta=1.927561975\cdots$, which are the real roots of $\beta^{4}-\beta^{3}-1=0$
and $\beta^{4}-\beta^{3}-\beta^{2}-1=0$ and $\beta^{4}-\beta^{3}-\beta^{2}-\beta-1=0$
respectively.

Following a suggestion by Dajani\cite{Dajani97}, the $\beta$ numbers
corresponding to the trouble spots may be called ``generalized golden
means''. Unfortunately, the term ``generalized golden mean'' is
in common use, and is applied to a variety of different systems. Not
all are relevant to the present situation; one that is, is given by
Hare \emph{et al.}\cite{Hare14} who provide series expansions for
the real roots of $\beta^{p}-\sum_{k=0}^{n-1}\beta^{k}=0$; these
are known as the n-bonacci constants (Fibonacci, tribonacci, tetranacci,
\emph{etc.}). Stakhov\cite{Stakhov05} considers $\beta^{p+1}-\beta^{p}-1=0$
in general settings. Some, but not all of these numbers are known
to be Pisot numbers or Salem numbers\cite{Maia07}. In what follows,
these will be referred to as the ``beta golden means'', since all
of the ones that appear here have explicit origins with the beta shift.

\subsection{Counting Orbits}

How many trouble spots are there? The table below shows the count
$M_{\nu}$ of the number of ``new'' trouble spots, as a function
of the orbit length $\nu$. 
\begin{center}
\begin{tabular}{|c||c|c|c|c|c|c|c|c|c|c|c|}
\hline 
$\nu$ & 2 & 3 & 4 & 5 & 6 & 7 & 8 & 9 & 10 & 11 & 12\tabularnewline
\hline 
$M_{\nu}$ & 1 & 2 & 3 & 6 & 9 & 18 & 30 & 56 & 99 & 186 & 335\tabularnewline
\hline 
\end{tabular}
\par\end{center}

This appears to be Sloane's OEIS A001037 which has a number of known
relationships to roots of unity, Lyndon words, and the number of orbits
in the tent map. The values are given by Moreau's necklace-counting
function. The trouble spots are the positive real roots of polynomials
of the form 
\[
p_{\left\{ b_{k}\right\} }\left(\beta\right)=\beta^{\nu}-\beta^{\nu-1}-b_{1}\beta^{\nu-2}-b_{2}\beta^{\nu-3}-\cdots-1=0
\]
with the $\left\{ b_{k}\right\} $ being certain binary bit sequences.
There is just one such (positive, real) root for each such polynomial.
These polynomials are relatively prime, in the sense that a bit-sequence
$b_{k}$ is disallowed if it has the same root as some lower-order
polynomial. For example, $\beta^{4}-\beta^{3}-\beta-1$ is disallowed;
it has the same root as $\beta^{2}-\beta-1$. Equivalently, the quadratic
is a factor of the quartic; the quartic is not relatively prime with
respect to the quadratic.

The reason for the appearance of the necklace-counting function is
straightforward: it is counting the number of distinct orbits of a
given length. An orbit of length $\nu$ is, by definition, a point
$x$ such that $x=T_{\beta}^{\nu}\left(x\right)$. Such an orbit generates
a binary string, of length $\nu$ corresponding to whether $T_{\beta}^{j}\left(x\right)<1/2$
is true or not. A cyclic permutation of such a string still corresponds
to the same orbit; a reversed permutation does not: thus, it is a
necklace without reversal. The necklace-counting function gives the
number of distinct, unique orbits of a given length that cannot be
factored into shorter orbits. The beta transformation has the property
that every possible orbit does occur; none are prohibited.

Yet, the bit-string defining the polynomial is not a necklace; it
cannot be rotated. Each bit-string corresponds to a unique polynomial,
having roots that differ from those of other polynomials. The polynomials
also have a canonical order, fixed as the integer that generates the
bit-string; they cannot be reordered. Ideas such as Lyndon words apply
to the orbits, but not to the defining polynomials. The ordering of
the polynomials is \emph{not} the lexicographic ordering of the Lyndon
words, and cannot be brought into this order.

The values of $M_{n}$ are given explicitly by Moreau's necklace-counting
function
\[
M_{n}=\frac{1}{n}\sum_{d|n}2^{d}\mu\left(\frac{n}{d}\right)
\]
where the sum runs over all integers $d$ that divide $n$ and $\mu$
is the Möbius function. The generating function is
\[
\frac{t}{\frac{1}{2}-t}=\sum_{n=1}^{\infty}nM_{n}\frac{t^{n}}{1-t^{n}}
\]
which has a radius of convergence of $\left|t\right|<1/2$. For large
$n$, the asymptotic behavior can be trivially deduced from the defining
sum:
\[
M_{n}=\frac{2^{n}}{n}-\mathcal{O}\left(\frac{2^{n/2}}{n}\right)
\]

The above counting function is for necklaces with only two colors.
In general, one can have necklaces with 3 or more colors; can that
happen here? Yes, of course: if one considers the general $\beta$-transform
for $2<\beta$, then, in general, it can be taken as a ``kneading
transform'' with $\left\lceil \beta\right\rceil $ branches or folds
in it. The analogous trouble-spots again appear, and they can appear
after an arbitrary finite-length orbit. Insofar as they correspond
to periodic orbits, they are necessarily counted by the necklace-counting
function. That is, one must consider all possible strings of $\left\lceil \beta\right\rceil $
letters, modulo a cyclic permutation: this is the very definition
of a necklace (or ``circular word''). The number of such necklaces
is given by the necklace-counting function. Each such orbit is necessarily
represented by a Lyndon word, which is a representative of the conjugacy
class of the orbit.

The isomorphism of different systems described by necklace polynomials
is a subject that gets some fair amount of attention. Golomb gives
an isomorphism between the irreducible polynomials over $\mathbb{F}_{p}$,
for $p$ prime and necklaces built from Lyndon words.\cite{Golomb69,Reutenauer93}
A number of other results exist, including \cite{Duzhin13,Cattell00}.
At any rate, a closer study of the beta-polynomials seems to be called
for.

\subsection{$\beta$-Golden Polynomials}

The ``trouble spots'' occur whenever the $k$'th iterate $m_{k}=T_{\beta}^{k}\left(m_{0}\right)$
of the midpoint $m_{0}=\beta/2$ lands on the starting midpoint $m_{k}=m_{0}$;
alternately, when $m_{k-1}=1/2$. Because of the piece-wise linear
form of $T_{\beta}$, the $k$'th iterate will be a piece-wise collection
of polynomials, each of order $k$, each of the form $p_{\left\{ b_{k}\right\} }\left(\beta\right)$.
These must be arranged such that $p_{\left\{ b_{k}\right\} }\left(\beta\right)=0$
at each discontinuity, as illustrated in figure \ref{fig:Location-of-Midpoints}.
This constrains the polynomials that can appear; it constrains the
possible coefficients $\left\{ b_{k}\right\} $; not all bit-sequences
appear. The sequences that do appear encode the orbit of the mid-point;
see below.

The table below explicitly shows the polynomials for the first few
orders. A polynomial is included in the table if it is an iterate
of a previous polynomial, and if it's real root is bracketed by the
roots of the earlier iterates. Adopting ordinal numbering, $\ensuremath{p_{n}\left(\beta\right)}$
must have the form
\begin{equation}
p_{n}\left(\beta\right)=\begin{cases}
\beta\left(p_{n/2}\left(\beta\right)+1\right)-1 & \mbox{ for \ensuremath{n} even }\\
\beta p_{\left(n-1\right)/2}\left(\beta\right)-1 & \mbox{ for \ensuremath{n} odd }
\end{cases}\label{eq:polynomial recursion}
\end{equation}
This recursion terminates at $p_{0}\left(\beta\right)=\beta-1$.

The positive real root $r_{n}$ satisfying $p_{n}\left(r_{n}\right)=0$
is unique; the other $n-1$ roots are complex; they are arranged in
a roughly evenly-spaced ring on the complex plane, not far from the
unit circle, reminiscent of roots of unity. There is always a positive
real root, which satisfies $1\le r_{n}<2$; the real roots and the
polynomials are in one-to-one correspondence. The roots must be bracketed
(to the left and right) by the roots of polynomials occurring earlier
in the sequence; if the root is not bracketed, then the corresponding
polynomial does not appear in the list.

The bracketing constraint can be represented by a recursive function
$\theta_{n}\left(\rho\right)$ returning a boolean true/false value,
as to whether a given polynomial is acceptable. It is
\begin{equation}
\theta_{n}\left(\rho\right)=\begin{cases}
\Theta\left(r_{n/2}-\rho\right)\cdot\theta_{n/2}\left(\rho\right) & \mbox{ for \ensuremath{n} even }\\
\theta_{\left\lfloor n/2\right\rfloor }\left(\rho\right) & \mbox{ for \ensuremath{n} odd }
\end{cases}\label{eq:theta mask}
\end{equation}
The Heaviside $\Theta\left(x\right)$ used here is one for strictly
positive $x>0$, and is zero otherwise. This is important, as using
$x\ge0$ will not work. In numerical work, the test should be bounded
away from zero. The recursion compares the candidate $\rho$ to some
root at each lower order. The recursion terminates at $\theta_{0}\left(\rho\right)=1$.
The index $n$ corresponds to a valid polynomial, and thus a valid
root, if and only if $\theta_{n}\left(r_{n}\right)=1$. Effectively,
this states that roots of higher-order polynomials must be less than
a certain sequence of lower-order roots. This is visible in the location
of the discontinuities in figure \ref{fig:Location-of-Midpoints}:
new discontinuities at higher orders must occur to the left of earlier
ones.

For example, the polynomial $\beta^{3}-\beta-1$ is excluded from
the list simply because it is not an iterate of an earlier polynomial,
even though it has the interesting real root $1.324717957244746\cdots$,
the ``silver constant''. The numbering scheme does not even have
a way of numbering this particular polynomial. Despite this, the silver
constant does appear, but a bit later, as the root of $p_{8}=\beta^{5}-\beta^{4}-1$,
which is an allowed polynomial.

The polynomial $p_{5}=\beta^{4}-\beta^{3}-\beta-1$ is excluded because
it has $\varphi=1.618\cdots$ as a root, which was previously observed
by $p_{1}$. The polynomial $p_{9}=\beta^{5}-\beta^{4}-\beta-1$ is
excluded because it's root, $r_{9}=1.497094048762796\cdots$ is greater
than its predecessor $r_{2}$; the recursive algorithm does not compare
it to $r_{4}$. Note that $p_{9}$ is relatively prime to the earlier
polynomials, so irreducibility is not a sufficient criterion; the
root must also be less.

The indexing has the property that, whenever $\theta_{n}\left(r_{n}\right)=1$,
the integer $2n+1$, expressed as a binary bitstring, encodes both
the coefficients of the polynomial, and also the orbit of the midpoint.
This can be taken as an alternate, non-recursive definition of $\theta_{n}$:
it is one if and only if the orbit of $r_{n}$ encodes the bitsequence
of $2n+1$.

The degree $\nu$ of the polynomial is identical to the length $\nu$
of the orbit; it is $\nu=\left\lceil \log_{2}\left(2n+1\right)\right\rceil $.
The bits $b_{i}$ of the bitstring $2n+1=b_{0}b_{1}b_{2}\cdots b_{\nu}$
correspond to the orbit as 
\[
b_{i}=\Theta\left(T_{\beta}^{i}\left(\frac{\beta}{2}\right)-\frac{1}{2}\right)=d_{i}\left(\frac{1}{2}\right)=k_{i}\left(\frac{\beta}{2}\right)
\]
where the $d_{i}$ are as given before, in eqn \ref{eq:rescaled-bitseq},
and the $k_{i}$ as in eqn \ref{eq:down-bits}. Note that $b_{0}=1$
always corresponds to $1/2<\beta/2$, always. By convention, the last
digit is always 1, also.

\medskip{}

\begin{center}
\begin{tabular}{|c|c|c|c|c|}
\hline 
order $\nu$ & $p_{n}\left(\beta\right)$ & $n$ & binary & root $r_{n}$\tabularnewline
\hline 
0 & $1$ &  &  & \tabularnewline
\hline 
\multirow{2}{*}{1} & $\beta$ &  & 0 & 0\tabularnewline
\cline{2-5} \cline{3-5} \cline{4-5} \cline{5-5} 
 & $\beta-1$ & 0 & 1 & 1\tabularnewline
\hline 
2 & $\beta^{2}-\beta-1$ & 1 & 11 & $\varphi=\frac{1+\sqrt{5}}{2}=1.618\cdots$\tabularnewline
\hline 
\multirow{2}{*}{3} & $\beta^{3}-\beta^{2}-1$ & 2 & 101 & $1.465571231876768\cdots$\tabularnewline
\cline{2-5} \cline{3-5} \cline{4-5} \cline{5-5} 
 & $\beta^{3}-\beta^{2}-\beta-1$ & 3 & 111 & $1.839286755214161\cdots$\tabularnewline
\hline 
\multirow{3}{*}{4} & $\beta^{4}-\beta^{3}-1$ & 4 & 1001 & $1.380277569097613\cdots$\tabularnewline
\cline{2-5} \cline{3-5} \cline{4-5} \cline{5-5} 
 & $\beta^{4}-\beta^{3}-\beta^{2}-1$ & 6 & 1101 & $1.754877666246692\cdots$\tabularnewline
\cline{2-5} \cline{3-5} \cline{4-5} \cline{5-5} 
 & $\beta^{4}-\beta^{3}-\beta^{2}-\beta-1$ & 7 & 1111 & $1.927561975482925\cdots$\tabularnewline
\hline 
\multirow{6}{*}{5} & $\beta^{5}-\beta^{4}-1$ & 8 & 10001 & $1.324717957244746\cdots$\tabularnewline
\cline{2-5} \cline{3-5} \cline{4-5} \cline{5-5} 
 & $\beta^{5}-\beta^{4}-\beta^{2}-1$ & 10 & 10101 & $1.570147312196054\cdots$\tabularnewline
\cline{2-5} \cline{3-5} \cline{4-5} \cline{5-5} 
 & $\beta^{5}-\beta^{4}-\beta^{3}-1$ & 12 & 11001 & $1.704902776041646\cdots$\tabularnewline
\cline{2-5} \cline{3-5} \cline{4-5} \cline{5-5} 
 & $\beta^{5}-\beta^{4}-\beta^{3}-\beta-1$ & 13 & 11011 & $1.812403619268042\cdots$\tabularnewline
\cline{2-5} \cline{3-5} \cline{4-5} \cline{5-5} 
 & $\beta^{5}-\beta^{4}-\beta^{3}-\beta^{2}-1$ & 14 & 11101 & $1.888518845484414\cdots$\tabularnewline
\cline{2-5} \cline{3-5} \cline{4-5} \cline{5-5} 
 & $\beta^{5}-\beta^{4}-\beta^{3}-\beta^{2}-\beta-1$ & 15 & 11111 & $1.965948236645485\cdots$\tabularnewline
\hline 
\end{tabular}
\par\end{center}

\medskip{}

The next table lists the acceptable polynomial indexes for order 5,
6 and 7. Again, the coefficients appearing in the polynomial are encoded
by the binary value of $2n+1$ in the sequence. This sequence is not
currently known to OEIS.

\medskip{}

\begin{center}
\begin{tabular}{|c|c|}
\hline 
order $\nu$ & valid indexes\tabularnewline
\hline 
\hline 
5 & 8,10,12,13,14,15\tabularnewline
\hline 
6 & 16,20,24,25,26,28,29,30,31\tabularnewline
\hline 
7 & 32,36,40,42,48,49,50,52,53,54,56,57,58,59,60,61,62,63\tabularnewline
\hline 
\end{tabular}
\par\end{center}

\medskip{}

The properties of this sequence are briefly reviewed in the next section.

\subsection{Properties of the theta sequence}

\begin{figure}
\caption{Indicator function $\theta$\label{fig:Theta-Indicator}}

\includegraphics[width=1\columnwidth]{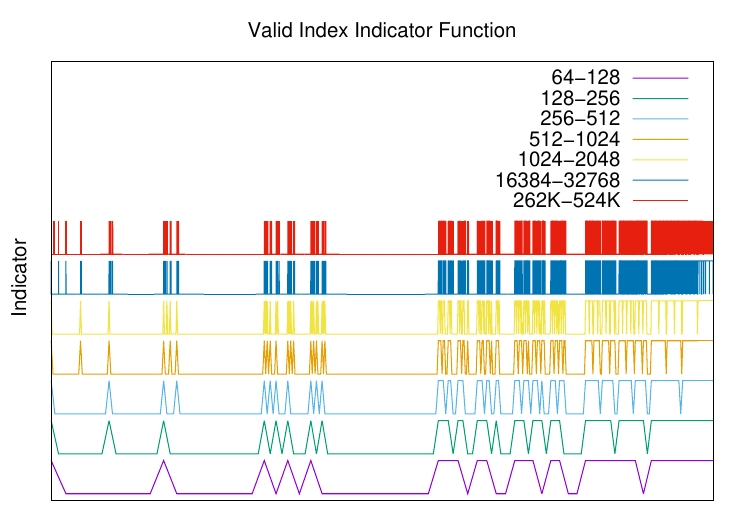}

This figure shows the indicator function $\theta_{n}$, as defined
in eqn \ref{eq:valid index set}, over the range of $2^{6}\le n\le2^{19}$.
The function is approximately periodic as a function of the order
$\nu=\left\lfloor \log_{2}n\right\rfloor +2$ of the corresponding
polynomial; the order is just the length of a finite orbit. This figure
just stacks seven such ranges, for orders $\nu$=6,7,8,9,10,14,18.
One can view each row as a comb. The teeth of the comb are the members
of $\Psi_{\nu}$. The number of teeth at each rank is given by Moreau's
necklace-counting function $M_{\nu}$. The width of the teeth go as
$2^{-\nu}$, but the number of teeth goes as $M_{\nu}\sim2^{\nu}/\nu$.
In the limit $\nu\to\infty$, the teeth converge onto a set of measure
zero. This can be understood as a representation of the set $\Psi$,
taken with a real-valued index, rather than an integer index.\rule[0.5ex]{1\columnwidth}{1pt}
\end{figure}

Define the ``validity set'' of valid or acceptable indexes as 
\begin{equation}
\Psi=\left\{ n\in\mathbb{N}:\theta_{n}\left(r_{n}\right)=1,p_{n}\left(r_{n}\right)=0\right\} \label{eq:valid index set}
\end{equation}
This is just the list of indexes from the previous two tables; it
is $\Psi$=\{0,1, 2,3, 4,6,7, 8,10,12,13,14,15, 16,20,24,25,26,28,29,30,31,
...\}. The $\theta_{n}$ is the acceptance function from eqn \ref{eq:theta mask}.
By abuse of notation, write $\theta_{n}=\mathds{1}_{\Psi}\left(n\right)=\theta\left(n\right)$
for the membership indicator function for this set. This function
is visualized in figure \ref{fig:Theta-Indicator}. 

The elements of the validity set $\Psi$ are ordered; this is the
validity sequence $\psi_{m}=\psi\left(m\right)$. It is convenient
to start the sequence at $\psi_{0}=0$. This corresponds to $\beta=2$
at the right; while, at the far left, for $\beta=1$, write $\psi_{-1}=\infty$.
The function $\psi_{m}$ encodes the locations of one-bits in the
bitmask $\theta$, and so $\theta\left(\psi_{m}\right)=1$. 

The summatory function of the indicator is $S\left(k\right)=\sum_{n=1}^{k}\theta_{n}$.
It counts the total number of one-bits below the location $k$. The
validity sequence is the pullback of the summatory function. Each
$\psi_{m}$ is the smallest integer $k$ for which $m=S\left(k\right)$
holds true, and so one has $m=S\left(\psi_{m}\right)$. The pullback
can be expressed as $\psi\left(m\right)=\psi\left(S\left(\psi\left(m\right)\right)\right)$.

A few additional properties may be noted:
\begin{itemize}
\item For all $m$, $\theta\left(2^{m}\right)=\theta\left(2^{m}-1\right)=1$.
\item If $\theta\left(m\right)=1$ then $\theta\left(2m\right)=1$. By recursion,
$\theta\left(2^{n}m\right)=1$ for all $n$.
\item If $m$ is odd, and if $\theta\left(m\right)=1$, then $\theta\left(\left(m-1\right)/2\right)=1$.
This is reminiscent of the Collatz conjecture.
\item If $\theta\left(m\right)=0$ then $\theta\left(2^{n}\left(2m+1\right)\right)=0$
for all $n$.
\end{itemize}
Each of these properties is visible in figure \ref{fig:Theta-Indicator}.
The second bullet accounts for the stability of the comb-teeth, once
they appear, while the last bullet accounts for the large spaces that
open up, and stay open, never filling in.

It is convenient to partition $\Psi$ into ranks $\nu$ that correspond
to the length of the orbits, or equivalently, the order of the defining
polynomial. Examining the earlier tables, the rank of $n$ is $\nu\left(n\right)=\left\lfloor \log_{2}n\right\rfloor +2$.
The partition is then $\Psi=\bigcup_{\nu=1}^{\infty}\Psi_{\nu}$ with
$\Psi_{\nu}=\left\{ n\in\Psi:\nu\left(n\right)=\nu\right\} $. As
before, it is convenient to extend the partition so that it can deal
with the endpoints $\beta=1$ and $\beta=2$; this is a kind of (two-point)
compactification of these and other various sequences and sets. The
compactification here is to write $\nu\left(0\right)=1$ and $\nu\left(-1\right)=\infty$,
which allows components $\Psi_{1}=\left\{ 0\right\} $ and $\Psi_{\infty}=\left\{ -1\right\} $.
The size of each component is $\left|\Psi_{\nu}\right|=M_{\nu}$ given
by Moreau's necklace counting function $M_{\nu}$.

The representation of $\theta$ as a real number is $\theta=\sum_{n=1}^{\infty}\theta_{n}2^{-n}=0.93258880035365\cdots$.
At this time, this does not appear in OEIS.

\subsubsection{Leaders}

An important subsequence consists of the leaders of the doubling sequences.
These can be defined in several equivalent ways. One property of the
bitmask is that if $\theta\left(m\right)=1$, then $\theta\left(2^{n}m\right)=1$
for all $n$. A leader $\lambda$ is the smallest such $m$ at the
front of such a doubling sequence: it is either an odd number $\lambda=\left(2k+1\right)$
satisfying $\theta\left(\lambda\right)=1$ or it is an even number
of the form $\lambda=2^{h}\left(2k+1\right),h>0$ such that $\theta\left(2^{h}\left(2k+1\right)\right)=1$
but $\theta\left(2^{h-1}\left(2k+1\right)\right)=0$. The minimal
power $h$ defining the leader will be called the ``height of the
leader''.

Equivalently, the leaders can be defined as a subsequence of $\psi_{m}$.
By definition, one has that $\theta\left(\psi_{m}\right)=1$. Thus,
if $\psi_{m}$ is odd, then $\psi_{m}$ is a leader. If $\psi_{m}$
is even, then it is a leader if and only if $\theta\left(\psi_{m}/2\right)=0$.

Given a valid index $\psi_{m}$, it is useful to define a function
$\Lambda$ that provides the map 
\begin{equation}
\Lambda\left(\psi_{m}\right)=2^{h}\left(2\psi_{m}+1\right)\label{eq:Leadership function}
\end{equation}
This map is not monotonic: $\Lambda\left(1\right)=3$ and $\Lambda\left(2\right)=10$
but $\Lambda\left(3\right)=7$. It is, however, one-to-one, as leaders
are always a product of an odd number times some power of two.

Sorting the leaders into ascending order, the start of the sequence
is 1,3,7,10,13,15,25,29,31,... This sequence is not currently known
to OEIS. It is handy to count from one, and to assign $\lambda_{0}=0$,
so that $\lambda_{1}=1$ and $\lambda_{2}=3$ and so on.

\subsection{Location of $\beta$-Golden Roots}

\begin{figure}
\caption{Normalization Integral\label{fig:Normalization-Integral}}

\begin{centering}
\includegraphics[width=1\columnwidth]{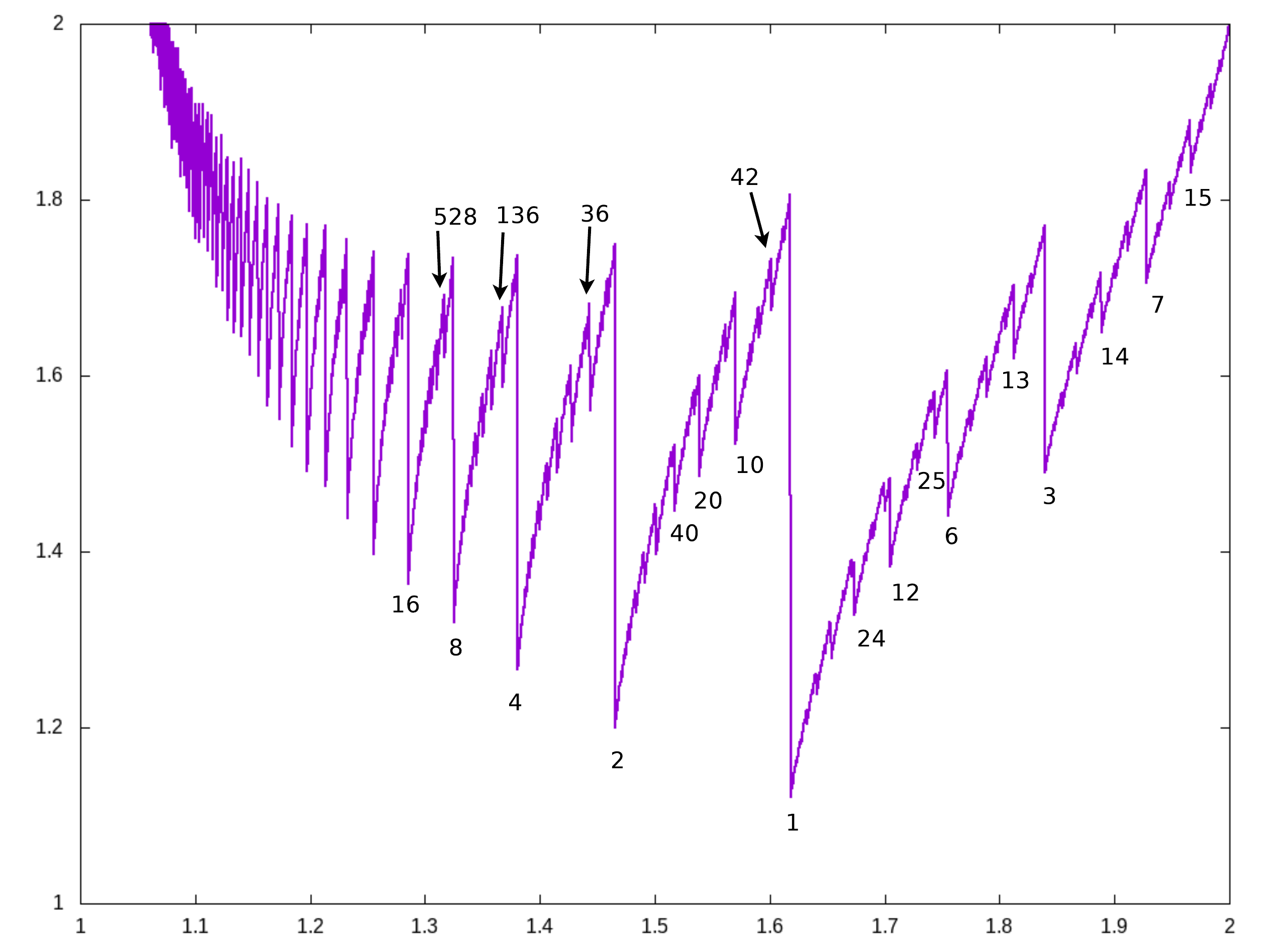}
\par\end{centering}
This figure shows the integral $I\left(\beta\right)=\sum_{n=0}^{\infty}\beta^{-n}T^{n}\left(\frac{\beta}{2}\right)$
with $1<\beta\le2$ running along the horizontal axis. Each discontinuity
corresponds to the location of a real root of one of the $\beta$-Golden
polynomials. Some of these are manually labeled by integers, corresponding
to the polynomial labels from the previous polynomial table. 

\rule[0.5ex]{1\columnwidth}{1pt}
\end{figure}

The location of the roots can be visualized by using the normalization
of the Parry–Gelfond measure. The function in eqn \ref{eq:invariant measure}
or more generally \ref{eq:generalized-eigenfunction} can be integrated
in a straightforward manner. One has
\[
I\left(\beta;z\right)=\int_{0}^{1}\nu_{\beta;z}\left(x\right)dx=\sum_{n=0}^{\infty}\frac{z^{n}}{\beta^{n}}\int_{0}^{1}d_{n}\left(x\right)dx=\sum_{n=0}^{\infty}\frac{z^{n}}{\beta^{n}}T^{n}\left(\frac{\beta}{2}\right)
\]
The result is a sawtooth, shown in figure \ref{fig:Normalization-Integral}.
Each discontinuity corresponds to the real root of one of the polynomials.
The first few are labeled by the integer labels from the previous
table. The doubling sequences and their leaders are easy to identify.

\subsubsection{Bracketing intervals}

That this figure is a self-similar fractal is presumably self-evident.
Thus, for example, the graph to to right of $r_{1}=\varphi=1.618...$
repeats again between $r_{2}$ and $r_{1}$ and again between $r_{4}$
and $r_{2}$. Each such bracketed interval contains a unique largest
discontinuity; it can be seen as being at the front of a doubling
sequence. Each discontinuity is in one-to-one correspondence with
a bracketing interval; the bracketing intervals are all self-similar
to one-another.

A distinct notation for bracketed intervals is useful. Write $\ell\Mapsto f\Mapsfrom\rho$
for the discontinuity $f$ bracketed on the left and right by $\ell$
and $\rho$. By ``left'' and ``right'', it is literally meant
that the three roots are in order, with $r_{\ell}<r_{f}<r_{\rho}$
with the inequalities being strict. Not all ascending sequences of
three roots form a valid bracket; valid brackets are obtained by recursive
subdivision; this is given in the next section. But first, some examples.

Taking the liberty to write $r_{0}=2$, the interval to the right
of $r_{1}$ is then $1\Mapsto3\Mapsfrom0$. The most prominent self-similar
intervals are then $2\Mapsto10\Mapsfrom1$ and $4\Mapsto36\Mapsfrom2$
and $8\Mapsto136\Mapsfrom4$. In each of these examples, the front
$f$ was also a leader, with leadership as defined in the previous
section. This is not always the case: the brackets $1\Mapsto6\Mapsfrom3$
and $1\Mapsto12\Mapsfrom6$ and $1\Mapsto24\Mapsfrom12$ are clearly
visible; they are a part of an index-doubling sequence.

The extreme left side can be assigned the index of $\infty$ so that
$r_{\infty}=1$. Thus, the entire interval $1\le\beta\le2$ corresponds
to the bracket $\infty\Mapsto1\Mapsfrom0$.

\subsubsection{The bracket tree}

The brackets can be arranged into a binary tree, recursively defined.
Any valid interval $\ell\Mapsto f\Mapsfrom\rho$ can be split into
two: the left side and the right side. These two pieces are $\ell\Mapsto2f\Mapsfrom f$
on the left, and $f\Mapsto\Lambda\left(f\right)\Mapsfrom\rho$ on
the right, where $\Lambda$ is the leader function given by eqn \ref{eq:Leadership function}
in the previous section. Recursion starts with $\infty\Mapsto1\Mapsfrom0$.
A bracketing relationship is valid if and only if it appears in this
recursive binary tree.

The left and right moves $L,\mathfrak{R}$ on the binary tree can
then be written as 
\begin{align}
L & :\left(\ell\Mapsto f\Mapsfrom\rho\right)\mapsto\left(\ell\Mapsto2f\Mapsfrom f\right)\label{eq:Bracket recursion}\\
\mathfrak{R} & :\left(\ell\Mapsto f\Mapsfrom\rho\right)\mapsto\left(f\Mapsto\Lambda\left(f\right)\Mapsfrom\rho\right)\nonumber 
\end{align}
The right-move is denoted with a fraktur $\mathfrak{R}$ instead of
a roman $R$ in order to distinguish between the action of the leader
function, and the conventional dyadic map 
\begin{align*}
L & :m\mapsto2m\\
R & :m\mapsto2m+1
\end{align*}
The issue is that not all dyadic $R$ moves result in a valid index:
having $\theta\left(m\right)=1$ does not generally imply that $\theta\left(2m+1\right)$
is one. However, the leader function does provide a successor that
is always valid: $\theta\left(\Lambda\left(m\right)\right)=1$ is
guaranteed by construction.

Any sequence of $L,\mathfrak{R}$ moves is guaranteed to produce a
valid interval; every location in the binary tree is mapped to a valid
index by the bracket recursion moves. A general location on the bracket
tree is

\[
G=L^{a_{1}}\mathfrak{R}^{a_{2}}L^{a_{3}}\cdots\mathfrak{R}^{a_{k}}
\]
and this has the property that if $\theta\left(m\right)=1$ then $\theta\left(G\left(m\right)\right)=1$.
Equivalently, if $m\in\Psi$ then $G\left(m\right)\in\Psi$. This
implies that the bracket recursion relations generate all of $\Psi$.
Starting at $m=1$, all left-right moves produce elements of $\Psi$;
conversely, every element of $\Psi$ can be expressed as some sequence
of left-right moves applied to $m=1$. It is convenient to take $G$
as the function that produces only the ``good'' indexes. Thus
\[
G:L^{a_{1}}R^{a_{2}}L^{a_{3}}\cdots R^{a_{k}}\mapsto L^{a_{1}}\mathfrak{R}^{a_{2}}L^{a_{3}}\cdots\mathfrak{R}^{a_{k}}
\]
Labeling the nodes of the binary tree with the natural numbers $\mathbb{N}$
provides a bijection $G:\mathbb{N}\to\Psi$.

\begin{figure}
\caption{Beta Bracket Map\label{fig:Beta-Bracket-Map}}

\begin{centering}
\includegraphics[width=0.9\columnwidth]{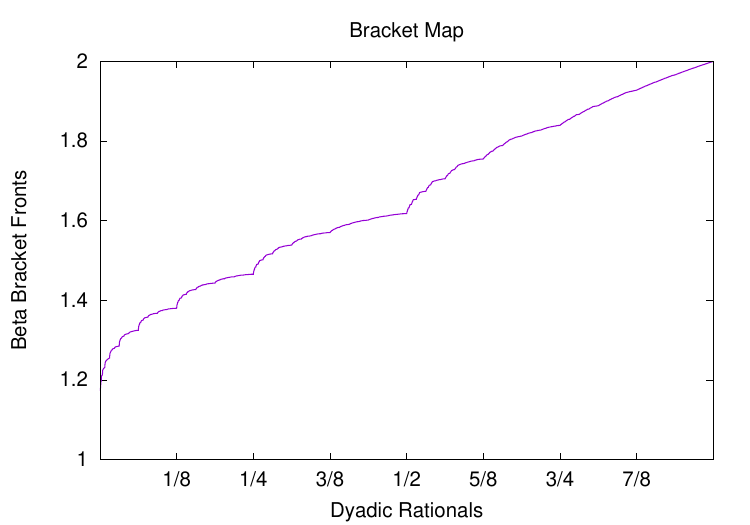}
\par\end{centering}
This figure shows the mapping from the dyadic rationals, interpreted
as strings of left-right moves, and the corresponding $\beta=r_{f}$
root for the front (center) of the corresponding bracket. Thus, 1/2
maps to $r_{1}=\varphi=1.61803...$ and 1/4 maps to $r_{2}=1.46557...$
while 3/4 maps to $r_{3}=1.83929...$ It would appear that the curve
is continuous. The curve does drop to $\beta=1$ at the left, but
it does so very slowly. The brackets $\ell\Mapsto f\Mapsfrom\rho$
remain quite wide, as the limit approaches $\beta=1$. It does get
there, though. An explicit derivation of this limit is given in a
later section. The sharp takeoff at the left is reminiscent of the
Minkowski question-mark mapping. The left generator for the bracket
map is is quite close to left generator $x/\left(1+x\right)$ of the
Minkowski fractal. The Minkowski fractal has $1/\left(2-x\right)$
as the right generator. Replacing this with $x+1/2$ generates a curve
resembling the bracket curve above, having Minkowski-like behavior
to the left, an linear-like behavior to the right. The de Rham curve
construction may be used to generate such curves, given arbitrary
L,R maps. The L,R maps for the bracket curve are both fractal themselves,
so finding affine generators would be a surprise.

\rule[0.5ex]{1\columnwidth}{1pt}
\end{figure}

The bracket recursion relations split the interval $r_{\ell}<r_{f}<r_{\rho}$
into left and right pieces, as well. That is, each node of the binary
tree is labeled by some (unique) $r_{f}$, with the root $r_{\ell}$
and $r_{\rho}$ appearing as predecessors in the tree. Of course,
the binary tree can also be mapped to the dyadic rationals; thus the
bracketing recursion relations give a bijection between the dyadic
rationals and the finite-length orbits of the $\beta$-map. The bracket
tree and the dyadic tree are in one-to-one correspondence. The correspondence
is shown in figure \ref{fig:Beta-Bracket-Map}. 

Perhaps the most interesting aspect to the figure is that it appears
to be continuous. That is, the roots $r_{n}$ appear to be dense in
the interval $1\le\beta\le2$. Although the midpoints $r_{f}$ in
an interval $r_{\ell}<r_{f}<r_{\rho}$ are not evenly spaced, the
midpoint always appears to be sufficiently far away from either endpoint
so that the convergents are dense in the reals, which implies in turn
that the map is continuous. It would appear that the map to infinite-length
sequences can be taken, without pathologies.

\subsubsection{The finite comb}

The figure \ref{fig:Theta-Indicator} shows each indicator rank $\Psi_{\nu}$
as a comb. A better understanding is gained by mapping this to the
infinite binary tree. Each rank $\nu$ corresponds to a single horizontal
row in the tree. The numbering that has been adopted is that $\nu=2$
is at the root of the tree. The dyadic left and right moves on the
comb are
\begin{align*}
L & :\theta_{m}\mapsto\theta_{2m}\\
R & :\theta_{m}\mapsto\theta_{2m+1}
\end{align*}
The indicator function marks each node in the tree with a zero or
a one.

One property of this marking is that when a node is marked with a
zero, all nodes in the right subtree are marked with a zero, as well:
the presence of a zero removes the entire right-hand branch under
that point. This follows from the property noted earlier, that if
$\theta\left(m\right)=0$, then $\theta\left(2^{n}\left(2m+1\right)\right)=0$
for all $n$. Restating in terms of moves, if $\theta_{m}=0$ then,
taking $n=0$, one has $R\theta_{m}=\theta_{2m+1}=0$. Of course it,
follows that $RR\theta_{m}=0$. Taking $n=1$, one has that $LR\theta_{m}=0$,
and so both left and right sub-branches are gone.

A converse marking is given by the property that, if $\theta\left(m\right)=1$
then $\theta\left(2m\right)=1$. On the tree, this means that left-branches
of a node marked with a one are never trimmed: if $\theta_{m}=1$,
then $L\theta_{m}=1$.

Write $2^{<\omega}$ for the collection of all finite strings in the
letters $L,R$ and $2^{\omega}$ as the set of all infinite strings.
The indicator function is then a map $\theta:2^{<\omega}\to\left\{ 0,1\right\} $
that indicates when a given finite, unbounded walk down the tree might
have more branches, or definitely does not. This definitely-maybe
marking defines a filter, and dually, an ideal. The filter/ideal can
be taken on the integers, or equivalently, on the reals. The second
possibility is what the figure \ref{fig:Theta-Indicator} is hinting
at.

The filter has the form that if $n\in\Psi_{\nu}$ then $2n\in\Psi_{v+1}$.
The converse is not true: $10\in\Psi_{5}$ but $5\notin\Psi_{4}$.
Written as a filter, this is $L\Psi_{\nu}\subset\Psi_{\nu+1}$ and
the subset relation is strict. Written as and ideal, if $n\notin\Psi_{\nu}$
then $\left(2n+1\right)\notin\Psi_{v+1}$. This is more clearly stated
with set complements at a given rank. Define the unit interval $I_{\nu}=\left\{ n:2^{\nu-2}\le n<2^{\nu-1}\right\} $
and the set complement $\overline{\Psi}_{\nu}=I_{\nu}\backslash\Psi_{\nu}$.
The ideal is then $R\overline{\Psi}_{\nu}\subset\overline{\Psi}_{\nu+1}$
together with $\left(LR\overline{\Psi}_{\nu}\cup RR\overline{\Psi}_{\nu}\right)\subset\overline{\Psi}_{\nu+2}$.

\subsubsection{The comb bijection}

The good-index bijection $G:\mathbb{N}\to\Psi$ provides a mechanism
to map the full binary tree into the trimmed tree, so that every node
in the trimmed tree is in one-to-one correspondence with the full
binary tree.

The leadership function provides the desired mapping. If $n\in\Psi_{\nu}$,
then $\Lambda\left(n\right)\in\Psi_{\nu+1+h}$ where $h$ was the
height of the leader. Any valid index $n$ gets ``kicked upstairs''
by the leadership function; one can write $\Lambda\Psi_{\nu}\subset\bigcup_{h=0}^{\infty}\Psi_{\nu+1+h}$.
This just corresponds to the bracket move $\mathfrak{R}:\theta_{m}\to\theta_{\Lambda\left(m\right)}$
on the trimmed tree: given any location $\theta\left(m\right)=1$
in the trimmed tree, the bracket right move returns the next valid
right branch in the trimmed tree: $\mathfrak{R}\left(m\right)=\Lambda\left(m\right)$,
since, by construction, $\theta\left(\Lambda\left(m\right)\right)=1$
whenever $\theta\left(m\right)=1$.

\subsubsection{The infinite comb}

The comb is mapped to the reals by defining open subsets bounded by
the dyadic rationals. Let 
\begin{equation}
I\left(m,\nu\right)=\left\{ x\in\mathbb{R}:m<x2^{\nu-2}<\left(m+1\right);0\le m<2^{\nu-2}\right\} \label{eq:topo-base}
\end{equation}
so that $I\left(0,2\right)=\left\{ x\in\mathbb{R}:0<x<1\right\} $.
The left and right moves are the obvious ones: $LI\left(m,\nu\right)=I\left(2m,\nu+1\right)$
and $RI\left(m,\nu\right)=I\left(2m+1,\nu+1\right)$. These are interpreted
as the left and right halves of the (fat) Cantor set, on the reals;
the fat Cantor set being taken as the reals with the dyadic rationals
removed.

The filters on $\Psi_{\nu}$ become filters on $I\left(m,\nu\right)$
in the obvious way, by means of a commuting diagram. Write $I\left(\Psi_{\nu}\right)=\bigcup_{m\in\Psi_{\nu}}I\left(m-2^{\nu-2},\nu\right)$
as the union of intervals that cover the nonempty runs in the binary
tree. This is effectively what is being graphed in the figure \ref{fig:Theta-Indicator}.
The limit set is then
\begin{equation}
\overline{\theta}=\bigcap_{\nu=2}^{\infty}I\left(\Psi_{\nu}\right)\label{eq:infinite comb}
\end{equation}
where the notation $\overline{\theta}$ is happily abused to denote
the limit of the indicator function $\theta_{n}$ as a subset of the
reals. This is possible precisely because $\theta_{n}$ can be mapped
to the binary tree, which can then be partitioned as filters and ideals.

Since $2^{\nu-2}\in\Psi_{\nu}$ for all $\nu$, i.e. the leftmost
branch is never trimmed, one easily concludes that $0\in\overline{\theta}$.
Since $\left(2^{\nu-1}-1\right)\in\Psi_{\nu}$ for all $\nu$, the
rightmost branch is never trimmed, and so one may conclude that $1\in\overline{\theta}$.

It should be clear that the set $\overline{\theta}$ is isomorphic
to the Cantor set. The comb bijection shows how to map points in $\overline{\theta}$
back into the untrimmed infinite tree; the full binary tree is isomorphic
to the Cantor set. From this, one concludes that the set $\overline{\theta}$
is uncountable, and can be placed in bijection with the reals. This
is explored in the next section. A measure can be assigned to $\overline{\theta}$.
This is given in the section after next. A more formal examination
of $\overline{\theta}$, tightening up some of the loose language
above, will be given in a later section on eventually-periodic orbits.

\subsubsection{Self-describing orbits}

What is the meaning of the finite and infinite combs? The finite comb
is a mapping of of the valid-index set $\Psi$ to the dyadics. The
infinite comb is the closure of the finite comb in the reals.

A polynomial index $n\in\Psi$ if and only if $n$ encodes a self-describing
finite orbit. That is, $n\in\Psi$ if and only if the real root $r_{n}$
of $p_{n}\left(r_{n}\right)=0$ iterates under the $\beta$-map such
that the iterate reproduces the bit-sequence of $n$. That is, $n\in\Psi$
if and only if the bits $b_{i}$ of the bitstring $2n+1=b_{0}b_{1}b_{2}\cdots b_{\nu}$
are given by the (finite) orbit as 
\[
b_{i}=\Theta\left(T_{\beta}^{i}\left(\frac{r_{n}}{2}\right)-\frac{1}{2}\right)=k_{i}\left(\frac{r_{n}}{2}\right)
\]
If $n$ does not have this property, then $n\notin\Psi$. At each
rank $\nu$, the elements of $\Psi_{v}$ correspond to finite orbits
of length $\nu$. In the limit of $\nu\to\infty$, one gets self-describing
orbits of unbounded (infinite) length.

By construction, if $n\in\Psi_{\nu}$ then $Ln\in\Psi_{\nu+1}$ and
$\Lambda n\in\Psi_{\nu+1+h}$. That is, if $n$ is a self-describing
orbit, then it appears as the prefix of longer self-describing orbits.

Polynomials of infinite order are holomorphic functions. Given an
infinite bitsequence $\left\{ b\right\} =b_{0}b_{1}\cdots$, define
a holomorphic function

\begin{equation}
q^{\left\{ b\right\} }\left(\zeta\right)=1-\sum_{j=0}^{\infty}b_{j}\zeta^{j+1}\label{eq:holomorphic-beta-function}
\end{equation}
If the bit-sequence is finite, in that all $b_{j}=0$ when $j>k$,
then this is related to the polynomials as 
\[
\zeta^{k+1}p_{n}\left(\frac{1}{\zeta}\right)=1-b_{0}\zeta-b_{1}\zeta^{2}-\cdots-b_{k}\zeta^{k+1}
\]
Given any arbitrary sequence $\left\{ b\right\} $, the holomorphic
function $q^{\left\{ b\right\} }\left(\zeta\right)$ will have a single,
unique real, positive root. To make contact with the polynomials,
write this as the reciprocal, so that the root $r$ satisfies $q^{\left\{ b\right\} }\left(1/r\right)=0$.
This root will have some orbit, given by
\[
a_{i}=\Theta\left(T_{\beta}^{i}\left(\frac{r}{2}\right)-\frac{1}{2}\right)=k_{i}\left(\frac{r}{2}\right)
\]
Such an orbit is self-describing if and only if $\left\{ a\right\} =\left\{ b\right\} $.

The claim being pursued here is that the infinite comb contains all
self-describing sequences, and conversely, every element of the comb
corresponds to a self-describing bit-sequence. That is, if
\[
x=\sum_{j=0}^{\infty}\frac{b_{j}}{2^{j}}\in\overline{\theta}
\]
then

\[
x=2\left(1-q^{\left\{ b\right\} }\left(\frac{1}{2}\right)\right)
\]
Every real number $x\in\overline{\theta}$ in the comb corresponds
to such a self-describing orbit. A proof of these claims will be given
in a later section, after the development of some formal definitions. 

In the meanwhile, it can be noted that every rational number corresponds
to a bitsequence that is ultimately periodic. After an initial unstable
finite sequence, the bitstring settles down to a cyclic orbit. One
task ahead is to examine the set $\mathbb{Q}\cap\overline{\theta}$:
this is the set of self-describing eventually-periodic infinite-length
orbits. It turns out these can be readily described as root of a finite
polynomial. The holomorphic function $q^{\left\{ b\right\} }\left(\zeta\right)$
to be factored into two finite polynomials, one describing the initial
aperiodic segment, and a second describing the cyclic segment. Such
orbits are examined in a later section.

\subsection{Formal definitions}

A sufficient number of distinct concepts have been introduced, that
some basic housekeeping is in order. The definitions that follow are
straight-forward and conventional. The goal is to provide a workable
vocabulary for further discussion.

Let $\mathbb{B}$ denote the finite but unbounded binary tree, and
$\overline{\mathbb{B}}$ it's closure as the infinite tree. The infinite
tree is, of course, isomorphic to the Cantor space $2^{\omega}$;
but this mechanism is not currently needed. A more careful definition
of the finite tree is needed. Thus, let $\mathbb{B}$ the the graph
of vertices and connecting edges $\mathbb{B}=\left\{ v_{i},e_{ij}:i\in\mathbb{N},j\in\left\{ 2i,2i+1\right\} \right\} $.
Let $\eta:\mathbb{N}\to\mathbb{B}$ denote the canonical labeling
of the binary tree by the positive integers, so that the root of the
tree is given the label 1, the first row is 2,3 and the next row is
4,5,6,7. This is a bijection: every finite walk down the tree can
be labeled with a positive integer. The walks are left and right moves,
in the canonical sense: $L:\mathbb{N}\to\mathbb{N}$ with $L:m\mapsto2m$
and likewise $R:\mathbb{N}\to\mathbb{N}$ with $R:m\mapsto2m+1$.
The pushforward of $L,R$ provide the canonical walks on the tree,
as a commuting diagram, so that $L\circ\eta=\eta\circ L$ and likewise
$R\circ\eta=\eta\circ R$. It is useful to adjoin the pre-root elements
$\left\{ 0,\infty\right\} $ so that $R:0\to1$ and $L:\infty\to1$.

Let $\mathbb{D}$ be the dyadic rationals, with the canonical bijection
to the natural numbers $\delta:\mathbb{D}\to\mathbb{N}$ given by
$\delta:\left(2n+1\right)/2^{m}\mapsto2^{m-1}+n$. This labels the
root of the tree with 1/2, and the first row under it as 1/4 and 3/4.
The left and right moves are $L:\left(2n+1\right)/2^{m}\mapsto\left(4n+1\right)/2^{m+1}$
and $R:\left(2n+1\right)/2^{m}\mapsto\left(4n+3\right)/2^{m+1}$.
This was set up so that $L\circ\delta=\delta\circ L$ and likewise
$R\circ\delta=\delta\circ R$. The map $\delta^{-1}:\mathbb{N}\to\mathbb{D}$
is equally familiar: it is $n=b_{0}b_{1}\cdots b_{k}\mapsto\sum_{n=0}^{k}b_{n}2^{-n-1}$
that interprets $n\in\mathbb{N}$ as a sequence of $L,R$ moves in
the binary tree, returning the fraction found at that location.

This allows the bracket tree and the good-index bijection to be specified
more precisely. The validity set $\Psi\subset\mathbb{N}$ is defined
in eqn \ref{eq:valid index set} as the collection of natural number
indexes corresponding to polynomials with self-describing roots. This
set came with two functions $L,\mathfrak{R}:\Psi\to\Psi$ given by
index-doubling $L:m\mapsto2m$ and leadership (eqn\ref{eq:Leadership function})
$\mathfrak{R}:m\mapsto\Lambda\left(m\right)$. The good-index bijection
$G:\mathbb{N}\to\Psi$ is then defined as the pullback $L\circ G=G\circ L$
and $\mathfrak{R}\circ G=G\circ R$. 

The trimmed tree $\mathfrak{B}\subset\mathbb{B}$ is the image of
$\Psi$ under the mapping $\eta$, so that $\mathfrak{B}=\eta\Psi$.
This consists of those nodes and edges in the finite binary tree that
are labeled by integers from the validity set $\Psi$. The good-index
function $G$ places the trimmed tree and the finite tree in a bijection,
so that $\mathfrak{B}=\eta G\eta^{-1}\mathbb{B}$.

The root map $r:\Psi\to\left[1,2\right]$ takes valid integer indexes,
and maps them to the roots of the corresponding (finite-orbit) polynomials,
so that $p_{n}\left(r_{n}\right)=0$. The bracket map, depicted in
figure \ref{fig:Beta-Bracket-Map}, can now be written more precisely
as the function $\varrho=r\circ G\circ\delta^{-1}:\mathbb{D}\to\left[1,2\right]$.
It maps the dyadics to the $\beta$ values that have finite orbits.

The notation for the finite comb $\theta$ is abused in several ways.
In eqn \ref{eq:theta mask} it is used to indicate whether a given
polynomial root has a self-describing orbit. It is then defined as
the indicator function for set membership $\theta=\mathds{1}_{\Psi}:\Psi\to\left\{ T,F\right\} $
and finally as the finite comb $\theta\subset\mathbb{D}$. Formally,
the finite comb can now be written as $\theta=\delta^{-1}\Psi$.

The good-index bijection $G:\mathbb{N}\to\Psi$ can be commuted with
the dyadic bijection to obtain the ``good dyadic map'' $\mathsf{G}=\delta\circ G\circ\delta^{-1}:\mathbb{D}\to\theta$.
Given any dyadic fraction, this map returns another dyadic that lies
within the finite comb. As all the other maps discussed so far, it
is a bijection. It is depicted in figure \ref{fig:Good-Map}.

The infinite comb $\overline{\theta}$ was the closure of the finite
comb, so that $\theta\subset\overline{\theta}$. The closure, defined
in eqn \ref{eq:infinite comb}, is constructed as the infinite intersection
of open sets. Corresponding to this is the closure $\overline{\mathfrak{B}}\subset\overline{\mathbb{B}}$
of infinitely long paths in the trimmed tree. It is taken as the limit
of the finite but unbounded-length paths in $\mathfrak{B}$. The dyadics
can be closed in several ways; one way is to the rationals $\mathbb{Q}$,
and then further to the reals $\mathbb{R}$. In the present case,
we restrict attention to the unit interval $I=\left[0,1\right]\subset\mathbb{R}$
and to $\mathbb{Q}_{I}=\mathbb{Q}\cap\left[0,1\right]$. Proper diligence
requires distinguishing $\overline{\mathfrak{B}}_{\mathbb{Q}}$ from
$\overline{\mathfrak{B}}_{\mathbb{R}}$ and also $\overline{\theta}_{\mathbb{Q}}$
from $\overline{\theta}_{\mathbb{R}}$.

The rationals correspond to infinite length orbits that are ultimately
periodic; these will be examined in a later section. A closure to
the root map to the rationals appears naturally; it can be written
as $\overline{r}:\overline{\theta}\to\left[1,2\right]$. Combined
with the closure $\mathsf{G}:\mathbb{Q}_{I}\to\overline{\theta}$
this gives a closure of the bracket map $\overline{\varrho}=\overline{r}\circ\mathsf{G}:\mathbb{Q}_{I}\to\left[1,2\right]$.
Closing to the reals gives $\mathsf{G}:\left[0,1\right]\to\overline{\theta}$
and $\overline{\varrho}=\overline{r}\circ\mathsf{G}:\left[0,1\right]\to\left[1,2\right]$.
This is the primary achievement of this section: the definition of
the bracket map as a continuous monotonic ascending bijection between
the reals in the unit interval, and $\beta$ values understood as
self-describing orbits. Proving continuity will require some theoretical
machinery and lemmas; these will be developed in a later section.
However, the net result can already be seen in figure \ref{fig:Beta-Bracket-Map}.

\subsection{The Good Map and Measures}

The formal definitions allow the description of the bracket map to
be completed. This is a matter of reviewing the remaining congruences.
The fundamental one, moving forward, is the ``good map'' $\mathsf{G}:\mathbb{D}\to\theta$
that is a correspondence between the dyadic rationals to the comb.
The comb, in turn, can be understood as inducing a measure. The measure,
combined with the good map, induces the ``beta measure'', which
is just the bracket map. This closes the circle of commuting diagrams.

\subsubsection{The good map}

\begin{figure}
\caption{The Good Map\label{fig:Good-Map}}

\includegraphics[width=1\columnwidth]{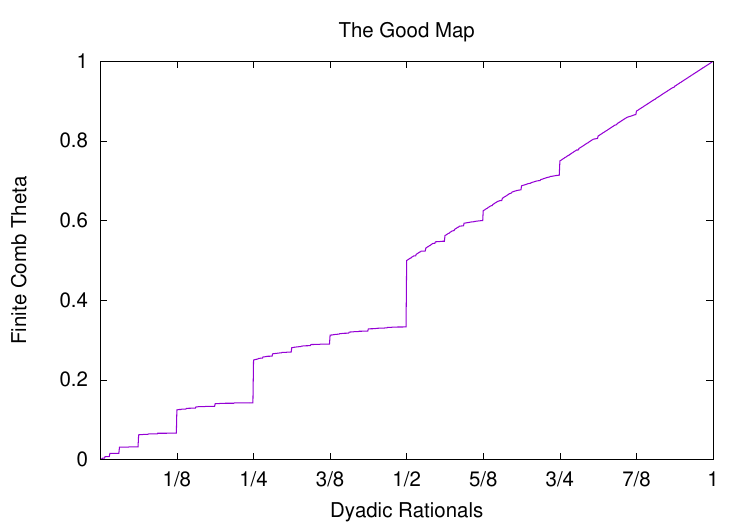}

A graph of the ``good map'' $\mathsf{G}:\mathbb{D}\to\theta$. It
is defined as the commutator $\mathsf{G}=\delta\circ G\circ\delta^{-1}$
of the valid-index map $G:\mathbb{N}\to\Psi$ by the canonical mapping
$\delta:\mathbb{N}\to\mathbb{D}$ between the natural numbers and
the dyadic rationals. The finite comb is just a mapping of the good
indexes $\Psi$ to the dyadics: $\theta=\delta^{-1}\Psi$.

\rule[0.5ex]{1\columnwidth}{1pt}
\end{figure}

The validity map $G:\mathbb{N}\to\Psi$ can be commuted with the canonical
mapping $\delta:\mathbb{N}\to\mathbb{D}$ between the natural numbers
and the dyadic rationals. This defines a map $\mathsf{G}=\delta\circ G\circ\delta^{-1}$
that is a bijection between the dyadic rationals and the finite comb:
$\mathsf{G}:\mathbb{D}\to\theta$. It is shown in figure \ref{fig:Good-Map}. 

\subsubsection{The comb measure}

\begin{figure}
\caption{Limit of Indicator Sum\label{fig:Indicator-Sum}}

\includegraphics[width=1\columnwidth]{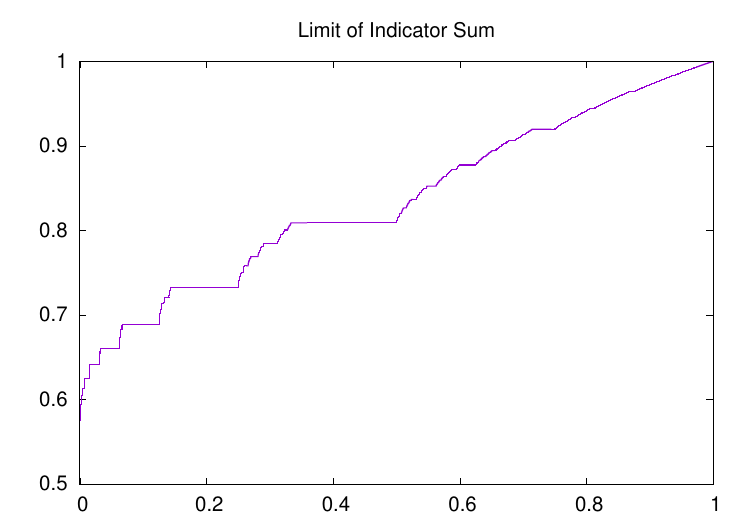}

This figure shows the limit of the indicator sum $A\left(x\right)=\lim_{\nu\to\infty}A_{\nu}\left(x\right)$
of eqn \ref{eq:indicator-sum}. More precisely, it shows $A_{\nu}\left(x\right)$
for $\nu=20$. By this point, convergence is sufficient that any differences
from the limit are not visible to the naked eye.

\rule[0.5ex]{1\columnwidth}{1pt}
\end{figure}

The infinite comb was constructed as a closure or limit of the finite
comb. An interesting trick is to interpret it as a measure, and so
to integrate over it. This can be obtained as a limit over sums of
ranks in the finite comb. The measure is depicted in figure \ref{fig:Indicator-Sum}.

The sum over the indicator function $S\left(k\right)=\sum_{n=1}^{k}\theta_{n}$
shows power-of-two periodicity, same as each rank in the finite comb.
Each rank $\nu$ can be separated out as

\[
S_{\nu}\left(k\right)=S\left(k\right)-S\left(2^{\nu-2}-1\right)
\]
The intent is to isolate the range $2^{\nu-2}\le k<2^{\nu-1}$, and
to split the sum $S\left(k\right)$ into a collection of ranks $S_{\nu}$. 

At the end of the range, the sum $S_{\nu}$ achieves Moreau's necklace-counting
function: $S_{\nu}\left(2^{\nu-1}-1\right)=M_{\nu}$. Dividing by
$M_{\nu}$ gives each rank the same vertical scale: zero to one. It
is also useful to rescale the horizontal range, to run zero-to-one
as well. This gives a normalized version
\[
F_{\nu}\left(x\right)=\frac{S_{\nu}\left(\left\lfloor 2^{\nu-2}\left(1+x\right)\right\rfloor \right)}{M_{\nu}}
\]
that runs from zero to one as $x$ runs from zero to one. This function
does not have an interesting limit as $\nu\to\infty$, as it slowly
drops to zero over the entire unit interval. However, it does so at
a fixed rate, and can be held constant with a radical. The sequence
of functions
\begin{equation}
A_{\nu}\left(x\right)=\exp\left(M_{\nu}2^{-\nu}\log F_{\nu}\left(x\right)\right)\label{eq:indicator-sum}
\end{equation}
converge rapidly and more-or-less uniformly to a limit $A\left(x\right)=\lim_{\nu\to\infty}A_{\nu}\left(x\right)$.
It appears to be well-defined on the unit real interval $0\le x\le1$.
This is the limit shown in figure \ref{fig:Indicator-Sum}. It is
perhaps useful to keep in mind the asymptotic limit of Moreau's function,
$M_{\nu}=2^{\nu}/\nu-\mathcal{O}\left(2^{\nu/2}/\nu\right)$ and so
the radical scales as $M_{\nu}2^{-\nu}=1/\nu-\mathcal{O}\left(2^{-\nu/2}/\nu\right)$.

The convergence to the limit appears to be uniform and rapid, except
at $x=0$, which proceeds slowly. This is easily demonstrated. The
$x=0$ limit is 
\[
A\left(0\right)=\lim_{\nu\to\infty}\left(F_{\nu}\left(0\right)\right)^{M_{\nu}2^{-\nu}}=\lim_{\nu\to\infty}\left(\frac{1}{M_{\nu}}\right)^{M_{\nu}2^{-\nu}}=\lim_{\nu\to\infty}\left(\frac{\nu}{2^{\nu}}\right)^{1/\nu}=\frac{1}{2}\lim_{\nu\to\infty}\nu^{1/\nu}=\frac{1}{2}
\]
The slow convergence is entirely due to the last limit, above.

The offsets to the sums and limit above were defined above, so as
to avoid having to debate the meaning of $\lim_{\nu\to\infty}\sqrt[\nu]{0}$.
Yet clearly, the intent is that $A\left(x\right)$ should provide
a measure for the infinite comb. Should the infinite comb be thought
of as having a large point-weight at $x=0$? Perhaps not; thus, perhaps
a more suitable measure is $\mu_{\theta}\left(x\right)=2A\left(x\right)-1$,
which runs from zero to one over the unit interval.

To summarize: The function $\mu_{\theta}\left(x\right)$ provides
a measure for the infinite comb given in eqn \ref{eq:infinite comb}.
The non-flat sections are where the infinite-length self-describing
sequences are accumulating.

\subsubsection{The Beta measure}

\begin{figure}
\caption{The Beta Measure\label{fig:Beta-measure}}

\includegraphics[width=1\columnwidth]{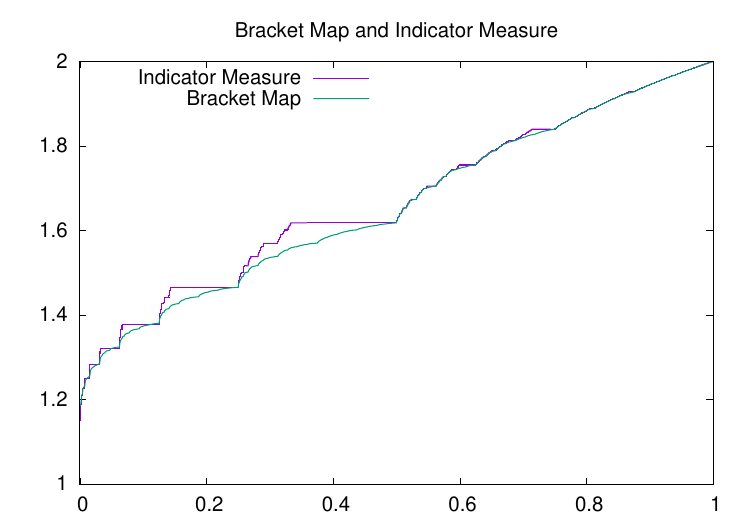}

This figure superimposes the two figures \ref{fig:Beta-Bracket-Map}
and \ref{fig:Indicator-Sum} into one. The indicator measure has been
rescaled, so that the y-axis aligns with the interval $1\le\beta\le2$,
as that is the nominal topic of discussion. The comb is also a map
through $\beta$ values, but taken sideways, as it were. This figure
indicates visually what that correspondence is. When the stair-treads
are removed with the ``good map'', the two curves are identical.
This is surprising, as they have entirely different origins: the indicator
measure is counting periodic orbits, while the bracket map is providing
the locations.

\rule[0.5ex]{1\columnwidth}{1pt}
\end{figure}

The significance of the comb measure is revealed by superimposing
it's graph \ref{fig:Indicator-Sum} on the bracket map \ref{fig:Beta-Bracket-Map},
so that both appear side by side. This is shown in \ref{fig:Beta-measure}.
The stair-treads line up with the blancmange dips in the bracket map.
The good map can be used to eliminate the stair treads. The identity
is 
\[
2A\circ\mathsf{G}=\varrho=r\circ G\circ\delta^{-1}=r\circ\delta^{-1}\circ\mathsf{G}
\]
Since $\mathsf{G}=\delta\circ G\circ\delta^{-1}$ is a bijection,
it can be peeled off, to give
\[
2A=r\circ\delta^{-1}
\]
This is surprising. The map $r$ is purely local: it takes integers
$n$ to the enumerated roots $p_{n}\left(r_{n}\right)=0$. It is just
specifying locations of roots. The map $A$ is global: it is counting
how many roots there are. It is the limit of a sum, a kind of peculiar
integral, that captures information about all other $\beta$ values,
and how they behaved.

An alternative interpretation is that this provides a way of estimating
the number of orbits of length $\nu$, satisfying $\beta<\alpha$
for some fixed $1\le\alpha\le2$. The total number of orbits of length
$\nu$ is given by Moreau's $M_{\nu}$. The counting function $S_{\nu}\left(k\right)$
returns the total number of orbits of length $\nu$ that occurred
at some $\beta$ with $\beta\le r_{k}$. The normalized version provides
this same number slightly more elegantly: the total number of orbits
of length $\nu$ with $\beta<\varrho\left(x\right)$ is given by $M_{\nu}F_{\nu}\left(x\right)$.
The indirection with $\varrho$ is annoying; define $f_{\nu}\left(\alpha\right)=F_{\nu}\left(\varrho^{-1}\left(\alpha\right)\right)$,
so that $f_{\nu}\left(\alpha\right)$ counts the fraction of all orbits
of length $\nu$ occuring at some (any) $\beta\le\alpha$. This fraction
is approximated as $f_{\nu}\left(\alpha\right)\approx\left(\alpha/2\right)^{\nu}$,
which holds exactly in the limit $\nu\to\infty$. 

This last again illustrates the local-global tie between these two:
$f_{\nu}\left(\alpha\right)$ is a counting function, an integral
of sorts, while $\alpha$ is just a number. It's quite rare to find
such specific analytic results in this project. Very rare: this is
the first.

\subsection{Complex Roots}

What are the complex roots? Numerical work clearly indicates that
they seem to be approximately cyclotomic in some sense or another.
They seem to be more-or-less uniformly distributed in an approximate
circle, always. The modulus of most of the complex roots appear to
be less than one. This is violated for the complex roots of $p_{2^{k}}\left(\beta\right)=\beta^{k+2}-\beta^{k+1}-1$,
where some of the roots in the right-hand quadrant have a modulus
larger than one. By contrast, the complex roots of $p_{2^{k}-1}\left(\beta\right)=\beta^{k+1}-\sum_{j=0}^{k}\beta^{j}$
seem to always have a modulus less than one. These two seem to be
the extreme cases: in general, the polynomials appear to be ``approximately
cyclotomic''. Its not clear how to make this statement more precise.

These numerical results can be argued heuristically: just divide the
polynomial by it's leading order. That is, a general polynomial of
this form is
\[
p_{n}\left(z\right)=z^{k+1}-\sum_{j=0}^{k}b_{j}z^{k-j}
\]
with the convention that $b_{0}=b_{k}=1$, and the bit-sequence $2n+1=b_{0}b_{1}b_{2}\cdots b_{p}$
corresponding to a terminating orbit. Dividing by $z^{k+1}$gives
a series 
\[
1-\frac{1}{z}-\frac{b_{1}}{z^{2}}-\frac{b_{2}}{z^{3}}-\cdots
\]
Clearly, this can have a zero only when $\left|z\right|<2$ as otherwise,
the terms get too small too quickly. 

\subsection{$\beta$-Golden $\beta$-Fibonacci Sequences}

It is well known that the golden ratio occurs as limit of the ratio
of adjacent Fibonacci numbers:
\[
\varphi=\lim_{m\to\infty}\frac{F_{m}}{F_{m-1}}
\]
where $F_{m}=F_{m-1}+F_{m-2}$. There is a generalization of this,
which also has received attention: the tribonacci, quadronacci, \emph{etc}.
sequences, whose limits are 
\[
\alpha_{n}=\lim_{m\to\infty}\frac{F_{m}^{(n)}}{F_{m-1}^{(n)}}
\]
where $F_{m}^{(n)}=F_{m-1}^{(n)}+F_{m-2}^{(n)}+\cdots+F_{m-n}^{(n)}$. 

Is it possible that the real roots of the polynomials $p_{n}(\beta)$
are also the roots of such sequences? But of course they are! Given
a finite string of binary digits $\left\{ b\right\} =b_{0},b_{1},\cdots,b_{k}$,
write the beta-Fibonacci sequence as
\[
F_{m}^{\{b\}}=b_{0}F_{m-1}^{\{b\}}+b_{1}F_{m-2}^{\{b\}}+\cdots b_{k}F_{m-k}^{\{b\}}
\]
The name ``beta-Fibonacci'' is needed because the term ``generalized
Fibonacci sequence'' is already in wide circulation for the special
case of all bits being one. The ratio of successive terms is 
\[
\alpha^{\{b\}}=\lim_{m\to\infty}\frac{F_{m}^{\{b\}}}{F_{m-1}^{\{b\}}}
\]
and is given as the (positive) real root of the polynomial
\[
p_{n}\left(\beta\right)=\beta^{k+1}-b_{0}\beta^{k}-b_{1}\beta^{k-1}-\cdots-b_{k}=0
\]
These polynomials and their roots were already enumerated and tabulated
in the previous section. 

The beta-Fibonacci sequences do not appear by accident: these sequences
have an ordinary generating function (OGF) given by the polynomial!
That is,
\[
\sum_{m=0}^{\infty}z^{m}F_{m}^{\{b\}}=\frac{z^{k}}{1-b_{0}z-b_{1}z^{2}-\cdots-b_{k}z^{k+1}}=\frac{1}{zp_{n}\left(\frac{1}{z}\right)}
\]
The factor of $z^{k}$ in the numerator serves only to initiate the
sequence so that $F_{0}^{\{b\}}=\cdots=F_{k-1}^{\{b\}}=0$ and $F_{k}^{\{b\}}=1$.

These sequences are generic: they indicate how many different ways
one can partition the integer $m$ into elements of the set $\left\{ b_{0},2b_{1},3b_{2},\cdots,\left(k+1\right)b_{k}\right\} $.
So, for example, the entry for $n=12$ in the table below corresponds
to OEIS A079971, which describes the number of ways an integer $m$
can be partitioned into 1, 2 and 5 (or that $5m$ can be partitioned
into nickels, dimes and quarters). This corresponds to the bit sequence
$\left\{ b\right\} =11001$; that is, $\left\{ b_{0},2b_{1},3b_{2},\cdots,\left(k+1\right)b_{k}\right\} =\left\{ 1\cdot1,2\cdot1,3\cdot0,4\cdot0,5\cdot1\right\} =\left\{ 1,2,5\right\} $.
From such partitions, it appears that one can build partitions of
the positive integers that are disjoint, and whose union is the positive
integers. This suggests a question: can these partitions be expressed
as Beatty sequences?

The previous table is partly repeated below, this time annotated with
the OEIS sequence references.

\medskip{}

\begin{center}
\begin{tabular}{|c|c|c|c|c|}
\hline 
$n$ & binary & root & root identity & sequence\tabularnewline
\hline 
0 & 1 & 1 &  & \tabularnewline
\hline 
1 & 11 & $\varphi=\frac{1+\sqrt{5}}{2}=1.618$ & golden ratio & Fibonacci A000045\tabularnewline
\hline 
2 & 101 & $1.465571231876\cdots$ & OEIS A092526 & Narayana A000930\tabularnewline
\hline 
3 & 111 & $1.839286755214\cdots$ & tribonacci A058265 & tribonacci A000073\tabularnewline
\hline 
4 & 1001 & $1.380277569097\cdots$ & 2nd Pisot A086106 & A003269, A017898\tabularnewline
\hline 
6 & 1101 & $1.754877666246\cdots$ & OEIS A109134 & A060945\tabularnewline
\hline 
7 & 1111 & $1.927561975482\cdots$ & tetranacci A086088 & tetranacci A000078\tabularnewline
\hline 
8 & 10001 & $1.324717957244\cdots$ & silver A060006 & A003520, A017899\tabularnewline
\hline 
10 & 10101 & $1.570147312196\cdots$ & Pisot A293506 & A060961\tabularnewline
\hline 
12 & 11001 & $1.704902776041\cdots$ &  & A079971\tabularnewline
\hline 
13 & 11011 & $1.812403619268\cdots$ &  & A079976\tabularnewline
\hline 
14 & 11101 & $1.888518845484\cdots$ &  & A079975\tabularnewline
\hline 
15 & 11111 & $1.965948236645\cdots$ & pentanacci A103814 & A001591\tabularnewline
\hline 
16 & 100001 & $1.285199033245\cdots$ &  & A005708, A017900\tabularnewline
\hline 
\end{tabular}
\par\end{center}

\medskip{}

All of these integer sequences and roots participate in a number of
curious relations having a regular form; this is, of course, the whole
point of listing them in the OEIS. This suggests a question: do the
known relationships generalize to the beta-shift setting?

For example, the polynomials of the form 
\[
\beta^{k+1}-\beta-1=0
\]
are the Lamé polynomials, they arise as solutions to Lamé's equation,
a kind of ellipsoidal harmonic differential equation. In the present
notation, these correspond to polynomials $p_{n}\left(\beta\right)=0$
for $n=2^{k}$.

Another example is the Fibonacci-tribonacci-tetranacci sequence of
``generalized golden means''. These are the roots of the series
for which all $b_{k}=1$, that is, the roots of 
\[
\beta^{k+1}-\beta^{k}-\beta^{k-1}-\cdots-1=0
\]
In the present notation, these would be the polynomials $p_{n}\left(\beta\right)=0$
for $n=2^{k}-1$. Such roots can be rapidly computed by a series provided
by Hare, Prodinger and Shallit\cite{Hare14}: 
\[
\frac{1}{\alpha_{k}}=\frac{1}{2}+\frac{1}{2}\sum_{j=1}^{\infty}\frac{1}{j}{j\left(k+1\right) \choose j-1}\frac{1}{2^{j\left(k+1\right)}}
\]
This series is obtained by making good use of the Lagrange inversion
formula. Here, $\alpha_{k}$ is the $k$'th generalized golden mean,
i.e. the solution $p_{2^{k}-1}\left(\alpha_{k}\right)=0$. Can the
Hare series be extended to provide the roots $r_{n}$ of $p_{n}\left(r_{n}\right)=0$
for general $n$?

Another set of observations seem to invoke the theory of complex multiplication
on elliptic curves, and pose additional questions. So:

The tribonacci root $r_{3}$ is given by
\[
r_{3}=\frac{1}{3}\left(1+\sqrt[3]{19+3\sqrt{33}}+\sqrt[3]{19-3\sqrt{33}}\right)\simeq1.839\cdots
\]

The silver number (plastic number) $r_{8}$ is given by 
\[
r_{8}=\frac{1}{6}\left(\sqrt[3]{108+12\sqrt{69}}+\sqrt[3]{108-12\sqrt{69}}\right)\simeq1.324\cdots
\]

The Narayana's cows number $r_{2}$ is given by 
\[
r_{2}=\frac{1}{6}\sqrt[3]{116+12\sqrt{93}}+\frac{2}{3\sqrt[3]{116+12\sqrt{93}}}+\frac{1}{3}\simeq1.645\cdots
\]

The root $r_{6}$ is related to the silver number $r_{8}$ as $r_{8}=r_{6}\left(r_{6}-1\right)$
and is given by
\[
r_{6}=\frac{1}{6}\sqrt[3]{108+12\sqrt{69}}+\frac{2}{\left(\sqrt[3]{108+12\sqrt{69}}\right)^{2}}\simeq1.754\cdots
\]

Do the other roots have comparable expressions? To obtain them, is
it sufficient to articulate the theory of ``complex multiplication''
on elliptic curves? Are these just disguised solutions to Lamé's ellipsoidal
harmonic differential equation? The appearance of only the cube and
square roots is certainly suggestive of an underlying process of points
on elliptic curves.

\subsection{$\beta$-Fibonacci sequences as shifts}

The nature of the $\beta$-Fibonacci sequences as shift sequences
can be emphasized by noting that they arise from the iteration of
the \href{https://en.wikipedia.org/wiki/Companion_matrix}{companion matrix}
for the polynomial $p_{n}\left(x\right)$. This is a $\left(k+1\right)\times\left(k+1\right)$
matrix in \href{https://en.wikipedia.org/wiki/Hessenberg_matrix}{lower-Hessenberg form}:
\begin{equation}
B=\left[\begin{array}{cccccc}
b_{0} & 1 & 0 & 0 & \cdots & 0\\
b_{1} & 0 & 1 & 0 & \cdots & 0\\
b_{2} & 0 & 0 & 1 & \cdots & 0\\
\vdots & \vdots & \vdots &  & \ddots & \vdots\\
b_{k-1} & 0 & 0 & 0 & \cdots & 1\\
b_{k} & 0 & 0 & 0 & \cdots & 0
\end{array}\right]\label{eq:golden matrix}
\end{equation}
Iteration produces a \href{https://en.wikipedia.org/wiki/Linear_recursive_sequence}{linear recursive sequence}
this is the $\beta$-Fibonacci sequence. The $m$'th element of the
sequence is obtained from the $m$'th iterate $B^{m}$.

Define the \href{https://en.wikipedia.org/wiki/Exchange_matrix}{exchange matrix}
as
\[
J=\left[\begin{array}{ccccc}
0 & 0 & \cdots & 0 & 1\\
0 & 0 & \cdots & 1 & 0\\
\vdots & \vdots & \iddots &  & \vdots\\
0 & 1 & \cdots & 0 & 0\\
1 & 0 & \cdots & 0 & 0
\end{array}\right]
\]
This can be used to write the above in the more conventional companion-matrix
form:
\[
C=\left[JBJ\right]^{\mathrm{T}}=\left[\begin{array}{cccccc}
0 & 1 & 0 & 0 & \cdots & 0\\
0 & 0 & 1 & 0 & \cdots & 0\\
0 & 0 & 0 & 1 & \cdots & 0\\
\vdots & \vdots & \vdots &  & \ddots & \vdots\\
0 & 0 & 0 & 0 & \cdots & 1\\
b_{k} & b_{k-1} & b_{k-2} & b_{k-3} & \cdots & b_{0}
\end{array}\right]
\]

Some explicit examples are in order. For the golden ratio, one has
\[
B=\left[\begin{array}{cc}
1 & 1\\
1 & 0
\end{array}\right]
\]
and the iterates are 
\[
B^{2}=\left[\begin{array}{cc}
2 & 1\\
1 & 1
\end{array}\right],\;B^{3}=\left[\begin{array}{cc}
3 & 2\\
2 & 1
\end{array}\right],\;B^{4}=\left[\begin{array}{cc}
5 & 3\\
3 & 2
\end{array}\right],\;B^{n}=\left[\begin{array}{cc}
F_{n} & F_{n-1}\\
F_{n-1} & F_{n-2}
\end{array}\right]
\]
with $F_{n}$ being the $n$'th Fibonacci number, as usual. For the
general case, after $m\ge k-1$ iterations, one gets the \href{https://en.wikipedia.org/wiki/Hankel_matrix}{Hankel matrix}
\[
B^{m}=\left[\begin{array}{cccccc}
F_{m}^{\left\{ b\right\} } & F_{m-1}^{\left\{ b\right\} } & F_{m-2}^{\left\{ b\right\} } & \cdots & F_{m-k+1}^{\left\{ b\right\} } & F_{m-k}^{\left\{ b\right\} }\\
F_{m-1}^{\left\{ b\right\} } & F_{m-2}^{\left\{ b\right\} } & F_{m-3}^{\left\{ b\right\} } & \cdots & F_{m-k}^{\left\{ b\right\} } & F_{m-k-1}^{\left\{ b\right\} }\\
F_{m-2}^{\left\{ b\right\} } & F_{m-3}^{\left\{ b\right\} } &  & \cdots\\
\vdots & \vdots & \vdots & \ddots & \vdots & \vdots\\
F_{m-k+1}^{\left\{ b\right\} } &  &  & \cdots\\
F_{m-k}^{\left\{ b\right\} } &  &  & \cdots &  & F_{m-2k}^{\left\{ b\right\} }
\end{array}\right]
\]
so that the top row consists of the latest sequence values. When multiplied
by the bits, this just generates the next iterate in the sequence.
The upper-diagonal 1's just serve to shift columns over by one, with
each iteration: that is why it's a shift!

The product $B^{m}J$ is a \href{https://en.wikipedia.org/wiki/Toeplitz_matrix}{Toeplitz matrix}.

The characteristic polynomial of this matrix is, of course, the polynomial
$p_{n}$:
\[
\det\left[B-xI\right]=\left(-1\right)^{k}p_{n}\left(x\right)
\]
Thus, we can trivially conclude that the eigenvalues of $B$ are given
by the roots of $p_{n}\left(x\right)$. This matrix is in lower-Hessenberg
form; this makes it obvious that it's a shift; a finite shift, in
this case.

For every Hankel matrix, there are corresponding \href{https://en.wikipedia.org/wiki/Moment_problem}{moment problems}
that ask if there is a measure that generates the corresponding matrix
entries (as moments). The moment problem on the unit interval is the
\href{https://en.wikipedia.org/wiki/Hausdorff_moment_problem}{Hausdorff moment problem}.
In the present case, it appears that the corresponding measure is
precisely the Parry–Gelfond measure. The only point of confusion is
that the measure is placed in the denominator, not the numerator,
as the measure is for the pushforward, rather than being a pullback.
There do appear to be various interrelationships between the moments,
these remain to be articulated.

\subsection{Infinite sequences and shifts}

Both of the last two sections extend in a natural way for infinite
sequences. The $\beta$-Fibonacci sequences are just as before,
\[
F_{m}^{\{b\}}=\sum_{j=1}^{m}b_{j-1}F_{m-j}^{\{b\}}
\]
starting with $F_{0}^{\left\{ b\right\} }=1$. The sum is always finite,
and all that one needs is the first $m$ bits of the (now infinite)
bit-sequence $\left\{ b\right\} $. The integer sequence still converges
as before,
\[
\beta=\lim_{m\to\infty}\frac{F_{m}^{\{b\}}}{F_{m-1}^{\{b\}}}
\]
with $\beta$ value is the one associated to $\left\{ b\right\} $.
The real number $\beta$ and the bit sequence $\left\{ b\right\} $
label exactly the same orbit.

The shift of eqn \ref{eq:golden matrix} is replaced by an infinite-dimensional
matrix $\mathcal{H}_{\beta}=B^{T}/\beta$. The transpose is take to
put it into upper-Hessenberg form; the subdiagonal provides the shift.
It is rescaled so that the largest eigenvalue is 1. To be explicit:
given the bit-sequence $\left\{ b\right\} $, the operator$\mathcal{H}_{\beta}$
has the matrix elements
\begin{align*}
\left\langle 0\left|\mathcal{H}_{\beta}\right|j\right\rangle  & =\frac{b_{j}}{\beta}\\
\left\langle j+1\left|\mathcal{H}_{\beta}\right|j\right\rangle  & =\frac{1}{\beta}
\end{align*}
with all other entries being zero. It is perhaps useful to visualize
this operator; in pictorial form, it is
\[
\mathcal{H}_{\beta}=\frac{1}{\beta}\left[\begin{array}{ccccc}
b_{0} & b_{1} & b_{2} & b_{3} & \cdots\\
1 & 0 & 0 & 0 & \cdots\\
0 & 1 & 0 & 0 & \cdots\\
0 & 0 & 1\\
\vdots & \vdots & \vdots &  & \ddots
\end{array}\right]
\]

This operator appears to be conjugate to the beta shift via a similarity
transform. There is an operator $S$ such that

\[
\mathcal{L}_{\beta}=S^{-1}\mathcal{H}_{\beta}S
\]
The Parry–Gelfond invariant measure $\mathcal{L}_{\beta}\rho=\rho$
is mapped to $\sigma=S\rho$, where $\mathcal{H}_{\beta}\sigma=\sigma$
is now the Frobenius–Perron eigenvector. This eigenvector is easy
to write down explicitly: it as $\sigma=\left(1,\beta^{-1},\beta^{-2},\cdots\right)$,
so that $\sigma_{j}=\beta^{-j}$. This is easy to verify: the subdiagonal
entries of $\mathcal{B}_{\beta}$ act as a shift on $\sigma$ and
the top row is just
\[
1=\sum_{j=0}^{\infty}\left\langle 0\left|\mathcal{H}_{\beta}\right|j\right\rangle \sigma_{j}=\sum_{j=0}^{\infty}b_{j}\beta^{-j-1}=1-q^{\left\{ b\right\} }\left(\frac{1}{\beta}\right)=1
\]
with $q^{\left\{ b\right\} }\left(\zeta\right)$ the same holomorphic
function as defined earlier, in eqn \ref{eq:holomorphic-beta-function}.

The similarity transformation can be given directly. If $w$ is a
vector satisfying $\mathcal{H}_{\beta}w=\lambda w$, with vector elements
$w_{j}$, then the function
\[
w\left(x\right)=\sum_{j=0}^{\infty}d_{j}\left(x\right)w_{j}
\]
is an eigenfunction of the transfer operator, satisfying $\mathcal{L}_{\beta}w=\lambda w$.
For $\lambda=1$, this is just $w=\sigma$ which is explicitly the
Parry–Gelfond invariant measure.

The conjugacy allows a strong statement: the \emph{only} solutions
to $\mathcal{H}_{\beta}w=\lambda w$ are necessarily of the form $w=\left(1,\left(\lambda\beta\right)^{-1},\left(\lambda\beta\right)^{-2},\cdots\right)$,
because the subdiagonal forces this shift. To satisfy the the top
row of $\mathcal{H}_{\beta}$, one must have that 
\[
\lambda=\sum_{j=0}^{\infty}\left\langle 0\left|\mathcal{H}_{\beta}\right|j\right\rangle v_{j}=\frac{1}{\beta}\sum_{j=0}^{\infty}\frac{b_{j}}{\left(\lambda\beta\right)^{j}}=\lambda\left(E+1\right)=\lambda\left(1-q^{\left\{ b\right\} }\left(\frac{1}{\lambda\beta}\right)\right)=\lambda
\]
and so the eigenvalue $\lambda$ is exactly the eigenvalue that solves
the $\beta$-series $q^{\left\{ b\right\} }\left(1/\lambda\beta\right)=0$,
equivalently, solves $E=1$ with $E$ as defined in eqn \ref{eq:holomorphic-disk}.
This effectively concludes a proof: the solutions to this series are
the only eigenvalues of the $\beta$-transfer operator; there are
no others.

And yet, there is an issue. For sequences $\left\{ b\right\} $ that
are not finite in length, the holomorphic function $q^{\left\{ b\right\} }\left(\zeta\right)$
has a countable number of zeros, arranged in a circular ring at $\lambda=1/\beta$.
These are dense on the circle, and so conventional analytic continuation
cannot be used to find smaller eigenvalues. Does this prove that there
are no smaller eigenvalues? No. Clearly there are eigenfunctions with
eigenvalues smaller than $\lambda=1/\beta$; an explicit example for
$\lambda=1/\beta^{2}$ was already given. Others can be numerically
approximated. Thus, this description, on it's face, is incomplete.
How? Why?

The issue can be localized to the mis-identification of $d_{k}\left(x\right)=\Theta\left(t_{k}-x\right)$
as the full and complete similarity transform; it is not. The $d_{k}\left(x\right)$
are flat plateaus, and are exactly what the Parry–Gelfond measure
is built out of. The example with $\lambda=1/\beta^{2}$ is a piece-wise
assembly of quadratics, not flat plateaus. The correct construction
requires replacing the $d_{k}\left(x\right)$ with fragments of polynomials.
The identification of the $d_{k}\left(x\right)$ as the conjugating
similarity transform is flawed: it gives only a part of the total
space; it cuts off a portion of the spectrum. It does pair Parry–Gelfond
to Frobenius-Perron, which is good. Just that parts are missing.

\subsection{Equivalent labels for orbits}

There are many equivalent ways of labeling the various expressions
and properties under consideration. These are recapped here. 

\subsubsection{Orbits}

For every given $1<\beta<2$ there is a unique orbit of midpoints
$\left\{ m_{p}\right\} $ given by $m_{p}=T_{\beta}\left(m_{p-1}\right)=T_{\beta}^{p}\left(m_{0}\right)$
and $m_{0}=\beta/2$. The orbits are in one-to-one correspondence
with $\beta$. The midpoints are the same as the Parry sequence; namely
$T_{\beta}^{p}\left(\beta/2\right)=\left(\beta/2\right)t_{\beta}^{p}\left(1\right)$,
recalling here the notation of eqn \ref{eq:beta transform} and \ref{eq:invariant measure}.
Some orbits are of finite length; the rest are either eventually periodic,
or are ergodic.

\subsubsection{Orbit encoding}

The midpoint generates a unique sequence of bits $\left\{ b\right\} =\left\{ b_{0},b_{1,},\cdots,b_{k},\cdots\right\} $
given by the left-right moves of the mid-point, as it is iterated.
That is, $b_{k}=\Theta\left(m_{k}-1/2\right)$ so that $b_{k}$ is
one if the midpoint is greater than half, else $b_{k}$ is zero. Each
bit-sequence is in one-to-one correspondence with $\beta$. Finite
orbits have finite-length sequences.

\subsubsection{Monotonicity}

The compressor function $w\left(\beta\right)=\sum_{k}b_{k}2^{-k}$
is a monotonically increasing function of $\beta$, so that values
of $w\left(\beta\right)$ are in one-to-one correspondence with $\beta$. 

\subsubsection{Polynomial numbering}

If the orbit is finite, then there exists a polynomial $p_{n}\left(z\right)=z^{k+1}-b_{0}z^{k}-b_{1}z^{k-1}-\cdots-b_{k-1}z-1$
with $k=1+\left\lfloor \log_{2}\left(2n+1\right)\right\rfloor $ being
the length of the orbit. The positive real root $r_{n}$ of $p_{n}\left(r_{n}\right)=0$
is $\beta=r_{n}$. That is, the iteration of $r_{n}$ will generate
the finite-length bit-sequence $\left\{ b\right\} =\left\{ b_{0},b_{1,},\cdots,b_{k}\right\} $.
The integer $n$ is in one-to-one correspondence with the bit sequence,
and with the value of $\beta$. The integer is explicitly given by
$2n+1=\sum_{j=0}^{k}2^{j}b_{j}$. 

If the orbit is not finite, there is a function $q^{\left\{ b\right\} }\left(\zeta\right)=1-\sum_{j=0}^{\infty}b_{j}\zeta^{j+1}$
holomorphic on the unit disk, having one unique positive real zero
$q^{\left\{ b\right\} }\left(r_{\left\{ b\right\} }\right)=0$ where
this $r_{\left\{ b\right\} }=1/\beta$ is the same $\beta$ that generated
the bit-sequence $\left\{ b\right\} $. Iterating $r_{\left\{ b\right\} }$
generates $\left\{ b\right\} $. If $\left\{ b\right\} $ is finite,
then $q^{\left\{ b\right\} }\left(\zeta\right)=\zeta^{k+1}p_{n}\left(1/\zeta\right)$,
so these functions agree on finite-length sequences.

\subsubsection{Brackets}

If the orbit is finite, then there exists a unique bracketing relationship
$\ell\Mapsto n\Mapsfrom\rho$ for which $n$ is the polynomial index.
The left and right bounds $\ell,\rho$ are strictly smaller indexes:
$\ell<n$ and $\rho<n$, and even more strongly, $2\ell\le n$ and
$2\rho\le n$ that have the property of bounding the positive real
roots of the corresponding polynomials: $r_{\ell}<r_{n}<r_{\rho}$,
with $p_{\ell}\left(r_{\ell}\right)=p_{n}\left(r_{n}\right)=p_{\rho}\left(r_{\rho}\right)=0$.

\subsubsection{Binary tree}

Each bracket is in one-to-one correspondence with a node in the full,
unbounded binary tree. Sub-brackets define left and right subintervals
that are disjoint, and whose union makes up the whole interval. Every
node in the full binary tree can be labeled with a unique sequence
of left-right moves to get to that node. This places the brackets
(and thus, the polynomials and the roots and the mid-point orbits)
in unique, one-to-one correspondences with finite-length strings of
L,R moves. Such strings are, in turn, in one-to-one correspondence
with the dyadic rationals. The L,R strings are in one-to-one correspondence
with the orbits $\left\{ b\right\} $ but they are \emph{not} numerically
the same! This are distinct sequences! In particular, all possible
L,R moves are allowed. Only a limited number of orbits $\left\{ b\right\} $
are possible, as limited by necklace-counting considerations.

\subsubsection{Baire sequences}

If the orbit is finite, then there exists a unique integer sequence
$\left[m_{1},m_{2},\cdots,m_{k}\right]$ such that the index is given
by $\eta\left[m_{1},m_{2},\cdots,m_{k}\right]$. This is a bijection
between all valid indexes and all possible finite-length sequences.
Due to the bounding property of the brackets, limits of $k\to\infty$
can be taken, and these limits are unique. Thus, all $\beta$ values,
including those with non-finite orbits, can be placed in a one-to-one
bijection with infinite-length sequences $\left[m_{1},m_{2},\cdots\right]\in\mathbb{N}^{\omega}$.

\subsubsection{Beta-Fibonacci sequences}

If the orbit is finite, then there exists a sequence of integers $F^{\{b\}}$,
the beta-Fibonacci sequence, that is in one-to-one correspondence
with the finite bit sequence $\left\{ b\right\} =b_{0},b_{1},\cdots,b_{k}$,
and with the value of $\beta$. There are also sequences for each
infinite-length orbit $\left\{ b\right\} $. These are briefly touched
on, below.

\subsubsection{Shift matrix}

If the orbit is finite, then the finite bit sequence $\left\{ b\right\} =b_{0},b_{1},\cdots,b_{k}$
defines a lower-Hessenberg ``golden shift'' matrix $B$, as shown
in eqn \ref{eq:golden matrix}. The limit of $k\to\infty$ can be
taken in a relatively straight-forward manner, given below.

\subsubsection{Summary}

To summarize: any one of these: the integer $n$, the polynomial $p_{n}\left(x\right)$,
the bracket location in the binary tree, a dyadic rational, a point
in Baire space, the integer sequence $F_{m}^{\left\{ b\right\} }$,
the orbit of midpoints $m_{p}=T^{p}\left(\beta/2\right)$, the orbit
encoding $\left\{ b\right\} $, the shift matrix $B$, the value of
the compressor function $w\left(\beta\right)$ and, of course, $\beta$
itself can each be used as a stand-in for the others, as they are
all in one-to-one correspondence. Specifying one determines the others;
all uniquely map to one-another. The formulas that provide maps between
each of these can all be given in closed form, except for the handful
of recursively-defined formulas. The recursive formulas are all invertable;
thus they are computable (decidable). They are all equivalent labels.
Fashionably abusing notation, $n\equiv p_{n}\left(x\right)\equiv r_{n}\equiv\left\{ b\right\} \equiv F_{m}^{\left\{ b\right\} }\equiv m_{p}\equiv w\left(\beta\right)\equiv\beta\equiv B$.

An explicit expression relating the orbit encoding and the orbit can
be read off directly from eqn \ref{eq:iterated shift}. Plugging in,
\begin{equation}
m_{p}=T_{\beta}^{p+1}\left(\frac{\beta}{2}\right)=\frac{\beta}{2}\left[\beta^{p+1}-\sum_{j=0}^{p}b_{j}\beta^{p-j}\right]\label{eq:midpoint poly}
\end{equation}
for $p<k$ the length of the bit sequence, and $m_{k}=T_{\beta}^{k+1}\left(\beta/2\right)=\beta p_{n}\left(\beta\right)/2=0$
terminating, since $\beta$ is the positive root of $p_{n}\left(x\right)$. 

Four of the correspondences given above ask for finite orbits. Three
of these can be extended to non-finite orbits in an unambiguous and
uncontroversial way. The extensions are covered in the next two sections.
The fourth is the numbering $n$ of the finite orbits. These are countable;
there is no way to extend the counting number $n$ to the non-finite
orbits. Indeed, there are too many: the non-finite orbits are uncountable.

\pagebreak{}

\section{Eventually-periodic orbits}

The previous section characterized the set of $\beta$ values which
have midpoint orbits of finite length. Another interesting class is
the set of eventually-periodic orbits: orbits of infinite length,
settling down to a stable, periodic cycle after an initial bout of
chaotic motion. These $\beta$ values occur as the roots of a slightly
different set of self-describing polynomials, as a sum of two parts:
one for the initial chaotic motion, and a second polynomial for the
cyclic motion. These can be (monotonically) paired with the rational
numbers, with infinite paths through the binary tree, and with locations
on the comb function. They naturally extend the bracket mapping from
the dyadic rationals to all rational numbers.

\subsection{Bitsequence polynomials}

Consider the set of all eventually-periodic bit-sequences. These consist
of a leading chaotic prefix of length $L$ followed by a periodic
orbit of length $N$. Such sequences can be placed in one-to-one correspondence
with the rationals, in the conventional fashion. Select such a sequence
$\left\{ b_{k}\right\} $. The cyclic condition has that $b_{k}=b_{k+N}$
for all $k\ge L$. Split this into two parts: a finite length-$L$
digit sequence $d_{k}=b_{k}$ for $k<L$ and a finite length-$N$
cyclic bit sequence $c_{k}=b_{k+L}$ for $k<N$. Associated to this
is a rational number $x=\sum_{k=0}^{\infty}b_{k}2^{-k-1}$. To get
it's value, write
\begin{align*}
S\left(\beta\right)= & \sum_{k=0}^{\infty}b_{k}\beta^{-k}\\
= & \sum_{k=0}^{L-1}b_{k}\beta^{-k}+\sum_{k=L}^{L+N-1}b_{k}\beta^{-k}+\sum_{k=L+N}^{L+2N-1}b_{k}\beta^{-k}+\cdots\\
= & \sum_{k=0}^{L-1}b_{k}\beta^{-k}+\frac{\beta^{N}}{\beta^{N}-1}\sum_{k=L}^{L+N-1}b_{k}\beta^{-k}\\
= & \frac{1}{\beta^{L-1}}\sum_{k=0}^{L-1}d_{k}\beta^{L-k-1}+\frac{1}{\beta^{L-1}\left(\beta^{N}-1\right)}\sum_{k=0}^{N-1}c_{k}\beta^{N-k-1}
\end{align*}
Plugging in $\beta=2$ gives ratios of whole numbers: a rational that
corresponds to this bit-sequence. 

Any given $1\le\beta\le2$ generates a mid-point orbit bit-sequence.
Sadly, we've introduced too many different but equivalent notations
for this. Starting with the mid-point $x=m_{0}=\beta/2$, the characteristic
bit-sequence is

\[
b_{n}=\Theta\left(m_{n}-\frac{1}{2}\right)=\Theta\left(T_{\beta}^{n}\left(\frac{\beta}{2}\right)-\frac{1}{2}\right)=d_{n}\left(\frac{1}{2}\right)=k_{n}\left(\frac{\beta}{2}\right)=\varepsilon_{n}\left(\frac{1}{\beta}\right)
\]
Such orbits have the self-describing property, that $\beta=S\left(\beta\right)$.
This follows from eqn \ref{eq:shift series}, as well as other identities.
For the eventually-periodic orbits, $S\left(\beta\right)$ is a ratio
of polynomials. It can clearly be placed in one-to-one correspondence
with the rationals. It is not hard to show that it has one unique
real root $1<\beta\le2$, and so these are in one-to-one correspondence
as well.

This can be related to earlier notation. Writing $\beta-S\left(\beta\right)=0$
and then multiplying through by the denominator, 

\[
\left(\beta^{N}-1\right)\left(\beta^{L}-\sum_{k=0}^{L-1}d_{k}\beta^{L-1-k}\right)-\sum_{k=0}^{N-1}c_{k}\beta^{N-1-k}=0
\]
This touches the earlier notation for the polynomials describing finite
orbits:
\[
p^{\left\{ d_{0}\cdots d_{\nu-1}\right\} }\left(\beta\right)=\beta^{\nu}-\sum_{k=0}^{\nu-1}d_{k}\beta^{\nu-1-k}
\]
and so the eventually-periodic orbits are the roots of
\[
\left(\beta^{N}-1\right)p^{\left\{ d_{0}\cdots d_{L-1}\right\} }\left(\beta\right)+p^{\left\{ c_{0}\cdots c_{N-1}\right\} }\left(\beta\right)-\beta^{N}=0
\]
Finite-length orbits have $N=1$ and $c_{0}=0$, so that the cyclic
term vanishes, and the earlier form is recovered, after ignoring an
extra factor of $\beta-1$. 

Finite-length orbits can be converted to infinite-length periodic
orbits simply by setting the last bit to zero, and then repeating
the bit sequence cyclically. Recall, finite orbits always have the
last bit equal to one, so this operation is always unambiguous. Restating
this explicitly: if $b_{0}b_{1}\cdots b_{\nu-1}$ is a finite-length
orbit, then set $L=0$ and $N=\nu$ and $c_{k}=b_{k}$ for $k<\nu-1$
and finally $c_{\nu-1}=0$. Plugging through just yields the usual
polynomial for finite orbits.

Collecting terms of the same order, write
\[
0=\beta^{N+L}-\sum_{k=0}^{N+L-1}a_{k}\beta^{N+L-1-k}
\]
The coefficients are then
\[
a_{k}=d_{k}+c_{k-L}+\delta_{k-N-1}^{0}-d_{k-N}
\]
with $\delta$ the Kronecker delta. The only coefficients that can
appear are in the set $a_{k}\in\left\{ -1,0,1,2\right\} $, and at
most one coefficient can be 2; it is specifically $a_{N-1}$, which
can only ever be 1 or 2. This can be best understood visually, lining
up columns. For $L>N$, get $a_{k}$ by summing the columns:\medskip{}

\begin{center}
\begin{tabular}{cccccccccc}
 & $d_{0}$ & $\cdots$ &  &  & $\cdots$ & $d_{L-1}$ & $c_{0}$ & $\cdots$ & $c_{N-1}$\tabularnewline
$+$ & $0$ & $\cdots$ & $0$ & 1 & -$d_{0}$ &  & $\cdots$ &  & -$d_{L-1}$\tabularnewline
\hline 
 & $a_{0}$ &  & $\cdots$ &  & $a_{N}$ &  & $\cdots$ &  & $a_{N+L-1}$\tabularnewline
\end{tabular}\medskip{}
\par\end{center}

\begin{flushleft}
The top row shows the prefix string and the first run of the cyclic
string. The second row shows minus the prefix string, shifted all
the way to the right, establishing alignment at the right side. Just
before it is a lone 1, and then padded on the left with zeros. The
bottom row is just the sum of the two rows above it. Clearly, the
coefficient $a_{N-1}$ can only ever be 1 or 2. It is exactly the
same for $L<N$, with the bottom row being shorter: \medskip{}
\par\end{flushleft}

\begin{center}
\begin{tabular}{cccccccccc}
 & $d_{0}$ & $\cdots$ & $d_{L-1}$ & $c_{0}$ &  &  & $\cdots$ &  & $c_{N-1}$\tabularnewline
$+$ & $0$ &  & $\cdots$ &  & 0 & 1 & -$d_{0}$ & $\cdots$ & -$d_{L-1}$\tabularnewline
\hline 
 & $a_{0}$ &  &  & $\cdots$ &  &  & $a_{N}$ & $\cdots$ & $a_{N+L-1}$\tabularnewline
\end{tabular}\medskip{}
\par\end{center}

\subsection{Examples}

Not all possible bit-sequences $\left\{ b_{k}\right\} $ are allowed.
They must be self-describing, so that the root of $\beta-S\left(\beta\right)=0$,
when iterated, generates $\left\{ b_{k}\right\} $. For orientation,
some examples are shown below.

\medskip{}
\begingroup
\def\arraystretch{1.3}
\begin{center}
\begin{tabular}{|c|c|c|c|c|c|c|c|}
\hline 
$\frac{p}{q}$ & $\frac{p+q}{2q}$ & $m,c$ & \negthinspace{}$b_{0}b_{1}\cdots${\footnotesize{}\negthinspace{}} & $L$ & $N$ & polynomial & $\beta$\tabularnewline
\hline 
\hline 
{\footnotesize{}0} & {\footnotesize{}1/2} & {\footnotesize{}\negthinspace{}1,0\negthinspace{}} & {\footnotesize{}$1\overline{\cdot}$} & {\footnotesize{}1} & {\footnotesize{}1} & {\footnotesize{}$\beta-1=0$} & {\footnotesize{}1}\tabularnewline
\hline 
{\footnotesize{}1} & {\footnotesize{}1} & {\footnotesize{}\negthinspace{}0,1\negthinspace{}} & {\footnotesize{}$\overline{1}$} & {\footnotesize{}0} & {\footnotesize{}1} & {\footnotesize{}$\beta-2=0$} & {\footnotesize{}2}\tabularnewline
\hline 
{\footnotesize{}1/2} & {\footnotesize{}3/4} & {\footnotesize{}\negthinspace{}3,0\negthinspace{}} & {\footnotesize{}$11\overline{\cdot}$} & {\footnotesize{}2} & {\footnotesize{}1} & {\footnotesize{}$\beta^{2}-\beta-1=0$} & {\footnotesize{}$\varphi=1.61803\cdots$}\tabularnewline
\hline 
{\footnotesize{}1/3} & {\footnotesize{}2/3} & {\footnotesize{}\negthinspace{}0,2\negthinspace{}} & {\footnotesize{}$\overline{10}$} & {\footnotesize{}0} & {\footnotesize{}2} & {\footnotesize{}above} & {\footnotesize{}above}\tabularnewline
\hline 
{\footnotesize{}2/3} & {\footnotesize{}5/6} & {\footnotesize{}\negthinspace{}1,2\negthinspace{}} & {\footnotesize{}$1\overline{10}$} & {\footnotesize{}1} & {\footnotesize{}2} & {\footnotesize{}$\beta^{3}-\beta^{2}-2\beta+1=0$} & {\footnotesize{}1.80...$=2\cos\frac{\pi}{7}$}\tabularnewline
\hline 
{\footnotesize{}5/6} & {\footnotesize{}\negthinspace{}11/12\negthinspace{}} & {\footnotesize{}\negthinspace{}3,2\negthinspace{}} & {\footnotesize{}$11\overline{10}$\negthinspace{}} & {\footnotesize{}2} & {\footnotesize{}2} & {\footnotesize{}$\beta^{4}-\beta^{3}-2\beta^{2}+1=0$} & {\footnotesize{}1.90516616775}\tabularnewline
\hline 
{\footnotesize{}\negthinspace{}11/12\negthinspace{}} & {\footnotesize{}\negthinspace{}23/24\negthinspace{}} & {\footnotesize{}\negthinspace{}7,2\negthinspace{}} & {\footnotesize{}\negthinspace{}$111\overline{10}$\negthinspace{}} & {\footnotesize{}3} & {\footnotesize{}2} & {\footnotesize{}$\beta^{5}-\beta^{4}-2\beta^{3}+1=0$} & {\footnotesize{}1.95468312004}\tabularnewline
\hline 
 &  &  &  &  &  & {\footnotesize{}$\cdots$} & \tabularnewline
\hline 
{\footnotesize{}1/4} & {\footnotesize{}5/8} & {\footnotesize{}\negthinspace{}5,0\negthinspace{}} & {\footnotesize{}$101\overline{\cdot}$} & {\footnotesize{}3} & {\footnotesize{}1} & {\footnotesize{}$\beta^{3}-\beta^{2}-1$} & {\footnotesize{}1.46557123187}\tabularnewline
\hline 
{\footnotesize{}1/7} & {\footnotesize{}4/7} & {\footnotesize{}\negthinspace{}0,4\negthinspace{}} & {\footnotesize{}$\overline{100}$} & {\footnotesize{}0} & {\footnotesize{}3} & {\footnotesize{}$\beta^{3}-\beta^{2}-1$} & {\footnotesize{}above}\tabularnewline
\hline 
{\footnotesize{}4/7\negthinspace{}} & {\footnotesize{}\negthinspace{}11/14\negthinspace{}} & {\footnotesize{}\negthinspace{}1,4\negthinspace{}} & {\footnotesize{}$1\overline{100}$} & {\footnotesize{}1} & {\footnotesize{}3} & {\footnotesize{}$\beta^{4}-\beta^{3}-\beta^{2}-\beta+1=0$} & {\footnotesize{}1.72208380573}\tabularnewline
\hline 
{\footnotesize{}\negthinspace{}11/14\negthinspace{}} & {\footnotesize{}\negthinspace{}25/28\negthinspace{}} & {\footnotesize{}\negthinspace{}3,4\negthinspace{}} & {\footnotesize{}\negthinspace{}11$\overline{100}$\negthinspace{}} & {\footnotesize{}2} & {\footnotesize{}3} & {\footnotesize{}$\beta^{5}-\beta^{4}-\beta^{3}-2\beta^{2}+\beta+1=0$\negthinspace{}\negthinspace{}} & {\footnotesize{}1.87134931301}\tabularnewline
\hline 
{\footnotesize{}\negthinspace{}25/28\negthinspace{}} & {\footnotesize{}\negthinspace{}53/56\negthinspace{}} & {\footnotesize{}\negthinspace{}7,4\negthinspace{}} & {\footnotesize{}\negthinspace{}111$\overline{100}$\negthinspace{}} & {\footnotesize{}3} & {\footnotesize{}3} & {\footnotesize{}$\beta^{6}-\beta^{5}-\beta^{4}-2\beta^{3}+\beta+1=0$\negthinspace{}\negthinspace{}} & {\footnotesize{}1.93992448793}\tabularnewline
\hline 
{\footnotesize{}2/7\negthinspace{}} & {\footnotesize{}9/14\negthinspace{}} & {\footnotesize{}\negthinspace{}1,2\negthinspace{}} & {\footnotesize{}$1\overline{010}$} & {\footnotesize{}1} & {\footnotesize{}3} & {\footnotesize{}$\beta^{4}-\beta^{3}-2\beta+1=0$} & {\footnotesize{}1.55897987798}\tabularnewline
\hline 
{\footnotesize{}9/14\negthinspace{}} & {\footnotesize{}\negthinspace{}23/28\negthinspace{}} & {\footnotesize{}\negthinspace{}3,2\negthinspace{}} & {\footnotesize{}\negthinspace{}11$\overline{010}$\negthinspace{}} & {\footnotesize{}2} & {\footnotesize{}3} & {\footnotesize{}$\beta^{5}-\beta^{4}-\beta^{3}-\beta^{2}+1=0$} & {\footnotesize{}1.77847961614}\tabularnewline
\hline 
 &  &  &  &  &  &  & \tabularnewline
\hline 
{\footnotesize{}4/5\negthinspace{}} & {\footnotesize{}9/10} & {\footnotesize{}\negthinspace{}1,12\negthinspace{}} & {\footnotesize{}\negthinspace{}1$\overline{1100}$\negthinspace{}} & {\footnotesize{}1} & {\footnotesize{}4} & {\footnotesize{}$\beta^{5}-\beta^{4}-\beta^{3}-\beta^{2}-\beta+1=0$\negthinspace{}\negthinspace{}} & {\footnotesize{}1.88320350591}\tabularnewline
\hline 
\end{tabular}
\par\end{center}

\endgroup\medskip{}

The rationals in the first column are from the set $\frac{0}{1}\le\frac{p}{q}\le\frac{1}{1}$.
The $\beta$ bit-sequences always have $b_{0}=1$ and thus, the second
column shows $S\left(2\right)=\left(p+q\right)/2q\ge1/2$. The third
column m,c shows the prefix and the cyclic part as integers. The fourth
column shows the actual bitsequence. An overline is drawn over the
repeating digits. If there is one repeating digit, and it is zero,
it is written as $\overline{\cdot}$ so as to keep things a bit more
readable. The prefix must always start with a 1, and so the integer
$m$ is unique. The length of the prefix is in the $L$ column. The
cycle might have leading zeros, and so specifying $c$ is not enough;
a cycle length is required, given in the $N$ column. The corresponding
polynomial and it's root are given in the last two columns.

Notes:
\begin{itemize}
\item The sequence $1\overline{10}$ is described in OEIS A160389.
\item The sequence $10\overline{10}$ is obviously not allowed, as it is
the same as $\overline{10}$.
\item The sequence $100\overline{10}=10\overline{01}$ is not allowed. In
the unreduced form, replacing $\overline{10}$ by $11\overline{0}$
gives $100\overline{10}\mapsto10011=19,$ which wasn't allowed as
a finite orbit.
\item The sequence $101\overline{10}$ is not allowed; it reduces to a finite
form $101\overline{10}\mapsto10111=23$ that is a disallowed finite
orbit.
\item These last two observations are a coincidence, and do not hold in
the general case. There are good periodic orbits that have disallowed
finite versions, and vice-versa.
\item The sequence $10\overline{100}$ is not allowed, because it is reducible:
$10\overline{100}=1\overline{010}$.
\item 4/7 gives $1\overline{100}$ seems to be OEIS A289917 which is $\left(1+\sqrt{13}+\sqrt{2\sqrt{13}-2}\right)/4$.
\item 4/5 gives 1$\overline{1100}$ seems to be OEIS A289915 which is $\left(1+\sqrt{2}+\sqrt{2\sqrt{2}-1}\right)/2$.
\item Rationals with prefix $m=0$ correspond to the finite orbits.
\item It seems the prefix $m=2^{L-1}$ never occurs, when $L>1$.
\item The prefix is often but not always odd. Counterexamples include $7/12=110\overline{01}$
and $13/20=110\overline{1001}$.
\item The first time that $m=5$ occurs is for $9/28=101\overline{010}$.
\end{itemize}
Not all rationals give valid bit-sequences. The set of valid rationals
are exactly the ones that are in the set $\mathbb{Q}\cap\overline{\theta}$
where $\overline{\theta}$ is the infinite comb, given in eqn \ref{eq:infinite comb},
and depicted visually in figure \ref{fig:Theta-Indicator}. The comb
is extremely fractal, and working directly with $\mathbb{Q}\cap\overline{\theta}$
would be a chore. Fortunately, the bracket map saves the day: it is
a map of the entire unit interval onto $\overline{\theta}$. The bracket
map, shown in figure \ref{fig:Beta-Bracket-Map}, maps rationals to
rationals; it provides a map $\mathbb{Q}\cap\left[0,1\right]\to\mathbb{Q}\cap\overline{\theta}$.
Thus, it provides a better, though more indirect way, of describing
the ultimately-periodic orbits. Some examples are reviewed in the
section after next.

\subsection{Rational to periodic binary}

As a practical matter, it is computationally useful to convert a given
fraction $p/q$ into the binary prefix $m$, it's length $L$, and
the cyclic part $c$, and its length $N$. This is not hard, but also
not entirely easy, and so is presented here. Write
\begin{align*}
\frac{p+q}{2q} & =\frac{m}{2^{L}}+\frac{1}{2^{L}}\left(\frac{c}{2^{N}}+\frac{c}{2^{2N}}+\frac{c}{2^{3N}}+\cdots\right)\\
 & \frac{1}{2^{L}}\left(m+\frac{c}{2^{N}}\left(1+\frac{1}{2^{N}}+\frac{1}{2^{2N}}+\cdots\right)\right)\\
 & \frac{1}{2^{L}}\left(m+\frac{c}{2^{N}-1}\right)\\
 & \frac{1}{2^{L}}\cdot\frac{m\left(2^{N}-1\right)+c}{2^{N}-1}
\end{align*}
In the last line, both numerator and denominator are integers. This
allows the following algorithm:
\begin{enumerate}
\item Find $\gcd\left(p+q,2q\right)$ and so reduce to lowest terms $a/b=\left(p+q\right)/2q$.
\item Factor $b=2^{L}b^{\prime}$ to obtain $L\ge0$.
\item Solve $\left(2^{N}-1\right)\mod b^{\prime}=0$ for the smallest positive
$N\ge1$.
\item Define $r=\left(2^{N}-1\right)a/b^{\prime}$. 
\item Solve $\left(r-c\right)\mod\left(2^{N}-1\right)=0$ for the smallest
$c$ such that $0\le c<2^{N}$.
\item Define $m=\left(r-c\right)/\left(2^{N}-1\right)$.
\end{enumerate}
This provides all four integers $m,L,c,N$ that define the cyclic
expansion.

\subsection{The Bracket Map}

The bracketing relationship provided a mapping to the dyadic rationals,
via moves on the binary tree. What do ultimately-periodic moves correspond
to, on that tree? 

The left and right moves $L,R$ on the valid-index binary tree were
$L:\left(\ell\Mapsto f\Mapsfrom\rho\right)\mapsto\left(\ell\Mapsto2f\Mapsfrom f\right)$
and $R:\left(\ell\Mapsto f\Mapsfrom\rho\right)\mapsto\left(f\Mapsto\Lambda\left(f\right)\Mapsfrom\rho\right)$.
(The symbol $R$ is used in place of $\mathfrak{R}$ in this section,
for ease of typography. They are isomorphic.) Consider starting at
$\infty\Mapsto1\Mapsfrom0$ and applying a sequence of alternating
left and right moves. This is shown in the table below:
\begin{center}
\medskip{}
\begingroup
\def\arraystretch{1.3}%
\begin{tabular}{|c|c|c|c|}
\hline 
position & move & bracket & $\beta$\tabularnewline
\hline 
\hline 
1/2 & - & $\infty\Mapsto1\Mapsfrom0$ & 1.61803398...\tabularnewline
\hline 
1/4 & $L$ & $\infty\Mapsto2\Mapsfrom1$ & 1.46557123...\tabularnewline
\hline 
3/8 & $LR$ & $2\Mapsto10\Mapsfrom1$ & 1.57014731...\tabularnewline
\hline 
5/16 & $LRL$ & $2\Mapsto20\Mapsfrom10$ & 1.53849659...\tabularnewline
\hline 
11/32 & $LRLR$ & $20\Mapsto82\Mapsfrom10$ & 1.56175206...\tabularnewline
\hline 
21/64 & $LRLRL$ & $20\Mapsto164\Mapsfrom82$ & 1.55392112...\tabularnewline
\hline 
43/128 & $LRLRLR$ & $164\Mapsto658\Mapsfrom82$ & 1.55970265...\tabularnewline
\hline 
85/256 & $LRLRLRL$ & $164\Mapsto1316\Mapsfrom658$ & 1.55767530...\tabularnewline
\hline 
171/512 & $LRLRLRLR$ & $1316\Mapsto5266\Mapsfrom658$ & 1.55917021...\tabularnewline
\hline 
1/3 & $\overline{LR}$ & — & 1.55897978...\tabularnewline
\hline 
\end{tabular}
\par\end{center}

\endgroup\medskip{}

The column labeled ``position'' is the location in the dyadic tree.
By convention, 1/2 is at the top, and 1/4 lies to the left, and 3/4
lies to the right. The column labeled ``move'' consists of the sequence
of left-right moves to get to a given tree position. For $n$ moves,
encoded as a binary integer $m$, the corresponding position in the
dyadic tree is $\left(2m+1\right)/2^{n}$. The bracket is likewise
a position in the valid-index binary tree. The string of moves to
get to that location are written in reverse applicative order, so
that the first letter in the string is the first move. If the moves
are treated as functions to be composed and applied, then the string
needs to be reversed to get the applicative order. The beta is the
beta value at the center of that bracket. Based on the beta value,
we conclude that $\mathsf{G}:1/3\mapsto2/7$.

A second example helps cement the idea.

\medskip{}

\begingroup
\def\arraystretch{1.3}
\begin{center}
\begin{tabular}{|c|c|c|c|}
\hline 
position & move & bracket & $\beta$\tabularnewline
\hline 
\hline 
1/2 & - & $\infty\Mapsto1\Mapsfrom0$ & 1.61803398...\tabularnewline
\hline 
1/4 & $L$ & $\infty\Mapsto2\Mapsfrom1$ & 1.46557123...\tabularnewline
\hline 
1/8 & $LL$ & $\infty\Mapsto4\Mapsfrom2$ & 1.38027756...\tabularnewline
\hline 
3/16 & $LLR$ & $4\Mapsto36\Mapsfrom2$ & 1.44326879...\tabularnewline
\hline 
5/32 & $LLRL$ & $4\Mapsto72\Mapsfrom36$ & 1.42705896...\tabularnewline
\hline 
11/64 & $LLRLR$ & $72\Mapsto580\Mapsfrom36$ & 1.43949911...\tabularnewline
\hline 
21/128 & $LLRLRL$ & $72\Mapsto1160\Mapsfrom580$ & 1.43591015...\tabularnewline
\hline 
43/256 & $LLRLRLR$ & $1160\Mapsto9284\Mapsfrom580$ & 1.43866733...\tabularnewline
\hline 
85/512 & $LLRLRLRL$ & $1160\Mapsto18568\Mapsfrom9284$ & 1.43784133...\tabularnewline
\hline 
1/6 & $L\overline{LR}$ & — & 1.43841656...\tabularnewline
\hline 
\end{tabular}
\par\end{center}

\endgroup\medskip{}

Based on the convergent, we conclude that $\mathsf{G}:1/6\mapsto2/15$,
deduced below. In this case, $1/6\notin\overline{\theta}$ and so
the raw fraction $1/6$ does not generate a self-describing orbit.
But $2/15$ is self-describing. Of course, that is the entire intent
of the bracket map. Originally formulated as the good-index map $G:\mathbb{N}\to\Psi$
which maps natural numbers to brackets, the extension $\mathsf{G}$
is the ``good rational'' map $\mathsf{G}:\mathbb{Q}\to\mathbb{Q}\cap\overline{\theta}$,
which maps rationals to those that are self-describing. This is the
same map as visualized in figure \ref{fig:Beta-Bracket-Map}.

\subsection{Bracketing Examples}

A curated collection of examples of rationals passed through the bracketing
function is given below. It serves mostly to give a sense of the patterns
that develop, as well as counterexamples that defy simple patterns.

\medskip{}
\begingroup
\def\arraystretch{1.3}
\begin{center}
{\small{}}%
\begin{longtable}[c]{|c|c|c|c|c|c|c|}
\hline 
{\small{}$\frac{a}{b}$} & {\small{}$\Psi$} & {\small{}moves} & {\small{}$\frac{p}{q}=\mathsf{G}\frac{a}{b}$} & {\small{}id} & {\small{}orbit} & {\small{}$\beta$}\tabularnewline
\hline 
\hline
\endhead
\hline 
{\small{}0/1} & {\small{}-1} & {\small{}$\epsilon$} & {\small{}0/1} & {\small{}y} & {\small{}$1\overline{\cdot}$} & {\small{}1}\tabularnewline
\hline 
{\small{}$1/2^{n}-\varepsilon$} &  & {\small{}$\approx RL\cdots L\overline{R}$} & {\small{}$1/\left(2^{n+1}-1\right)-\delta$} & {\small{}-} & {\small{}$\approx1\overline{0\cdots01}$} & \tabularnewline
\hline 
{\small{}$1/2^{n}$} & {\small{}$2^{n-1}$} & {\small{}$RL\cdots L\overline{R}$} & {\small{}$1/2^{n}$} & {\small{}y} & {\small{}$1\overline{0\cdots01}$} & \tabularnewline
\hline 
{\small{}$1/8-\varepsilon$} &  & {\small{}$\approx RLLL\overline{R}$} & {\small{}$1/15-\delta$} & {\small{}-} & {\small{}$\approx1\overline{0001}$} & {\small{}$<$1.380277}\tabularnewline
\hline 
{\small{}1/8} & {\small{}4} & {\small{}$RLL$} & {\small{}1/8} & {\small{}y} & {\small{}$1001$} & {\small{}1.38027756}\tabularnewline
\hline 
{\small{}1/7} &  & {\small{}$R\overline{LLR}$} & {\small{}4/31} & {\small{}-} & {\small{}$1\overline{00100}$} & {\small{}1.42109608}\tabularnewline
\hline 
{\small{}1/6} &  & {\small{}$RL\overline{LR}$} & {\small{}2/15} & {\small{}-} & {\small{}$1\overline{0010}$} & {\small{}1.43841656}\tabularnewline
\hline 
{\small{}1/5} &  & {\small{}$R\overline{LLRR}$} & {\small{}12/85} & {\small{}-} & {\small{}$1\overline{00100100}$} & {\small{}1.45394278}\tabularnewline
\hline 
{\small{}$1/4-\varepsilon$} &  & {\small{}$\approx RLL\overline{R}$} & {\small{}$1/7-\delta$} & {\small{}-} & {\small{}$\approx1\overline{001}$} & {\small{}$<$1.4655712}\tabularnewline
\hline 
{\small{}1/4} & {\small{}2} & {\small{}$RL$} & {\small{}1/4} & {\small{}y} & {\small{}$101\overline{\cdot}$} & {\small{}1.46557123}\tabularnewline
\hline 
{\small{}2/7} &  & {\small{}$R\overline{LRL}$} & {\small{}4/15} & {\small{}-} & {\small{}$1\overline{0100}$} & {\small{}1.52626195}\tabularnewline
\hline 
{\small{}1/3} &  & {\small{}$R\overline{LR}$} & {\small{}2/7} & {\small{}-} & {\small{}$1\overline{010}$} & {\small{}1.55897987}\tabularnewline
\hline 
{\small{}2/5} &  & {\small{}$R\overline{LRRL}$} & {\small{}20/63} & {\small{}-} & {\small{}$1\overline{010100}$} & {\small{}1.58925391}\tabularnewline
\hline 
{\small{}3/7} &  & {\small{}$R\overline{LRR}$} & {\small{}10/31} & {\small{}-} & {\small{}$1\overline{01010}$} & {\small{}1.60022189}\tabularnewline
\hline 
{\small{}$1/2-\varepsilon$} &  & {\small{}$\approx RL\overline{R}$} & {\small{}$1/3-\delta$} & {\small{}-} & {\small{}$\approx1\overline{01}$} & {\small{}$<1.618\cdots$}\tabularnewline
\hline 
{\small{}1/2} & {\small{}1} & {\small{}$R$} & {\small{}1/2} & {\small{}y} & {\small{}$11\overline{\cdot}$} & {\small{}$\varphi=1.618\cdots$}\tabularnewline
\hline 
{\small{}6/11} &  &  & {\small{}8900/16383} & {\small{}-} &  & {\small{}1.69971346}\tabularnewline
\hline 
{\small{}5/9} &  & {\small{}$R\overline{RLLLRR}$} & {\small{}40/73} & {\small{}-} &  & {\small{}1.70348856}\tabularnewline
\hline 
{\small{}4/7} &  & {\small{}$R\overline{RLL}$} & {\small{}4/7} & {\small{}y} & {\small{}$1\overline{100}$} & {\small{}1.72208380}\tabularnewline
\hline 
{\small{}3/5} &  & {\small{}$R\overline{RLLR}$} & {\small{}25/42} & {\small{}-} & {\small{}$11\overline{001100}$} & {\small{}1.74720863}\tabularnewline
\hline 
{\small{}9/14} &  & {\small{}$RR\overline{LRL}$} & {\small{}9/14} & {\small{}y} & {\small{}11$\overline{010}$} & {\small{}1.77847961}\tabularnewline
\hline 
{\small{}2/3} &  & {\small{}$R\overline{RL}$} & {\small{}2/3} & {\small{}y} & {\small{}$1\overline{10}$} & {\small{}1.80193773}\tabularnewline
\hline 
{\small{}7/10} &  & {\small{}$RR\overline{LRRL}$} & {\small{}43/62} & {\small{}-} &  & {\small{}1.82000973}\tabularnewline
\hline 
{\small{}5/7} &  & {\small{}$R\overline{RLR}$} & {\small{}7/10} & {\small{}-} & {\small{}$11\overline{0110}$} & {\small{}1.82651577}\tabularnewline
\hline 
{\small{}3/4} & {\small{}3} & {\small{}$RR$} & {\small{}3/4} & {\small{}y} & {\small{}$111\overline{\cdot}$} & {\small{}1.83928675}\tabularnewline
\hline 
{\small{}7/9} &  & {\small{}$R\overline{RRLLLR}$} & {\small{}227/292} & {\small{}-} &  & {\small{}1.86625406}\tabularnewline
\hline 
{\small{}11/14} &  & {\small{}$RR\overline{RLL}$} & {\small{}11/14} & {\small{}y} & {\small{}11$\overline{100}$} & {\small{}1.87134931}\tabularnewline
\hline 
{\small{}4/5} &  & {\small{}$R\overline{RRLL}$} & {\small{}4/5} & {\small{}y} & {\small{}$1\overline{1100}$} & {\small{}1.88320350}\tabularnewline
\hline 
{\small{}5/6} &  & {\small{}$RR\overline{RL}$} & {\small{}5/6} & {\small{}y} & {\small{}$11\overline{10}$} & {\small{}1.90516616}\tabularnewline
\hline 
{\small{}6/7} &  & {\small{}$R\overline{RRL}$} & {\small{}6/7} & {\small{}y} & {\small{}$1\overline{110}$} & {\small{}1.92128960}\tabularnewline
\hline 
{\small{}8/9} &  & {\small{}$R\overline{RRRLLL}$} & {\small{}8/9} & {\small{}y} & {\small{}$1\overline{111000}$} & {\small{}1.93762945}\tabularnewline
\hline 
{\small{}25/28} &  & {\small{}$RRR\overline{RLL}$} & {\small{}25/28} & {\small{}y} & {\small{}111$\overline{100}$} & {\small{}1.93992448}\tabularnewline
\hline 
{\small{}9/10} &  & {\small{}$RR\overline{RRLL}$} & {\small{}9/10} & {\small{}y} & {\small{}$11\overline{1100}$} & {\small{}1.94470383}\tabularnewline
\hline 
{\small{}10/11} &  &  & {\small{}10/11} & {\small{}y} &  & {\small{}1.94988340}\tabularnewline
\hline 
{\small{}11/12} &  & {\small{}$RRR\overline{RL}$} & {\small{}11/12} & {\small{}y} & {\small{}$111\overline{10}$} & {\small{}1.95468312}\tabularnewline
\hline 
{\small{}1/1} & {\small{}0} & {\small{}$\overline{R}$} & {\small{}1/1} & {\small{}y} & {\small{}$\overline{1}$} & {\small{}2}\tabularnewline
\hline 
\end{longtable}{\small\par}
\par\end{center}

\endgroup\medskip{}

Table legend.
\begin{itemize}
\item The first column shows selected interesting rationals. Every possible
rational is allowed here, and will generate an ultimately-periodic
orbit. 
\item The second column shows the corresponding index $\Psi$; only seven
are shown, for the finite orbits that correspond to the dyadic rationals;
all others have an infinite limit. 
\item The third column shows the bracket moves generated by that rational.
These are obtained as the binary expansion of the fraction $a/b$,
or, more properly, the expansion for $\left(a+b\right)/2b$. Thus,
the first move is always $R$, to arrive from the far left to the
center of the tree. The expansions for 6/11 and 10/11 were left blank;
these are pointlessly long strings that would have cluttered the table.
\item The fourth column shows the result of the ``good rational'' map
$\mathsf{G}:\mathbb{Q}\to\mathbb{Q}\cap\overline{\theta}$ that takes
$\mathsf{G}:\frac{a}{b}\mapsto\frac{p}{q}$. Thus, by definition,
$p/q\in\overline{\theta}$. The fractions are those obtained by extending
the ``good index'' map $G:\mathbb{N}\to\Psi$ to it's limit points.
A graph of the fourth column, relative to the first, is given in figure
\ref{fig:Good-Map}. Although that figure was prepared for the finite
orbits, it is exactly the same for the periodic orbits.
\item Some of the discontinuities visible in \ref{fig:Good-Map} are shown
here, in the rows with $\delta,\epsilon\ll1$ in them. For example,
the top of the tree is 1/2, which has a finite orbit $1\overline{0}$.
It can also be written as an infinite cyclic orbit $1\overline{01}=1/3$
and these two distinct binary strings are mapped to distinct fractions;
thus, the discontinuities. Such a discontinuity appears at every finite
orbit. 
\item Notable rows are 6/11, which has a horribly large denominator; and
7/9 is the runner-up. The large denominator is due to leader heights
that are greater than one, that are encountered during the expansion.
Compare rows 7/9 to 8/9. Both have long move strings, but $\mathsf{G}\left(7/9\right)=7/9$.
Similarly, 6/11 and 10/11 have intolerably long move strings, but
$\mathsf{G}\left(10/11\right)=10/11$.
\item The fifth column has a 'y' in it, whenever $\mathsf{G}\left(a/b\right)=a/b$;
when the first and fifth column are equal. Note that most 'y''s occur
at the end of the table, rather than at the beginning.
\item The sixth column shows the bitsequence $b_{0}b_{1}\cdots$ for the
binary expansion of $p/q$ . Some rows are blank; these rows had horribly
long expansions that provide no intuition.
\item The seventh and rightmost column shows the corresponding $\beta$
value for the orbit in the previous column.
\end{itemize}
Perhaps most notable is that the good map $\mathsf{G}$ maps rationals
to rationals. Reassuringly, it extends from the dyadic rationals,
to all rationals. Formal lemmas and proofs follow in a later section.

\pagebreak{}

\section{Formalities}

Theorems, lemmas, proofs. A number of observations, claims and assumptions
have been made. They seem ``obviously true'' from numerical work,
but are lacking a formal proof. These are sketched below. As sketches,
this chapter is still quite informal; the difference is that it uses
a more sophisticated vocabulary to say things in a more abstract manner.
\begin{itemize}
\item Theorem: The finite orbits are dense in $1\le\beta\le2$. Proof: Provided
by the bracketing relation, which states that between any two endpoints,
there's another orbit strictly in the middle. We even have more: the
estimate $\nu^{1/\nu}$ from the limit diagram, showing how they don't
want to accumulate near $\beta=1$.
\item Lemma: The bracketed roots are always in strict ascending order, with
$r_{\ell}<r_{f}<r_{\rho}$.
\item Lemma: The valid-index map $G$ is a bijection. Proof: Leadership
function is monotonically increasing.
\item Theorem: The good map $\mathsf{G}$ always maps rationals to rationals.
More precisely, those rationals are always inside of $\overline{\theta}$
so that $\mathsf{G}:\mathbb{Q}_{I}\to\overline{\theta}\cap\mathbb{Q}$
is a bijection, except at the dyadics, where we have a countable freedom
of finite orbits and one infinite-periodic orbit to choose from. 
\item Theorem: $\mathsf{G}:\left[0,1\right]\to\overline{\theta}$ is a bijection
for all reals $\left[0,1\right]\subset\mathbb{R}$. That is, the closure
to reals works as expected. This follows because the roots form a
countable dense subset, the reals are separable, and the function
is continuous.
\item Proper diligence requires distinguishing $\overline{\mathfrak{B}}_{\mathbb{Q}}$
from $\overline{\mathfrak{B}}_{\mathbb{R}}$ and also $\overline{\theta}_{\mathbb{Q}}$
from $\overline{\theta}_{\mathbb{R}}$. 
\end{itemize}
The remainder of this section consists of fragments of proofs of the
above theorems and lemmas. If these were stitched together properly,
whole proofs should emerge. As it is, its more an outline or sketch. 

\subsection{The Rational Comb}

The dyadic comb, of figure \ref{fig:Theta-Indicator}, shows where
valid finite orbits can occur. It's limit to infinite rank gives the
infinite comb $\theta$ of eqn \ref{eq:infinite comb}. The only valid
rationals, that is, those giving self-describing orbits, are those
that belong to $\mathbb{Q}\cap\theta$. 

This can be understood as a limiting procedure. Every rational can
be approximated by a dyadic rational; each set is dense in the other.
The dyadic rationals are finite walks down the binary tree; the non-dyadic
rationals are infinite walks. The infinite comb corresponds to the
trimmed tree, described earlier. The trimming always maintains branches
of unbounded length. These have a limit; in the limit, some of these
will correspond to rationals. These are exactly the rationals that
correspond to the ultimately-periodic orbits. The correspondence goes
in both directions: given an infinite path down the tree, it is sufficient
to truncate it to be of finite length. By construction, the truncated
path corresponds to a valid dyadic rational. To recap: going in one
direction, there is a sequence of finite orbits that converge, in
the limit, onto an eventually-periodic orbit. Conversely, given such
an orbit, every truncated version thereof is valid. 

\subsection{Closures}

The ultimately-periodic orbits live inside the closure $\overline{\mathfrak{B}}\subset\overline{\mathbb{B}}$.
Write $\mathfrak{C}\subset\overline{\mathbb{B}}$ for the cyclic (ultimately
periodic) orbits; these are infinite-length strings describing moves
down the infinite binary tree. As already noted, $\mathfrak{C}\cap\mathbb{B}\subset\mathfrak{B}$;
that is, every infinite-length cyclic orbit, when truncated to finite
length, is a valid finite orbit. The closure $\overline{\mathfrak{B}}$
consists of all self-describing bit-strings obtained as solutions
to $\beta=S\left(\beta\right)$ where $S\left(\beta\right)=\sum_{k=0}^{\infty}b_{k}\beta^{-k}$.
The cyclic orbits are just a special case: $\mathfrak{C}\subset\overline{\mathfrak{B}}$.
There are presumably many more, uncountably many chaotic orbits in
$\overline{\mathfrak{B}}$. 

Just a little bit more machinery is needed. Let $2^{\omega}$ be Cantor
space, the space of infinitely-long binary strings. Let $\chi:2^{\omega}\to\left[0,1\right]$
be the canonical mapping of Cantor space to the unit interval: $\chi:\left(b_{0}b_{1}\cdots\right)\mapsto\sum_{k}b_{k}2^{-k-1}$.
Note that this mapping manages to miss all the dyadic rationals, as
these correspond to finite strings, of which there aren't any in the
Cantor space. Let $2^{<\omega}=2^{*}$ be the set of finite-length
binary strings, and allow the symbol $\chi$ to perform double-duty,
by writing $\chi:2^{*}\to\mathbb{D}$ to map finite strings to dyadics.
These binary strings can also be interpreted as strings in the symbols
$L,R$, so that they are moves on the tree. Write $\iota:2^{*}\to\mathbb{B}$
that maps the empty string to the root of the tree, and likewise $\iota:2^{\omega}\to\overline{\mathbb{B}}$
that maps infinite strings to the leaves of the infinite tree.

The comb was a subset of the unit interval: $\overline{\theta}\subset\left[0,1\right]$.
The claim is that $\overline{\theta}=\chi\iota\overline{\mathfrak{B}}$.
The mapping is onto, but not one-to-one: The finite orbits all have
corresponding cyclic orbits, all of which are distinct elements in
$\overline{\mathfrak{B}}$ but map to the (finite-orbit) polynomial,
and thus the same root.

The good-index function provided a bijection $\mathfrak{B}=\eta G\eta^{-1}\mathbb{B}$
between the trimmed and untrimmed finite, unbounded trees. The claim
is that this can be extended to a bijection $\overline{\mathfrak{B}}=\overline{G}\overline{\mathbb{B}}$.
This requires a short detour into topology. The standard (weak) topology
on Cantor space is the product topology. The base for the topology
may be taken as the set of all finite strings, followed by an infinite
number of don't-care markers. These are the open sets; the full topology
is the finite intersection and infinite union of the open sets in
this base.

Every set in the base of the topology is represented by a finite string.
The function $\eta^{-1}$ maps this to an integer. The function $G$
is defined on all integers; it returns an integer, which is mapped
by $\eta$ to a finite string and thus an open set. To conclude, $\eta G\eta^{-1}$
maps open sets to open sets, and it is defined on \emph{every} open
set in the topology. Thus, it is safe to write $\overline{G}=\eta G\eta^{-1}$,
as it is defined everywhere; all of $\overline{\mathbb{B}}$ is in
it's domain. The range of $\overline{G}$ can be taken to be the definition
of the closure $\overline{\mathfrak{B}}$. This has several benefits:
it avoids having to make funny arguments that pass through the reals
via the comb, and it also makes clear that $\overline{G}$ is a bijection.
Note $\overline{G}$ is not continuous, in the context that $\overline{\mathfrak{B}}\subset\overline{\mathbb{B}}$,
since $\overline{G}^{-1}$ is defined almost nowhere on $\overline{\mathbb{B}}$.

The points $r_{n}$ are dense in the interval $\left[1,2\right]$.
Within each rank $\nu$ they are monotonically increasing. The bracketing
relation guarantees that, for all $r_{n}$ in the subtree, these are
all contained within the endpoints of the bracket, and that they are
always totally ordered, by the natural tree-ordering. The notion of
open sets is compatible. Given a basic open set in $\overline{\mathbb{B}}$,
the function $\rho=r\circ G\circ\eta^{-1}$ maps it to an open set
$U\subset\left[1,2\right]$. Claim that the inverse map also maps
open sets to open sets, and thus $\varrho$ is continuous. This is
established as follows. Since the $r_{n}$ are dense in the interval
$\left[1,2\right]$, any open set $U\subset\left[1,2\right]$ can
be written as a countable union of brackets. The function $\rho^{-1}$
is defined on all brackets, and by bracketing, it is defined on all
$U\subset\left[1,2\right]$ and by bracketing $\rho^{-1}U$ is an
open set in the product topology on $\overline{\mathbb{B}}$.

The goal is to extend this reasoning to $\varrho:\mathbb{D}\to\left[1,2\right]$,
given by $\varrho=\rho\circ\eta\circ\delta^{-1}=r\circ G\circ\delta^{-1}$
so as to obtain a continuous, monotonic function $\overline{\varrho}:\left[0,1\right]\to\left[1,2\right]$
on the unit interval. And we are done, more or less.

The mapping $\chi$ maps the basic open sets of $\overline{\mathbb{B}}$
to open subintervals of the unit interval; more precisely, to subintervals
with dyadics at each end. These are precisely the sets $I\left(m,\nu\right)$
defined in eqn \ref{eq:topo-base}. They allowed the infinite comb
to be built up from unbounded-length but finite strings in $\mathfrak{B}$.

The base of the product topology was mapped to open intervals on the
real-number line; and so $\overline{\varrho}$ is continuous on the
reals. This can also be seen in a different way. The standard measure
$\mu$ on the reals says that $\mu I\left(m,\nu\right)=2^{2-\nu}$.
The function $G$ is built from left and right moves. The left moves
bump the rank by one, and so always map to sets that are exactly half
the size on the real number line. The right moves are given by the
leadership function: $\Lambda\Psi_{\nu}\subset\bigcup_{h=0}^{\infty}\Psi_{\nu+1+h}$
for each rank $\nu$. The resulting sets are (at a minimum) half the
size, and possibly larger, as the union of additional small intervals.
This allows a conventional delta-epsilon proof of continuity to go
forward: for each epsilon, one can choose a rank $\nu$ where the
basis sets $I\left(m,\nu\right)$ are smaller than epsilon. The function
$G$ maps them to other open sets that are strictly smaller. Strictly,
because the height of a leader is always a finite number; the union
is written as a union over all possible heights, but the union is
always finite.

\subsection{Proof that roots are correctly bounded}

Theorem: The bracket relationship gives roots with $r_{\ell}<r_{f}<r_{\rho}$
being strict inequalities.

Proof: Recursion starts with $\infty\Mapsto1\Mapsfrom0$. The bracket
moves are 
\begin{align*}
L & :\left(\ell\Mapsto f\Mapsfrom\rho\right)\mapsto\left(\ell\Mapsto2f\Mapsfrom f\right)\\
\mathfrak{R} & :\left(\ell\Mapsto f\Mapsfrom\rho\right)\mapsto\left(f\Mapsto\Lambda\left(f\right)\Mapsfrom\rho\right)
\end{align*}
So $p_{1}=\beta^{2}-\beta-1$ and $p_{\infty}=p_{1/0}=\beta-1$ and
$p_{1/\infty}=\beta-2$. Must show that bracket bounds are respected
in all four cases. The proofs depend on having $p_{n}^{\prime}\left(\beta\right)>0$
for all $1\le\beta\le2$.
\begin{description}
\item [{Case}] 1: Show that $r_{2f}<r_{f}$. Proof: The $L$ move gives
$p_{2f}=\beta\left(p_{f}+1\right)-1$. So $1=r_{2f}\left(p_{f}\left(r_{2f}\right)+1\right)$
but since $r_{2f}>1$, must have $p_{f}\left(r_{2f}\right)<0$ and
since $p_{n}^{\prime}>0$, conclude that $r_{2f}<r_{f}$.
\item [{Case}] 2: Show that $r_{\ell}<r_{2f}$. Proof: For any $\ell\Mapsto f\Mapsfrom\rho$
it is always the case that $\rho$ is in the tree rooted by $\ell$.
This can be seen by examining the two possibilities arising from $\infty\Mapsto1\Mapsfrom0$
and then keeping in mind all trees are isomorphic. The path to $\rho$
from it's root $\ell$ is $\rho=L^{n}R\ell$ for some $n\ge0$. Thus
$p_{f}=\beta^{n+1}p_{\ell}-1$ and so $p_{2f}=\beta^{n+2}p_{l}-1$.
But $p_{\ell}\left(r_{\ell}\right)=0$ and so $p_{2f}\left(r_{\ell}\right)=-1<0$.
Since $p_{n}^{\prime}\left(\beta\right)>0$, conclude that $r_{\ell}<r_{2f}$. 
\item [{Case}] 3: Show that $r_{f}<r_{\Lambda\left(f\right)}$. Proof:
Proceeds as a modified proof of case 1. For height zero, $\Lambda\left(f\right)=2f+1$
so $p_{\Lambda\left(f\right)}=\beta p_{f}-1$. Thus, $1=r_{\Lambda\left(f\right)}p_{f}\left(r_{\Lambda\left(f\right)}\right)$
so $0<p_{f}\left(r_{\Lambda\left(f\right)}\right)<1$ and since $p_{n}^{\prime}>0$,
conclude that $r_{f}<r_{\Lambda\left(f\right)}$. What about positive
heights? For height $h$, the moves are $L^{h}\mathfrak{R}$ and so
$p_{\Lambda\left(f\right)}=\beta^{h+1}p_{f}-1$. Again, $1=\left(r_{\Lambda\left(f\right)}\right)^{h+1}p_{f}\left(r_{\Lambda\left(f\right)}\right)$
and again $0<p_{f}\left(r_{\Lambda\left(f\right)}\right)$ and so
again $r_{f}<r_{\Lambda\left(f\right)}$.
\item [{Case}] 4: Show that $r_{\Lambda\left(f\right)}<r_{\rho}$. Proof:
This proceeds as a modified version of case 2, with the role of left
and right reversed. In this case, $\rho$ roots a tree, and the path
to $f$ is $f=R^{n}L\rho$ for some $n\ge0$. For $n=0$, one has
$p_{f}=\beta\left(p_{\rho}+1\right)-1$ and so $p_{\Lambda\left(f\right)}=\beta^{h+1}\left(p_{f}\right)-1$.
Expanding, $p_{\Lambda\left(f\right)}=\beta^{h+2}p_{\rho}+\beta^{h+2}\left(\beta-1\right)-1$,
and since $p_{\rho}\left(r_{\rho}\right)=0$, one concludes that $p_{\Lambda\left(f\right)}\left(r_{\rho}\right)>0$
and therefore $r_{\Lambda\left(f\right)}<r_{\rho}$. For $n>0$, the
proof proceeds the same way, with $p_{\Lambda\left(f\right)}=ap_{\rho}+b$
for some positive numbers $a,b$. Thus, $p_{\Lambda\left(f\right)}\left(r_{\rho}\right)=b>0$
and again $r_{\Lambda\left(f\right)}<r_{\rho}$.
\end{description}
To conclude: the roots are ordered as $r_{\ell}<r_{f}<r_{\rho}$ with
strict inequalities holding. QED.

\pagebreak{}

\section{Islands of Stability as Arnold Tongues\label{subsec:Islands-of-Stability}}

The classical Feigenbaum bifurcation diagram, shown in figure \ref{fig:Logistic-Map-Bifurcation},
manifests two distinct behaviors: the ``islands of stability'',
in which there are periodic orbits, and the ``chaotic regions''.
These regions alternate and interleave as a parameter $\lambda$ appearing
in the iterated equation $\lambda x\left(1-x\right)$ is swept through
a range of values. By contrast, the equivalent diagram for the beta-map,
shown in figure \ref{fig:Undershift-Bifurcation-Diagram}, does not
seem to have regions of stability. This is only an illusion: they
are there, they are only infinitely thin. This chapter focuses on
how to crowbar them open, to finite size.

In the previous chapters, it was demonstrated that there is a countable
set of $\beta$ values, dense in the range $1\le\beta\le2$, for which
orbits are finite and terminate after a fixed number of iterations.
Alternately, they can be made periodic, simply by changing a less-than
sign to a less-than-or-equals sign. The $\beta$ values for which
orbits are periodic can be placed in correspondence with the periodic
orbits of the logistic map; the $\beta$ values for which orbits are
chaotic correspond to the chaotic orbits of the logistic map. The
open problem is to demonstrate this correspondence explicitly. This
problem is not tackled here; it is just brushed up against.

The circle map $x_{n+1}=x_{n}+\theta+K\sin2\pi x_{n}$ provides one
possible mechanism for taking a set of measure zero, and crowbaring
it open to a set of finite size. For $K=0$, this is just the rotation
map $x_{n+1}=x_{n}+\theta$ which has only periodic orbits, when $\theta$
is rational, and chaotic orbits, when $\theta$ is not rational. As
one sweeps $\theta$ through a range, the subset of periodic orbits
is countable, and is a set of measure zero: thus, the rotation map
can be said to be chaotic for almost all $\theta$. Setting $K$ to
a non-zero value expands the regions of periodic orbits to finite
size, termed Arnold tongues. These are the mode-locking regions that
are generically visible in driven oscillator systems. The perturbation
by the kick $K$ displaces what would have been chaotic orbits into
mode-locked regions. For small $K$, this perturbation is soft, in
that one might say ``it shouldn't change things much.'' But even
a whisper of a miniscule $K$ is enough to convert the set of periodic
orbits from a measure of zero to a finite measure.

Another possibility is to just crowbar open the periodic regions with
a ``hard'' perturbation, localized at a point. Take the natural
saw-tooth shape of the $\beta$-map, widen the middle, and insert
a slanting downward line, to create a zig-zag. That is, connect the
two endpoints in the middle of the beta shift, ``widening'' it so
that it has a finite, not infinite slope, thereby converting the iterated
function from a discontinuous to a continuous one. This can be constructed
directly: given some ``small'', real $\varepsilon>0$, define the
piecewise-linear $\varepsilon$-generalization of the map \ref{eq:downshift}
as 
\begin{equation}
T_{\beta,\varepsilon}(x)=\begin{cases}
\beta x & \mbox{ for }0\le x<\frac{1}{2}-\varepsilon\\
\frac{\beta}{4}-\beta\left(\frac{1}{4}-\varepsilon\right)w & \mbox{ for }\frac{1}{2}-\varepsilon\le x<\frac{1}{2}+\varepsilon\\
\beta\left(x-\frac{1}{2}\right) & \mbox{ for }\frac{1}{2}+\varepsilon\le x\le1
\end{cases}\label{eq:zig-zag map}
\end{equation}
where $w$ is just a handy notation for a downward sloping line:
\[
w=\frac{2x-1}{2\varepsilon}
\]
Observe that $w=1$ when $x=\frac{1}{2}-\varepsilon$ and that $w=-1$
when $x=\frac{1}{2}+\varepsilon$ so that $w$ just smoothly interpolates
between +1 and -1 over the middle interval. The additional factors
of $\frac{\beta}{4}-\beta\left(\frac{1}{4}-\varepsilon\right)w$ just
serves to insert the downward slope smack into the middle, so that
the endpoints join up. The results is the zig-zag map, illustrated
in the figure below
\begin{center}
\includegraphics[width=0.9\columnwidth]{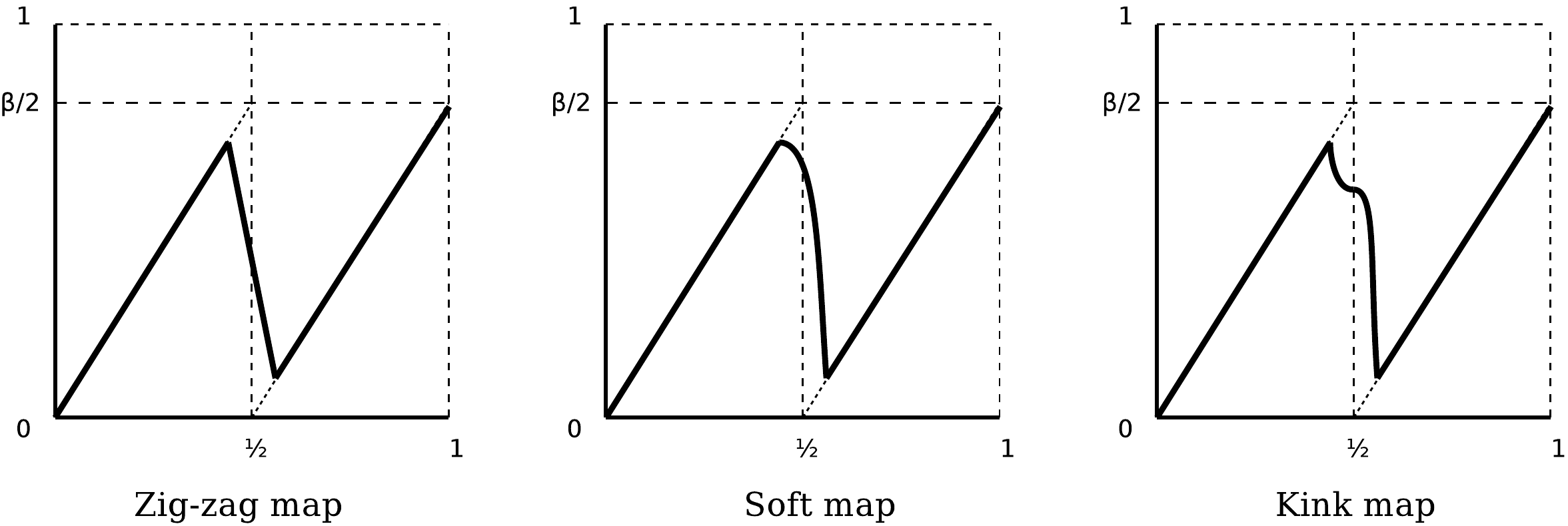}
\par\end{center}

\begin{figure}
\caption{Z-shaped Map\label{fig:Islands}}

\begin{centering}
\includegraphics[width=0.32\columnwidth]{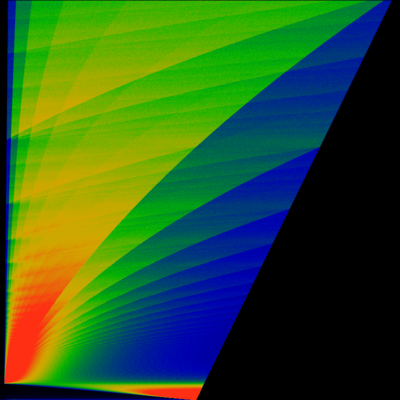}\,\includegraphics[width=0.32\columnwidth]{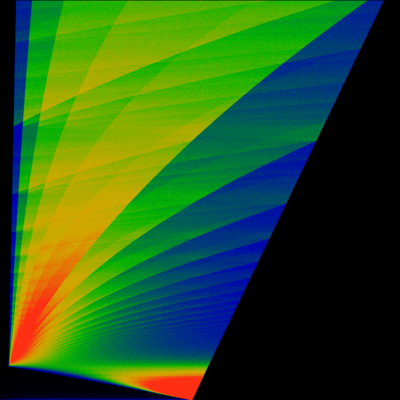}\,\includegraphics[width=0.32\columnwidth]{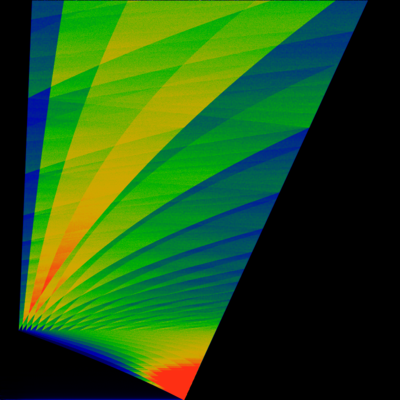}\vspace{0.003\columnwidth}
\par\end{centering}
\begin{centering}
\includegraphics[width=0.32\columnwidth]{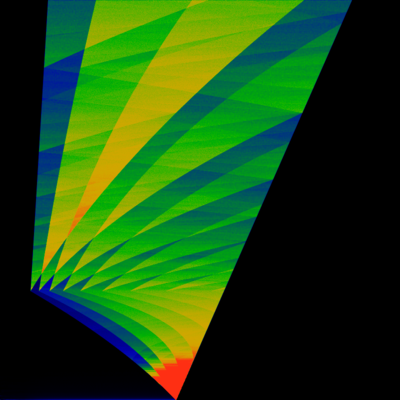}\,\includegraphics[width=0.32\columnwidth]{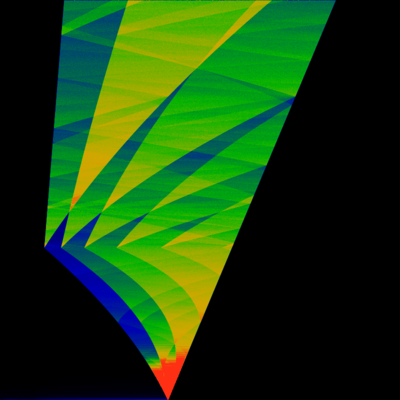}\,\includegraphics[width=0.32\columnwidth]{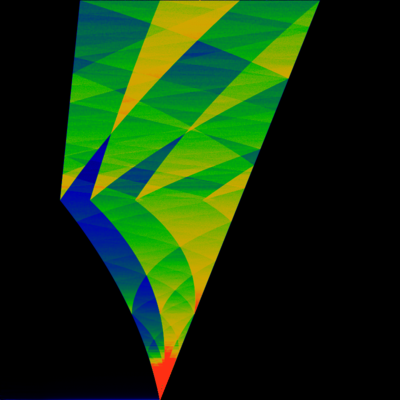}\vspace{0.003\columnwidth}
\par\end{centering}
\begin{centering}
\includegraphics[width=0.32\columnwidth]{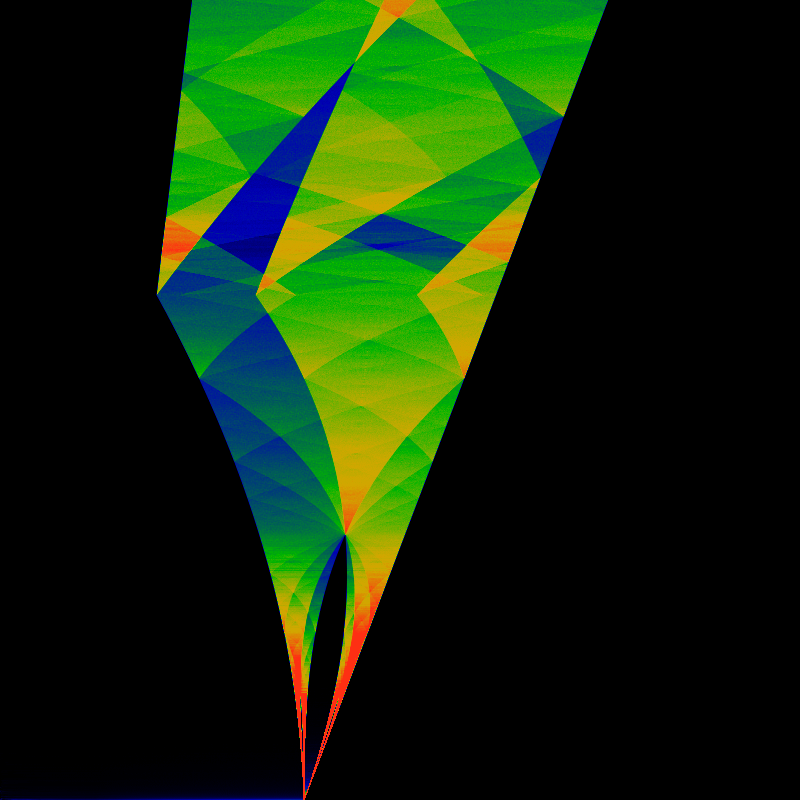}\,\includegraphics[width=0.32\columnwidth]{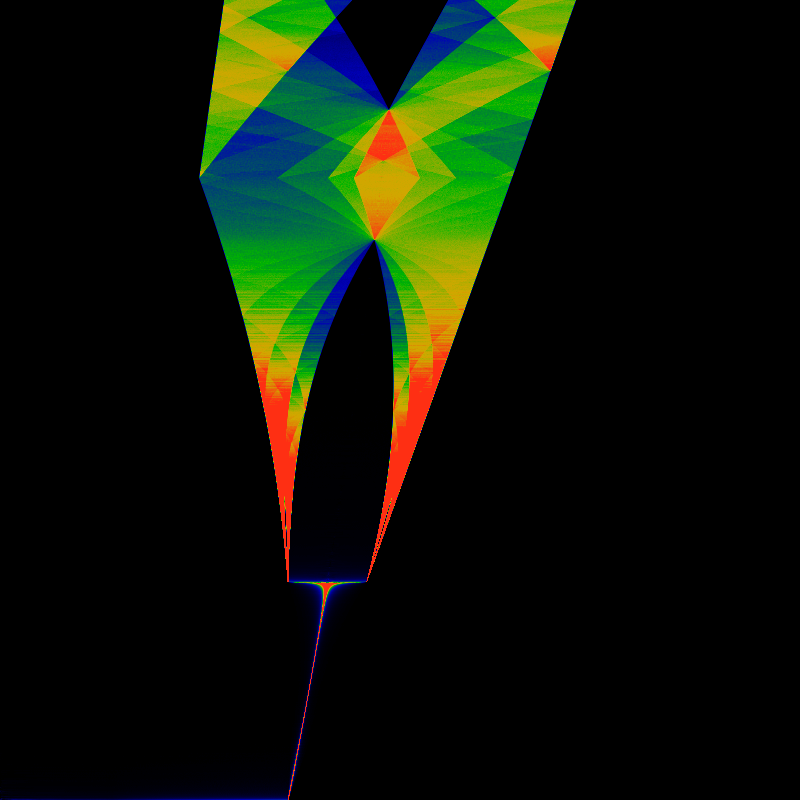}\,\includegraphics[width=0.32\columnwidth]{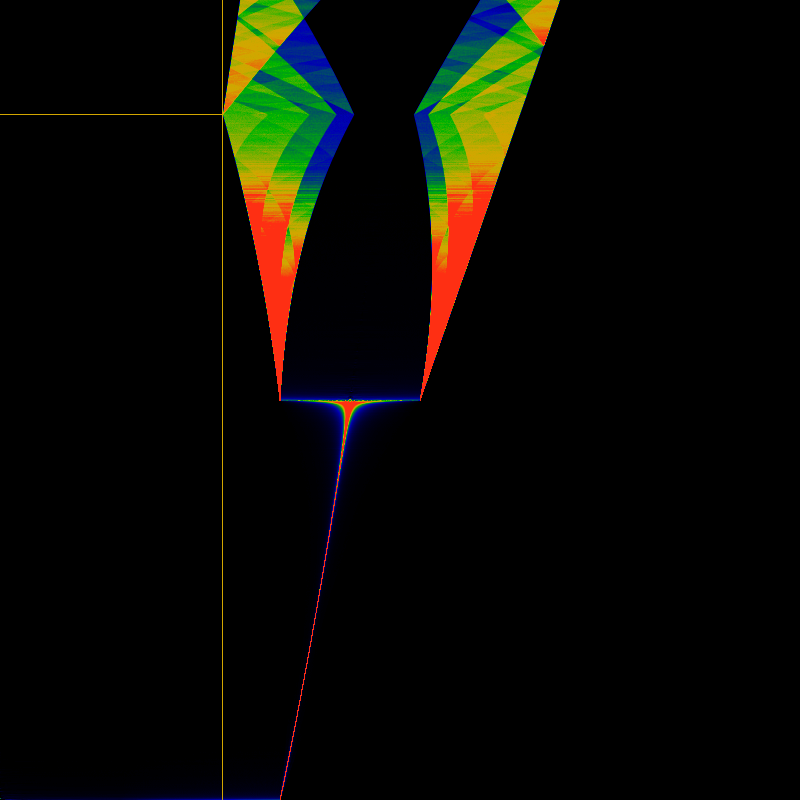}
\par\end{centering}
This illustrates a sequence of iterated maps, obtained from eqn \ref{eq:zig-zag map}.
Shown are $\varepsilon=0.01,$ 0.02, 0.04 in the first row, 0.06,
0.08, 0.10 in the second row and 0.12, 0.14, 0.15 in the third row.
The image for $\varepsilon=0$ is, of course, figure \ref{fig:Undershift-Bifurcation-Diagram}.
The parameter $\beta$ runs from 1 at the bottom to 2 at the top.
Thus, a horizontal slice through the image depicts the invariant measure
of the iterated map, black for where the measure is zero, and red
where the measure is largest. The sharp corner at the lower-left is
located $\beta=\left(1+2\varepsilon\right)/\left(1-2\varepsilon\right)$
and $x=\varepsilon\left(1+2\varepsilon\right)/\left(1-2\varepsilon\right)$.
A yellow horizontal and vertical line in the last image indicate the
location of this corner. 

\rule[0.5ex]{1\columnwidth}{1pt}
\end{figure}

In the limit of $\epsilon\to0$, one regains the earlier beta shift:
$\lim_{\varepsilon\to0}T_{\beta,\varepsilon}=T_{\beta}$, as the slope
of the middle bit becomes infinite. The middle segment is a straight
line; it introduces another folding segment into the map. This segment
introduces a critical point only when $\varepsilon$ is sufficiently
large, and $\beta$ is sufficiently small, so that its slope is less
than 45 degrees (is greater than -1). When this occurs, a fixed point
appears at $x=1/2$. A sequence of images for finite $\varepsilon$
are shown in figure \ref{fig:Islands}.

The appearance of islands of stability in the Feigenbaum attractor
is due to the presence of a fixed point at any parameter value. In
order to ``surgically add'' islands of stability to the beta transform,
the middle segment interpolation must also have a critical point at
``any'' value of $\varepsilon$. To achieve this, consider the curve

\begin{equation}
D_{\beta,\varepsilon}(x)=\begin{cases}
\beta x & \mbox{ for }0\le x<\frac{1}{2}-\varepsilon\\
\frac{\beta}{4}-\beta\left(\frac{1}{4}-\varepsilon\right)g\left(w\right) & \mbox{ for }\frac{1}{2}-\varepsilon\le x<\frac{1}{2}+\varepsilon\\
\beta\left(x-\frac{1}{2}\right) & \mbox{ for }\frac{1}{2}+\varepsilon\le x\le1
\end{cases}\label{eq:soft map}
\end{equation}
where the straight line has been replaced by a soft shoulder 
\[
g(w)=1-2\cos\frac{\pi}{4}\left(1+w\right)
\]
and $w$ is the same as before. This is scaled so that its a drop-in
replacement for the straight line: $g\left(\frac{1}{2}-\varepsilon\right)=1$
and $g\left(\frac{1}{2}+\varepsilon\right)=-1$. A cosine was used
to create this soft shoulder, but a parabola would have done just
as well. It is illustrated above, with the label ``soft map''.

\begin{figure}
\begin{centering}
\caption{Critical-Z map\label{fig:Critical-Z-map}}
\includegraphics[width=0.32\columnwidth]{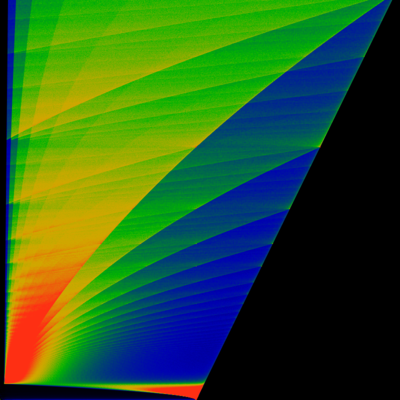}\,\includegraphics[width=0.32\columnwidth]{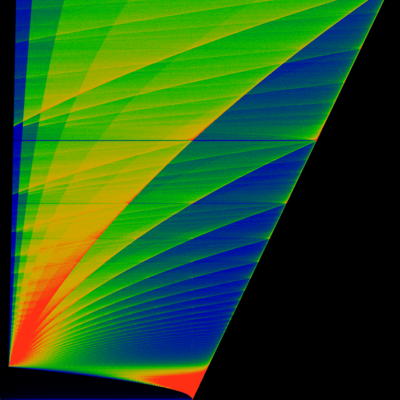}\,\includegraphics[width=0.32\columnwidth]{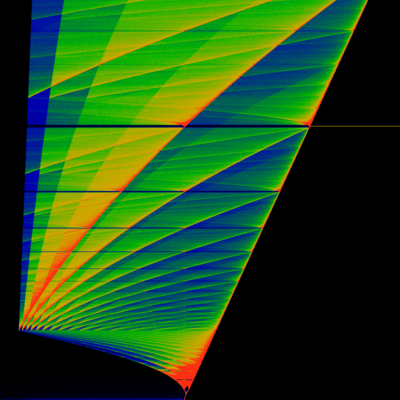}\vspace{0.003\columnwidth}
\par\end{centering}
\begin{centering}
\includegraphics[width=0.32\columnwidth]{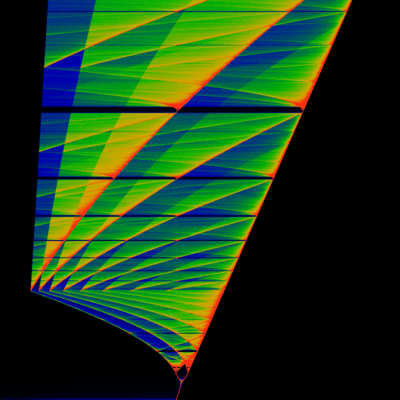}\,\includegraphics[width=0.32\columnwidth]{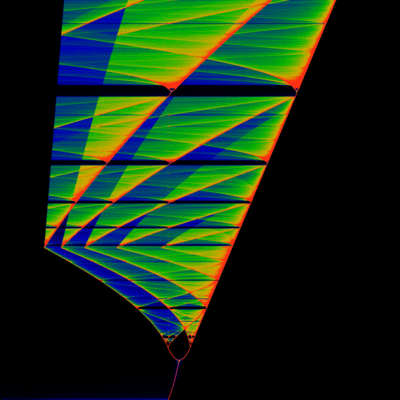}\,\includegraphics[width=0.32\columnwidth]{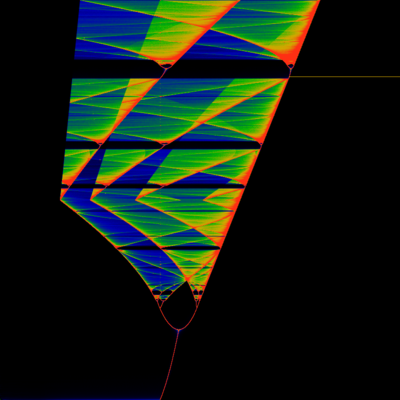}\vspace{0.003\columnwidth}
\par\end{centering}
\begin{centering}
\includegraphics[width=0.32\columnwidth]{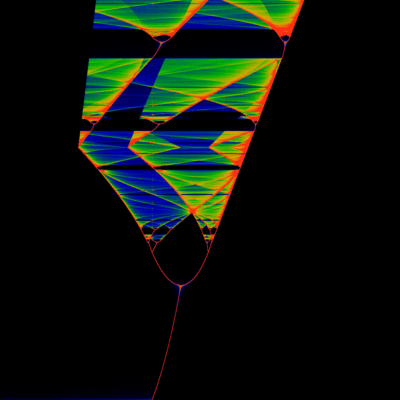}\,\includegraphics[width=0.32\columnwidth]{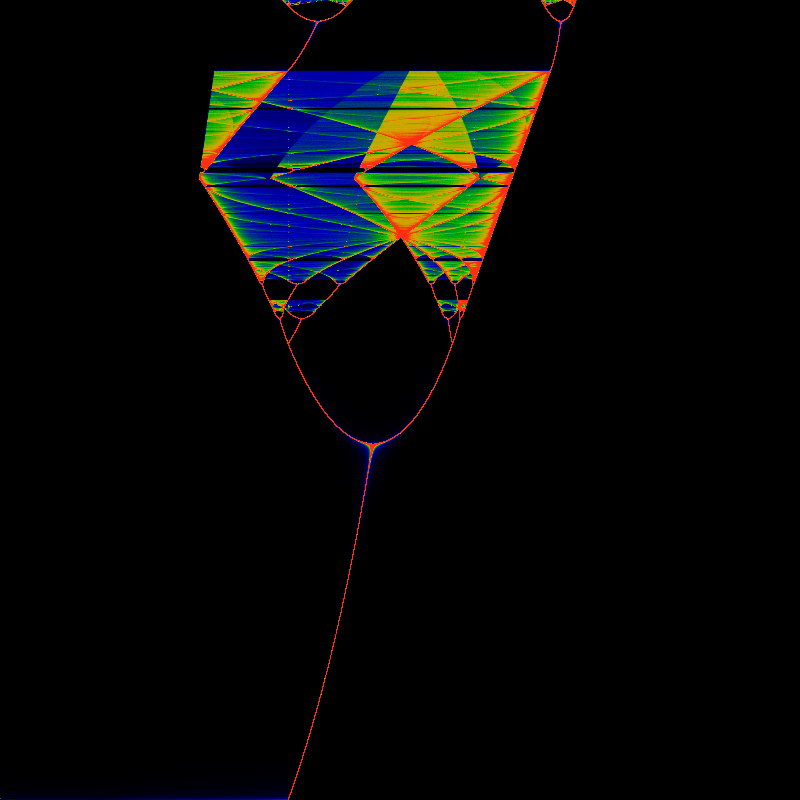}\,\includegraphics[width=0.32\columnwidth]{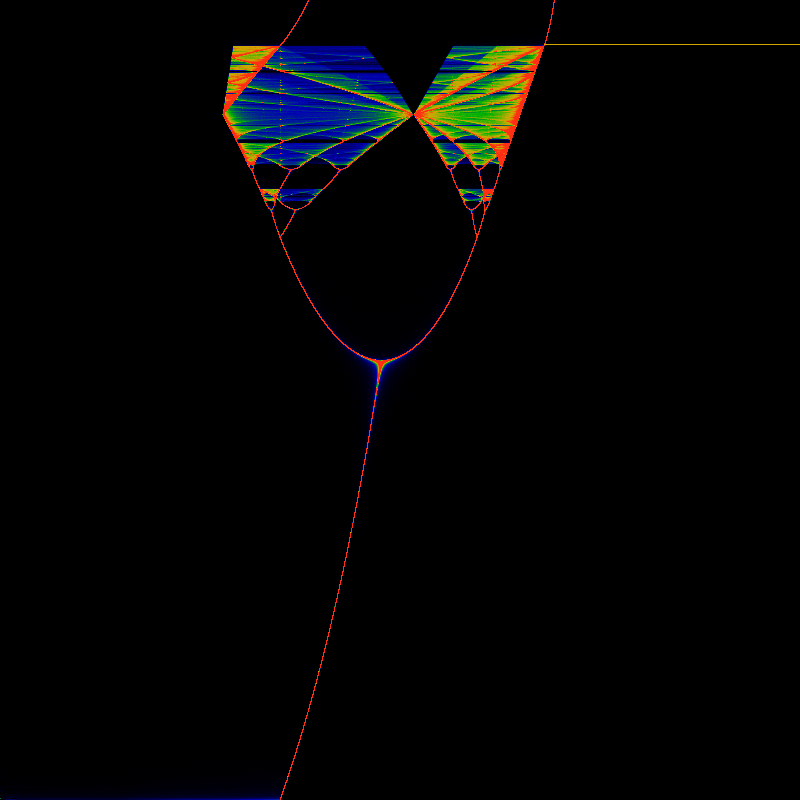}
\par\end{centering}
This illustrates a sequence of iterated maps, obtained from eqn \ref{eq:soft map}.
The sequence of depicted $\varepsilon$ values are the same as in
figure \ref{fig:Islands}. The top row shows $\varepsilon=0.01,$
0.02, 0.04, with 0.06, 0.08, 0.10 in the second row and 0.12, 0.14,
0.15 in the bottom row. The image for $\varepsilon=0$ is, of course,
figure \ref{fig:Undershift-Bifurcation-Diagram}. The parameter $\beta$
runs from 1 at the bottom to 2 at the top. Working from bottom to
top, one can see islands of stability forming in the $\varepsilon=0.02$
and 0.04 images. The largest island, one third from the top, corresponds
to $\beta=\varphi=1.618\cdots$ the golden ratio. Moving downwards,
the other prominent islands correspond to the ``trouble spots''
101, 1001 and 10001, which are the Narayana's Cows number, an unnamed
number, and the Silver Ratio, at $\beta=1.4655\cdots$ and so on.
Moving upwards, one can see a faint island at the tribonacci number.
Due to the general asymmetry of the map, these islands quickly shift
away from these limiting values. For example, the primary island appears
to start near $\beta=\delta+\left(2-\delta\right)\left(\varphi-1\right)$,
where $\delta=\left(1+2\varepsilon\right)/\left(1-2\varepsilon\right)$.
This location is indicated by a horizontal yellow line in the images
in the right column. The other islands shift away in a more complicated
fashion.

\rule[0.5ex]{1\columnwidth}{1pt}
\end{figure}

This map also interpolates between the left and right arms of the
beta transform, forming a single continues curve. The curve is smooth
and rounded near $\frac{1}{2}-\varepsilon\apprle x$, having a slope
of zero as $x$ approaches $\frac{1}{2}-\varepsilon$ from above.
This introduces a critical point near $\frac{1}{2}-\varepsilon$.
Notice that there is a hard corner at $\frac{1}{2}+\varepsilon$.
The interpolation is NOT an S-curve! A sequence of images for finite
$\varepsilon$ are shown in figure \ref{fig:Critical-Z-map}.

\begin{figure}
\begin{centering}
\caption{Interpolating Sine Map\label{fig:Interpolating-Sine-Map}}
\includegraphics[width=0.32\columnwidth]{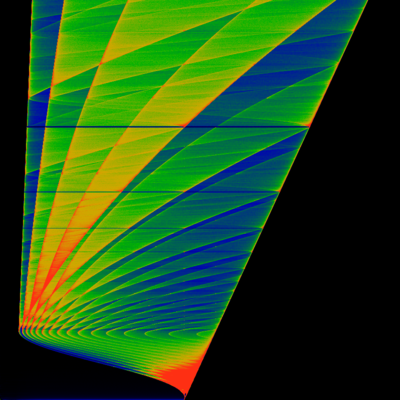}\,\includegraphics[width=0.32\columnwidth]{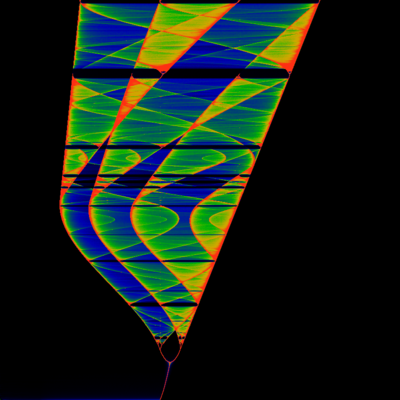}\,\includegraphics[width=0.32\columnwidth]{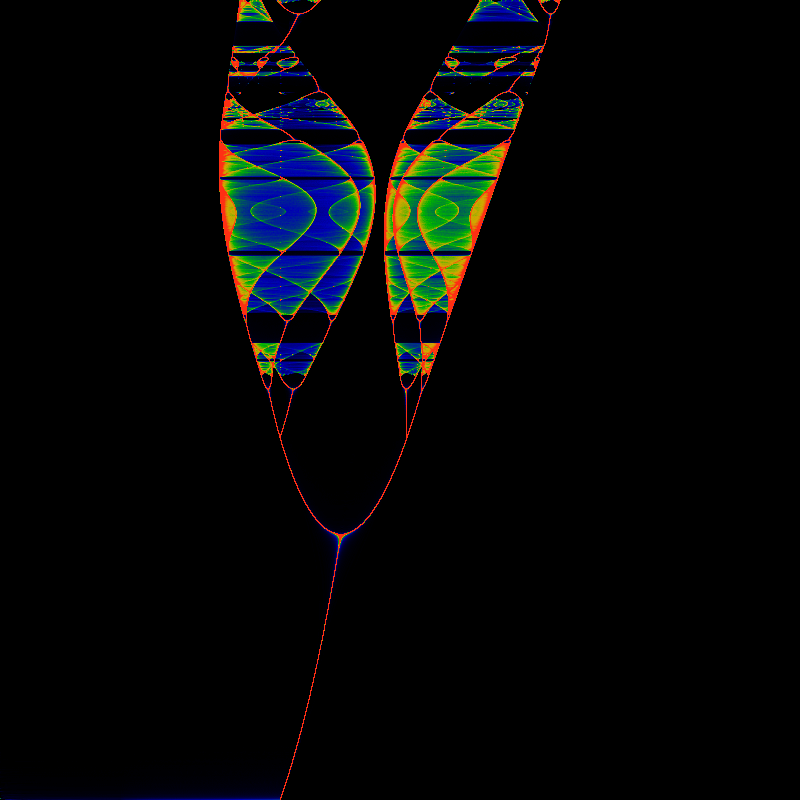}\vspace{0.003\columnwidth}
\par\end{centering}
\begin{centering}
\includegraphics[width=0.32\columnwidth]{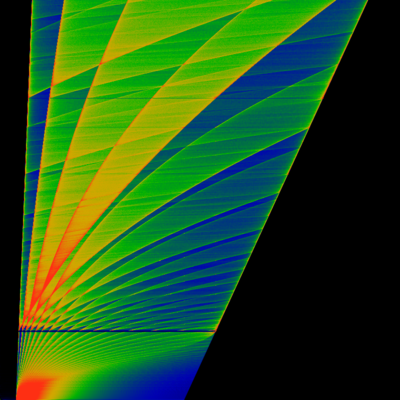}\,\includegraphics[width=0.32\columnwidth]{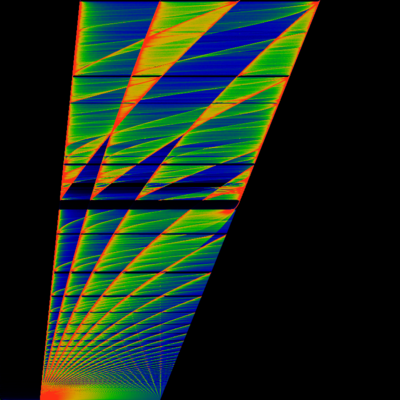}\,\includegraphics[width=0.32\columnwidth]{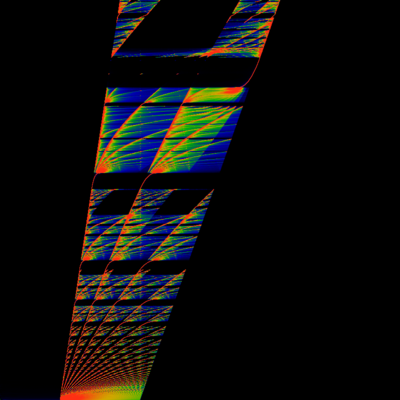}\vspace{0.003\columnwidth}
\par\end{centering}
\begin{centering}
\includegraphics[width=0.32\columnwidth]{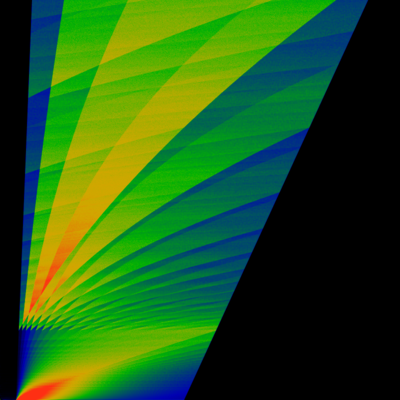}\,\includegraphics[width=0.32\columnwidth]{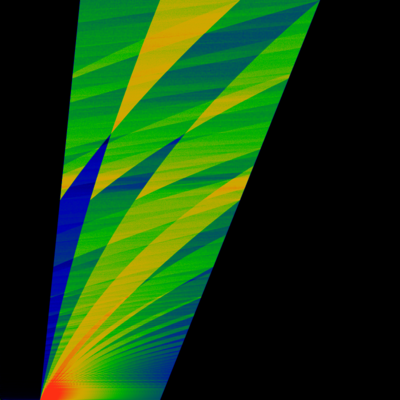}\,\includegraphics[width=0.32\columnwidth]{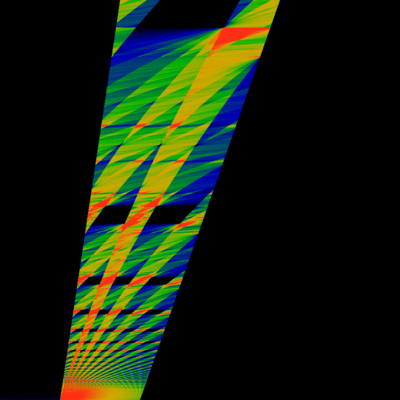}
\par\end{centering}
This illustrates a sequence of iterated maps, obtained from eqn \ref{eq:sine map}.
The sequence in the upper row shows $\varepsilon=0.04,$ 0.10 and
0.15; with $\sigma=+1$. The upper row is much like the sequence shown
in figure \ref{fig:Critical-Z-map}, except that its made sinuous,
thanks to symmetrical S-shape. The middle row shows the same $\varepsilon$
values, but for $\sigma=-1$. The bottom row shows eqn \ref{eq:hard map}
with $p=1$ and $\sigma=-1$; thus, because $p=1$ gives a straight-line
segment in the middle, this bottom row is directly comparable to the
zig-zag map. It should make clear that the islands appear in the middle
row due to critical points in the S-curve, and not due to the tripartite
map. The lower right diagram exhibits islands, but only because the
middle segment has a slope of less than 45 degrees, resulting in a
critical point at the middle of the map. As usual, the parameter $\beta$
runs from 1 at the bottom to 2 at the top. 

\rule[0.5ex]{1\columnwidth}{1pt}
\end{figure}

Two more variant maps can be considered. Both replace the center piece
with symmetrical sinuous S-shaped curves, but in different ways. Consider
\begin{equation}
S_{\beta,\varepsilon,\sigma}(x)=\begin{cases}
\beta x & \mbox{ for }0\le x<\frac{1}{2}-\varepsilon\\
\frac{\beta}{4}-\sigma\beta\left(\frac{1}{4}-\varepsilon\right)\sin\frac{\pi}{2}w & \mbox{ for }\frac{1}{2}-\varepsilon\le x<\frac{1}{2}+\varepsilon\\
\beta\left(x-\frac{1}{2}\right) & \mbox{ for }\frac{1}{2}+\varepsilon\le x\le1
\end{cases}\label{eq:sine map}
\end{equation}
and
\begin{equation}
H_{\beta,\varepsilon,p,\sigma}(x)=\begin{cases}
\beta x & \mbox{ for }0\le x<\frac{1}{2}-\varepsilon\\
\frac{\beta}{4}-\sigma\beta\left(\frac{1}{4}-\varepsilon\right)\mbox{sgn}\left(x-\frac{1}{2}\right)\left|w\right|^{p} & \mbox{ for }\frac{1}{2}-\varepsilon\le x<\frac{1}{2}+\varepsilon\\
\beta\left(x-\frac{1}{2}\right) & \mbox{ for }\frac{1}{2}+\varepsilon\le x\le1
\end{cases}\label{eq:hard map}
\end{equation}
The $S_{\beta,\varepsilon}(x)$ replaces the central segment with
a softly-rounded segment, containing two critical points: near $\frac{1}{2}-\varepsilon$
and near $\frac{1}{2}+\varepsilon$, where the curve flattens out
to a zero slope. When $\sigma=+1$, the map as a whole is continuous.
When $\sigma=-1$, the map consists of three discontinuous pieces.
Different values are explored in figure \ref{fig:Interpolating-Sine-Map}.

\begin{figure}
\begin{centering}
\caption{Interpolating Kink Map\label{fig:Interpolating-Kink-Map}}
\includegraphics[width=0.24\columnwidth]{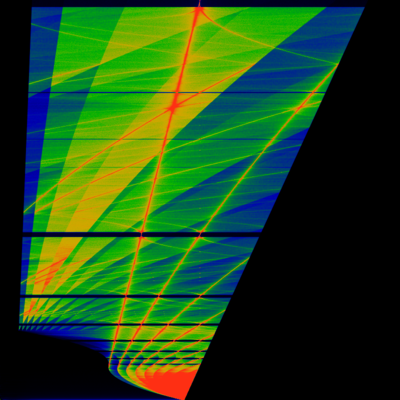}\,\includegraphics[width=0.24\columnwidth]{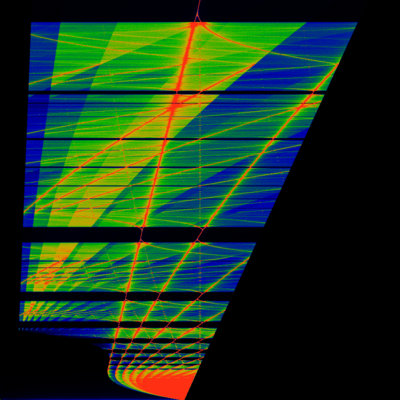}\,\includegraphics[width=0.24\columnwidth]{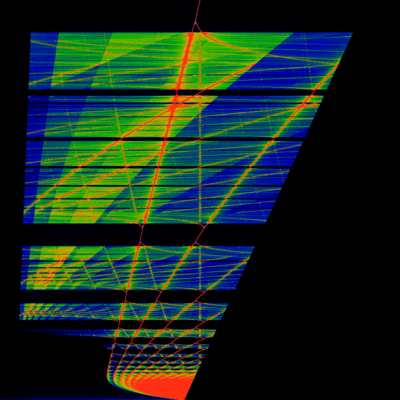}\,\includegraphics[width=0.24\columnwidth]{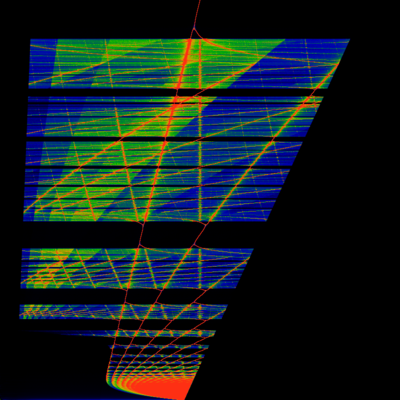}\vspace{0.003\columnwidth}
\par\end{centering}
\begin{centering}
\includegraphics[width=0.24\columnwidth]{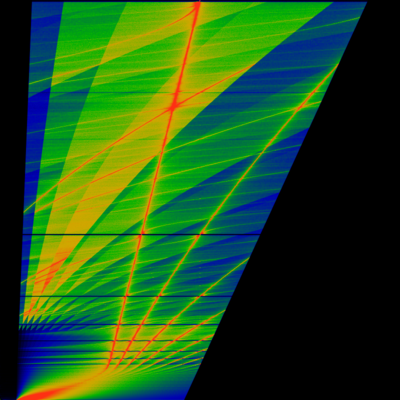}\,\includegraphics[width=0.24\columnwidth]{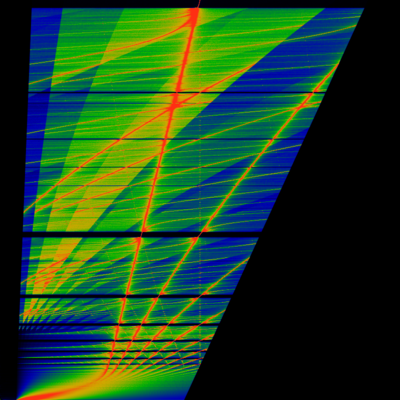}\,\includegraphics[width=0.24\columnwidth]{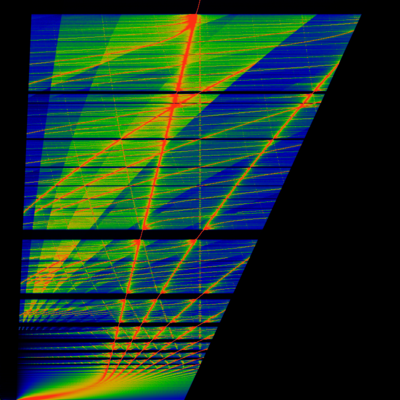}\,\includegraphics[width=0.24\columnwidth]{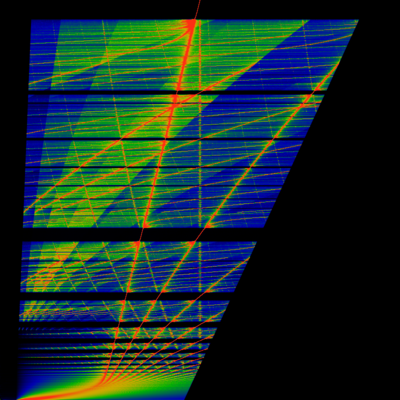}
\par\end{centering}
This illustrates a sequence of iterated maps, obtained from eqn \ref{eq:hard map}.
All eight images are held at $\varepsilon=0.04$. The top row has
$\sigma=+1$ (and thus the map is continuous) while the bottom row
has $\sigma=-1$ (and thus the map has three disconnected branches.
Left to right depicts the values $p=2,3,4,5$. As usual, the parameter
$\beta$ runs from 1 at the bottom to 2 at the top. In all cases,
islands appear, and numerous common features are evident. Perhaps
most interesting is that the islands do NOT contain period-doubling
sequences. The primary sequence of islands, starting from the central
largest, proceeding downwards, are located the inverse powers of two,
\emph{viz} at $\beta=\sqrt[k]{2}$. Why are the islands located at
inverse powers of two, instead or, for example, the golden means?
The short answer: it depends on the location of the kink in the map,
as explored in the main text. 

\rule[0.5ex]{1\columnwidth}{1pt}
\end{figure}

The $H_{\beta,\varepsilon,p,\sigma}(x)$ replaces the central segment
with a segment that has a kink in the middle, when $p>1$. Note that
$H_{\beta,\varepsilon,1,1}(x)=T_{\beta,\varepsilon}(x)$. Here, $\mbox{sgn}x$
is the sign of $x$. The general shape of $H_{\beta,\varepsilon,p,\sigma}(x)$
is shown above, labeled as the ``kink map''. The location of the
kink in $H$ is always centered; an off-center kink, as depicted in
the figure, is explored below. The bifurcation diagrams for $H$ are
illustrated in figure \ref{fig:Interpolating-Kink-Map}.

\begin{figure}
\caption{No Period Doubling\label{fig:No-Period-Doubling}}

\begin{centering}
\includegraphics[width=1\columnwidth]{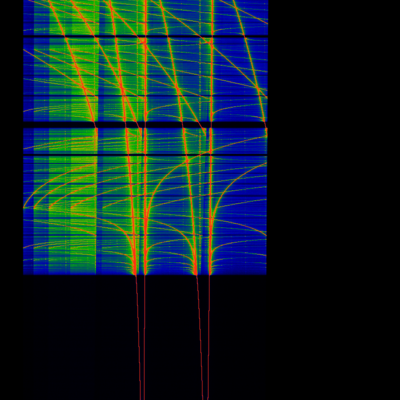}
\par\end{centering}
This figure is a zoom, confirming a lack of period doubling in the
map $H_{\beta,\varepsilon,p,\sigma}(x)$ of eqn \ref{eq:hard map}.
The explored region is $0\le x\le1$, viz no zoom in the horizontal
direction. Vertically, the image is centered on $\beta=1.45$, having
a total height of $\varDelta\beta=0.015625$. This uses the quintic
kink, so $p=5$ and $\sigma=+1$, making the the continuous variant.
The value of $\varepsilon=0.04$ makes this directly comparable to
other images.

\rule[0.5ex]{1\columnwidth}{1pt}
\end{figure}

To summarize: the ``trouble spots'' aren't ``just some periodic
orbits'' at certain values of $\beta$: they are more ``fundamental''
than that: they indicate the regions where (``phase-locked'') periodic
orbits can be made to appear. And conversely: bifurcations can only
appear here, and not elsewhere! The last sequence of images, shown
in figure \ref{fig:Interpolating-Kink-Map} indicate that the islands
of stability need NOT consist of the period-doubling sequences seen
in the Feigenbaum map. This is made explicit in figure \ref{fig:No-Period-Doubling},
which shows a zoom by a factor of thirty.

\begin{figure}
\begin{centering}
\caption{Poincaré recurrence\label{fig:Poincar=0000E9-recurrence}}
\includegraphics[width=0.49\columnwidth]{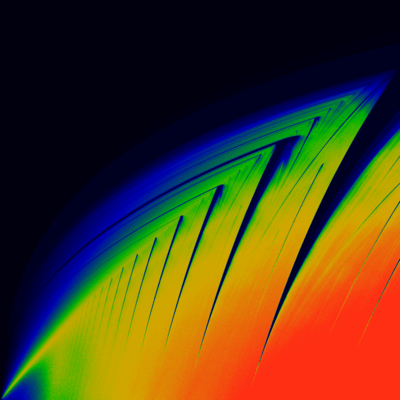}\,\includegraphics[width=0.49\columnwidth]{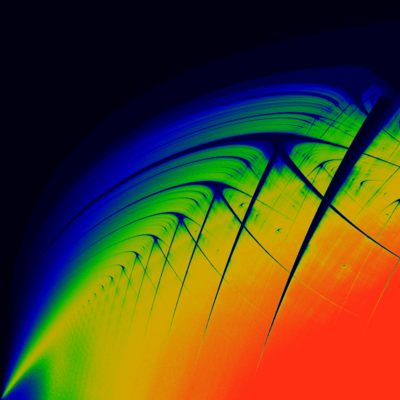}
\par\end{centering}
The above visualize the Poincaré recurrence times for the map $D_{\beta,\varepsilon}(x)$
of eqn \ref{eq:soft map} on the left, and the map $S_{\beta,\varepsilon,1}(x)$
of eqn \ref{eq:sine map} on the right. In both cases, the parameter
$\beta$ runs from 1 to 2, left to right. The parameter $\varepsilon$
runs from 0 to 0.2, bottom to top. The Poincaré recurrence time is
obtained by iterating on the maps, and then counting how many iterations
it takes to get near an earlier point. The color coding is such that
yellow/red indicates large recurrence times; green is intermediate
time, blue a short time, and black corresponds to $n$ less than 3
or 4 or so. The vertical black spikes are the Arnold tongues; they
correspond to parameter regions which lie in an island of stability.
That is, the recurrence time is low, precisely because the the point
$x$ is bouncing between a discrete set of values. The yellow/red
regions correspond to chaos, where the iterate $x$ is bouncing between
all possible values. The largest right-most spike is located at $\beta=\varphi=1.618\cdots$,
with the sequence of spikes to the left located at the other primary
golden means (\emph{viz}, $1.3803\cdots$ and the silver mean$1.3247\cdots$
and so on). As noted earlier, the general curve of that spike appears
to follow $\beta=\delta+\left(2-\delta\right)\left(\varphi-1\right)$,
where $\delta=\left(1+2\varepsilon\right)/\left(1-2\varepsilon\right)$.
The dramatic swallow-tail shapes in the right-hand image are identical
to those that appear in the \href{https://en.wikipedia.org/wiki/Arnold_tongue}{classic iterated circle map}.

\rule[0.5ex]{1\columnwidth}{1pt}
\end{figure}

\begin{figure}
\caption{Arnold Tongues\label{fig:Arnold-Tongues}}

\begin{centering}
\includegraphics[width=0.49\columnwidth]{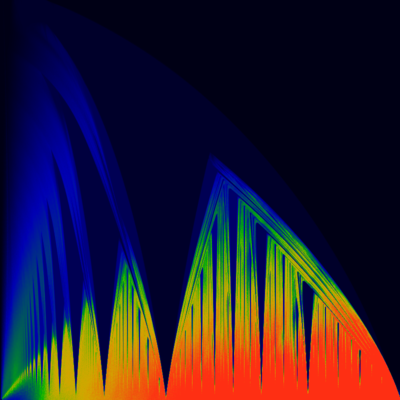}\,\includegraphics[width=0.49\columnwidth]{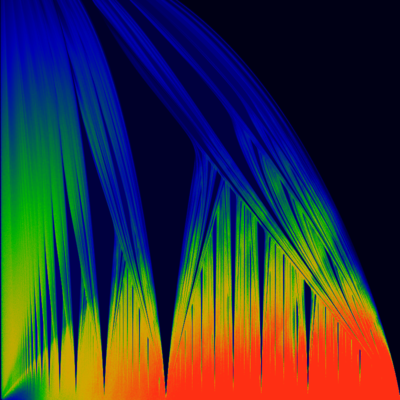}
\par\end{centering}
The above visualize the Poincaré recurrence times for the map $H_{\beta,\varepsilon,p,\sigma}(x)$
of eqn \ref{eq:hard map}. The parameter $\beta$ runs from 1 to 2,
left to right. The parameter $\varepsilon$ runs from 0 to 0.2, bottom
to top. The power $p$ is held fixed at $p=5$. The left image shows
$\sigma=$-1; the right shows $\sigma=+1$. The Poincaré recurrence
time is obtained by iterating on $H_{\beta,\varepsilon,p,\sigma}(x)$
and counting how many iterations it takes until $\left|x-H_{\beta,\varepsilon,p,\sigma}^{n}(x)\right|<0.009$.
The shapes depicted are not sensitive to the recurrence delta 0.009;
this value is chosen primarily to make the colors prettier. The color
coding is such that yellow/red indicates large recurrence times $n$;
green is intermediate time, blue a short time, and black corresponds
to $n$ less than 3 or 4 or so. The vertical blue spikes are the Arnold
tongues; they correspond to parameter regions which lie in an island
of stability. That is, the recurrence time is low, precisely because
the the point $x$ is bouncing between a discrete set of values. The
yellow/red regions correspond to chaos, where the iterate $x$ is
bouncing between all possible values. The central spike is located
at $\beta=\sqrt{2}$ with the sequence of spikes to the left located
at $\sqrt[k]{2}$ for increasing $k$. In that sense, the large black
region dominating the right side of the figures corresponds to $\beta=2$.
These correspond to the black bands in figure \ref{fig:Interpolating-Kink-Map}.

\rule[0.5ex]{1\columnwidth}{1pt}
\end{figure}

Another interesting visualization is a Poincaré recurrence plot. The
islands of stability should manifest as \href{https://en.wikipedia.org/wiki/Arnold_tongue}{Arnold tongues}.
These are shown in figures \ref{fig:Poincar=0000E9-recurrence} and
\ref{fig:Arnold-Tongues}.

To intuitively understand the location of the islands (the location
of the Arnold tongues), its easiest to examine a map with a kink in
it, whose location is adjustable. 

\[
H_{\beta,\varepsilon,\alpha,\sigma}(x)=\begin{cases}
\beta x & \mbox{ for }0\le x<\frac{1}{2}-\varepsilon\\
\frac{\beta}{4}-\sigma\beta\left(\frac{1}{4}-\varepsilon\right)h_{\alpha,p} & \mbox{ for }\frac{1}{2}-\varepsilon\le x<\frac{1}{2}+\varepsilon\\
\beta\left(x-\frac{1}{2}\right) & \mbox{ for }\frac{1}{2}+\varepsilon\le x\le1
\end{cases}
\]
with 
\[
h_{\alpha,p}\left(x\right)=\begin{cases}
\alpha+\left(1-\alpha\right)\left|w\right|^{p} & \mbox{ for }x<\frac{1}{2}\\
\alpha-\left(1+\alpha\right)\left|w\right|^{p} & \mbox{ for }\frac{1}{2}\le x
\end{cases}
\]
As before, $h_{\alpha,p}\left(x\right)$ is designed to interpolate
appropriately, so that $h_{\alpha,p}\left(\frac{1}{2}-\varepsilon\right)=1$
and $h_{\alpha,p}\left(\frac{1}{2}+\varepsilon\right)=-1$. The location
of the kink is now adjustable: $h_{\alpha,p}\left(\frac{1}{2}\right)=\alpha$.
Iterating on this map results in figures that are generically similar
to those of figure \ref{fig:Interpolating-Kink-Map}, except that
this time, the location of the islands is controllable by the parameter
$\alpha$. Roughly, to first order, the primary series of islands
are located at $\sqrt[k]{2/\left(1-\alpha\right)}$; as before, these
islands do not allow period-doubling to take place. 

To get islands with period doubling, one needs to recreate the ``soft
shoulder'' of eqn \ref{eq:soft map}, but at a variable location.

Thus, the above presents a general surgical technique for controlling
both the general form of the chaotic regions, the location of the
islands of stability, and what appears within the islands.

Conjectures are fun! The above arguments should be sufficient to fully
demonstrate that the circle map, which is well-known to exhibit phase
locking regions called Arnold tongues, is topologically conjugate
to the fattened beta shift $T_{\beta,\varepsilon}$. Or something
like that. In a certain sense, this can be argued to be a ``complete''
solution, via topological conjugacy, of the tent map, the logistic
map and the circle map. This is a worthwhile exercise to actually
perform, i.e. to give explicit expressions mapping the various regions,
as appropriate. 

Essentially, the claim is straight-forward: topologically, all chaotic
parts of a map correspond to folding (as per Milnor, 1980's on kneading
maps), into which one may surgically insert regions that have cycles
of finite length. The surgical insertion can occur only at the discontinuities
of the kneading map. It almost sounds trivial, expressed this way;
but the algebraic articulation of the idea would be worthwhile.

\pagebreak{}

\section{Conclusion}

The idea of analytic combinatorics takes on a whole new meaning in
the computational age. Historically, the ability to provide an ``exact
solution'' in the form of an analytic series has been highly prized;
the ultimate achievement in many cases. Being able to expression a
solution in terms of the addition and multiplication of real numbers
is very comforting. Every school student eventually comes to feel
that arithmetic on the real numbers is very natural and normal. It's
more than that: Cartesian space is smooth and uniform, and all of
differential geometry and topology are founded on notions of smoothness.

The inner workings of computers expose (or hide!) a different truth.
The most efficient algorithm for computing $sine(x)$ is not to sum
the analytic series. Arbitrary precision numerical libraries open
the rift further: neither addition nor multiplication are simple or
easy. Both operations have a variety of different algorithms that
have different run-times, different amounts of memory usage. In the
effort to minimize space and time usage, some of these algorithms
have grown quite complex. The root cause of the complexity is bewildering:
it is the use of the binary digit expansion to represent a real number.
Computers use the Cantor space $\left\{ 0,1\right\} ^{\omega}$ or
at least a subset thereof, under the covers.

Different representations of the real numbers potentially offer different
algorithms and performance profiles. One could represent reals by
rationals, but then several other issues arise. One is that the rationals
are not evenly distributed across the real number line: rationals
with small denominators cluster about in a fractal fashion. This is
easily exhibited by considering continued fractions. As a result,
one promptly gets stuck in a quagmire of trying to understand what
a ``uniform distribution'' should be. Binary expansions are more
``obviously'' uniform. A more basic issue is that, if working with
rationals, one must somehow accomplish the addition or multiplication
of two integers. To accomplish this, one has to represent the integers
as sequences of bits, which only takes us back to where we started.
There is no computational oracle that automatically knows the sum
or product of integers: it has to be computed.

Compare this situation to that of iterated functions and fractals.
At first impression, these seem pathological in almost every respect:
differentiable nowhere, unbounded and nonuniform: somehow they feel
like the quintessential opposite of the analytic series, of the smoothness
of Cartesian space, of the smoothness of addition and multiplication.
The place where these two worlds come together is that both are attempts
to approach countable infinity, and both are attempts to harness the
first uncountable infinity. The real number number is an infinite
string of binary digits. The analytic series is an infinite sum. The
iterated function is recursively infinite. The historic labor of finding
``exact solutions'' to problems can perhaps be better views as the
discovery of correspondences between finite structures (``the problem
to be solved'') and infinite structures (``the solution'').

The situation here is more easily illustrated in a different domain.
The hypergeometric series was presented and studied by Gauss; then
Kummer, Pfaff and Euler observed various identities yoking together
different series. By the 1950's, thousands of relations were known,
along with some algorithms that can enumerate infinite series of relations.
The current situation is that there is no known algorithm that can
enumerate all such relations; there is no systematic way to classify
them. There is an interplay between infinite series and algorithmic
relationships between them. Stated a different way: hypergeometric
series have a class of self-similarities, and the identities relating
them are expressions of that self-similarity. What is that class of
self-similarities? For the hypergeometric series, it remains unknown.

For Cantor space, that place where we represent real numbers, the
situation is much better. The Cantor space itself has the structure
of an infinite binary tree; the tree and it's subtrees are obviously
self-similar; the class of similarities is described by the dyadic
monoid. The dyadic monoid embeds naturally into the modular group;
this in turn is a gateway to vast tracts of modern mathematics. The
recursive aspects, the shadow that the Cantor space seems to leave
behind everywhere appears to be ``explained'' by Ornstein theory.

Yet, the picture remains incomplete. The $\beta$-transform provides
a simple, silly model for multiplying two real numbers together: $\beta$
and $x$. The ``extra complication'' of taking mod 1 after multiplication
just reveals how complex multiplication really is. After all, mod
1 is just the subtraction of 1; how hard can that be? Moving in one
direction: the fastest, most efficient-possible algorithm for multiplying
two numbers is not known. Moving in another direction, the simple
iterated maps, shown in figures \ref{fig:Undershift-Bifurcation-Diagram},
\ref{fig:Tent-Map-Bifur} and \ref{fig:Logistic-Map-Bifurcation}
are obviously not only self-similar, but also are surely topologically
conjugate to one-another, and in all cases are presumably described
by the dyadic monoid; likewise the Mandelbrot set and it's exterior.
Yet the details remain obscure.

The meta-question is: what is the correct framework by which one can
best understand the interplay between symmetries, infinite series,
infinite recursion and algorithms? Until modern times, mathematical
practice has reified addition and multiplication into oracular operations
that magically obtain ``the right answer''. Modern computers have
put a lie to this: the theory of numerical methods has made clear
that addition and multiplication are necessarily algorithmic operations
performed on finite truncations of infinite series. What other algorithms
are hiding nearby, and what is their relationship to analytic series?

\pagebreak{}

\bibliographystyle{tufte}
\bibliography{solvit}

\end{document}